\documentclass[reqno,11pt]{article}
\usepackage{amsmath, latexsym, amsfonts, amssymb, amsthm}
\usepackage{graphicx, color,hyperref,dsfont,epsfig,caption,wrapfig,subfig}
\usepackage{pstricks}

\hypersetup{
    linktoc=page,
    linkcolor=red,          
    citecolor=blue,        
    filecolor=blue,      
    urlcolor=cyan,
    colorlinks=true           
}

\numberwithin{equation}{section}

\newtheorem{theorem}{Theorem}[section]
\newtheorem{lemma}[theorem]{Lemma}
\newtheorem{proposition}[theorem]{Proposition}

\newtheorem{rem}[theorem]{Remark}
\newtheorem{definition}[theorem]{Definition}

\newtheorem{conjecture}[theorem]{Conjecture}

\newtheorem{example}[theorem]{Example}

\renewcommand{\tilde}{\widetilde}          
\DeclareMathSymbol{\leqslant}{\mathalpha}{AMSa}{"36} 
\DeclareMathSymbol{\geqslant}{\mathalpha}{AMSa}{"3E} 
\DeclareMathSymbol{\eset}{\mathalpha}{AMSb}{"3F}     
\renewcommand{\leq}{\;\leqslant\;}                   
\renewcommand{\geq}{\;\geqslant\;}                   

\newcommand{\C}{\mathbb{C}}
\newcommand{\R}{\mathbb{R}}
\newcommand{\Z}{\mathbb{Z}}
\newcommand{\N}{\mathbb{N}}

\newcommand{\E}{\mathds{E}}
\newcommand{\Pb}{\mathds{P}}
\newcommand{\ind}{\mathds{1}}
\newcommand{\ps}[1]{\langle #1 \rangle}

\title{Gaussian multiplicative chaos and applications: a review}
\author{}

\begin{document}

\maketitle
\begin{center}
{R\'emi Rhodes\\
\footnotesize \noindent
Universit{\'e} Paris-Dauphine, Ceremade, F-75016 Paris, France}

{\footnotesize \noindent e-mail: \texttt{rhodes@ceremade.dauphine.fr}}

\bigskip

{Vincent Vargas \\
\footnotesize 
 
 CNRS, UMR 7534, F-75016 Paris, France \\
  Universit{\'e} Paris-Dauphine, Ceremade, F-75016 Paris, France} \\

{\footnotesize \noindent e-mail: \texttt{vargas@ceremade.dauphine.fr}}
\end{center}

\begin{abstract}
In this article, we review the theory of Gaussian multiplicative chaos initially introduced by Kahane's seminal work in 1985. Though this beautiful paper faded from memory until recently, it already contains ideas and results that are nowadays under active investigation, like the construction of the Liouville measure in $2d$-Liouville quantum gravity or thick points of the Gaussian Free Field. Also, we mention important extensions and generalizations of this theory that have emerged ever since and discuss a whole family of applications, ranging from finance, through the Kolmogorov-Obukhov model of turbulence  to $2d$-Liouville quantum gravity. This review also includes new results like the convergence of discretized Liouville measures on isoradial graphs (thus including the triangle and square lattices) towards the continuous Liouville measures (in the subcritical and critical case) or multifractal analysis of the measures in all dimensions.
\end{abstract}
\vspace{1cm}
\footnotesize


\noindent{\bf Key words or phrases:} Random measures, Gaussian processes, Multifractal processes.

\noindent{\bf MSC 2000 subject classifications: 60G57, 60G15, 60G25, 28A80}

\normalsize

\tableofcontents

\section{Introduction}
Log-normal multiplicative martingales were introduced by Mandelbrot \cite{cf:Man} in order to  build random measures describing energy dissipation and contribute explaining intermittency effects in  Kolmogorov's theory of fully developed turbulence (see \cite{cf:Castaing,cf:Sch,cf:Sto,cf:Cas,cf:Fr} and references therein). However, his model was difficult to define mathematically and  this is why he proposed  in \cite{mandelbrot} the simpler model of random multiplicative cascades whose detailed study started with Kahane's and Peyri\`ere's notes~\cite{kahane74,pey74}, improved and gathered in their joint paper~\cite{KP}.

From that moment on, multiplicative cascades have been widely used as reference or toy models in many applications as they feature beautiful stochastic  scaling relations, modeling a phenomenon that is commonly called intermittency. Let us roughly explain this point. Consider a  multiplicative cascade $M$ constructed on a dyadic tree. It is a  random measure over the interval $[0,1]$. If you look at the measure at a dyadic scale $2^{-n}$, you observe the same object as $M$ up to an independent stochastic factor:
$$M[\frac{k}{2^n},\frac{k+1}{2^n}]\stackrel{law}{=} e^{\Omega_n}M[0,1],$$
where $\Omega_n$ is a random variable independent of $M$, the law of which depends on the scale $2^{-n}$. However, multiplicative cascades are constructed on a dyadic (or p-adic) tree and therefore possess many drawbacks: they do not possess stationary fluctuations and present discrete (p-adic) scaling relations.

Gaussian multiplicative chaos, introduced by Kahane \cite{cf:Kah} in 1985, is born from the need of making rigorous Mandelbrot's initial model of energy dissipation \cite{cf:Man}, the so-called Kolmogorov-Obukhov model. It is about constructing a continuous parameter theory of suitable multifractal random measures. Kahane's efforts were followed by several authors \cite{Rnew1,cf:BaMu,Bar,Fan,Rnew8,Rnew9,cf:RoVa,cf:Sch} coming up with various generalizations at different scales. These measures have found many applications in various fields of science, especially in mathematical finance (or equivalently boundary Liouville quantum gravity), $2d$-Liouville quantum gravity and $3d$-turbulence.

In dimension $d$, a standard Gaussian multiplicative chaos is a random measure on a given domain $D$ of $\R^d$ that can be formally written,  for any Borelian set $A$ as:
\begin{equation}\label{measintrorev}
M_\gamma(A)=\int_Ae^{\gamma X(x)-\frac{\gamma^2}{2}\E[X^2(x)]}\,\sigma(dx)
\end{equation}
 where $X$ is a centered Gaussian "field" and $\sigma$ is a Radon measure on $D$. 
  In the situations of interest,   $X$ is rather badly behaved and cannot be defined as a random function: it is a random {\em distribution} (in the sense of Schwartz), like the Gaussian Free Field (GFF for short, see \cite{She07} for an overview about the GFF) for instance. In his seminal work, Kahane focused on the case where $X$ possesses a covariance kernel of the form:
\begin{equation}\label{Kintrorev}
\E[X(x)X(y)]=\ln _+\frac{1}{|x-y|}+g(x,y),
\end{equation} with $\ln _+(u)=\max(\ln u , 0)$ and $g$  a continuous bounded function over $D\times D$. Surprisingly, it turns out that this is the only situation of interest since this family of kernels can be thought of as a transition, separating the family of kernels for which \eqref{measintrorev} is trivially converging from the family of kernels  for which \eqref{measintrorev} is trivially vanishing.  The covariance kernel thus possesses a singularity and it is now clear that giving sense to  \eqref{measintrorev} is not straightforward (how do you define the exponential of a distribution?). The standard approach consists in applying a "cut-off" to the distribution $X$, that is in regularizing the field $X$ in order to get rid of the singularity of the covariance kernel and get a nicer field. The regularization usually depends on a small parameter that stands for the extent to which the field has been regularized. The measure \eqref{measintrorev} is naturally understood as the limit that you get when the regularization parameter goes to $0$. Kahane's paper \cite{cf:Kah} is about making  this sketch of construction  rigorous.  When $\sigma$ is the Lebesgue measure for instance, it turns out that it produces non trivial limiting objects when the parameter $\gamma$ is less than some critical value $\gamma_c=\sqrt{2d}$. These are the foundations of Kahane's theory. Several questions are then raised like:\\
- Does the limiting measure $M_\gamma$ depend on the chosen cut-off procedure?\\
- What are the geometrical and statistical properties of the measure $M_\gamma$?\\
- What are the regularity properties of the measure $M_\gamma$?\\
- How can we characterize the measure  $M_\gamma$?\\
- What happens at $\gamma_c$? and what about $\gamma>\gamma_c$?

In this review, we will discuss the above questions as well as possible generalizations and applications of   
 the theory of Gaussian multiplicative chaos   in the light of recent progresses. Among the main points that we will address are the convolution techniques used to produce a Gaussian multiplicative chaos. The situation may be summarized as follows: given a Gaussian distribution $X$ with covariance kernel of the type \eqref{Kintrorev}, what if we regularize $X$ by convolution with a smoothing family of functions $(\theta_\epsilon)_\epsilon$, converging towards the Dirac mass at $0$ when $\epsilon\to 0$? If we formally define: 
 $$X_\epsilon(x)=\int X(y)\theta_\epsilon(x-y)\,dy,$$ it turns out that, under weak conditions, the family of random measures:
 \begin{equation}\label{measintrorevapproch}
M_{\epsilon,\gamma}(A)=\int_Ae^{\gamma X_\epsilon(x)-\frac{\gamma^2}{2}\E[X_\epsilon^2(x)]}\,dx
\end{equation}
weakly converges in law towards the same measure as that produced by Kahane's theory. Convolution techniques were first introduced in \cite{cf:RoVa1,cf:RoVa}, where convergence in law is established. In the particular case when $X$ is a Gaussian Free Field (GFF), a particular convolution technique was also studied in \cite{cf:DuSh}, i.e. convolution of the field by a circle, where almost sure convergence along subsequences is established (see section \ref{sec:extent}).
 
We  will also discuss multifractality of Gaussian multiplicative chaos  and related scaling relations. Roughly speaking, multifractal analysis is the study of objects, like measures or functions, possessing several levels of local regularity: for instance, the local H\"older exponent may vary spatially. In particular, we will explain why the (nowadays called) thick points of the GFF are very closely related to a general theory called multifractal analysis, which at least goes back to Kahane's paper when applied to Gaussian multiplicative chaos (see subsection \ref{sec:multi}). Nowadays, there is a huge amount of literature on multifractal analysis and it is far beyond the scope of this review to cite or discuss all the related mathematical achievements (in the case of Mandelbrot's multiplicative cascades see \cite{Barral1,barralJTP,BMI} and references therein).

The applications that we will mention range from mathematical finance to fully developed turbulence through $2d$-Liouville quantum gravity or decaying Burgers turbulence. On the one hand, some of them are rather well established so that we will just recall the basic framework and give references. On the other hand, we feel important to devote a considerable part of the paper to $2d$-Liouville quantum gravity. Roughly, it can be seen as an attempt to construct a canonical "Riemannian" random metric on the sphere. Since Polyakov's original work \cite{Pol}, physicists have understood that such a metric takes on the form of the exponential of a GFF \cite{cf:Da,Pol}. Interestingly, they have also understood that such metrics can be discretized by randomly triangulated surfaces (see \cite{amb} for instance), a recent field of mathematical research that culminated with the construction of the so-called Brownian map (see \cite{LeGallS,legall,LGM,miermont}) in the special case of pure gravity. While constructing properly a metric is presently out of reach, Gaussian multiplicative chaos theory gives straightforwardly a way of defining the associated volume form, called the Liouville measure. Non specialists of Gaussian multiplicative chaos will find in subsection \ref{LQG} the different constructions suggested in the literature, i.e. white noise decomposition/circle average/$H^1$-expansion of the GFF, as well as a proof that they all produce the same measure in law. Let us mention here the remarkable works of Duplantier-Sheffield \cite{cf:DuSh} and Sheffield \cite{QZ} where the authors state an impressive series of conjectures, which can be seen as a starting point to understand  the physicist picture (see also the nice review for mathematicians \cite{garban}). In particular, this ambitious program could give a rigorous geometrical framework to the Knizhnik-Polyakov-Zamolodchikov formula  (KPZ for short, see \cite{cf:KPZ}). For instance, such a KPZ formula has already been used by physicists \cite{DupKwo} to predict the exact values of the Brownian intersection exponents and has been rigorously checked to hold in some special cases (see \cite{BenjCur}). We will explain the geometrical KPZ formulae rigorously proved in \cite{cf:DuSh,Rnew10} (see also \cite{Benj,Rnew4}). 

 Furthermore, as pointed out in \cite{Rnew10}, the theory of Gaussian multiplicative chaos allows us to deal with much more general situations concerning constructions of measures or the KPZ formula. We stress here that the theory of Gaussian multiplicative chaos in \cite{cf:Kah} and the KPZ relation proved in \cite{Rnew10} (or \cite{Rnew4}) is valid in any dimension when applied to log-correlated Gaussian fields. In particular, boundary Liouville measures are discussed since they are nothing but Gaussian multiplicative chaos along $1d$ Riemannian manifolds. Also, in dimension $d$ one can consider situations where the field $X$ has correlations given by the kernel $(m^2-\Delta)^{-d/2}$  (with possibly $m=0$) since such correlations are logarithmic (see \cite{DSRV3}). We do not  detail this situation here, first because we cannot explain in great details all the situations where the theory of Gaussian multiplicative chaos applies and second, because it could be instructive  for the reader to check that the framework drawn in \cite{cf:Kah} for the construction of measures or in \cite{Rnew10,Rnew4} for the KPZ formula applies (the reader may also consult \cite{chen} on this topic). Finally, as a new result, we explain how to combine the results in \cite{chelkak,cf:RoVa} to prove that the discrete Liouville measures on isoradial graphs converge towards the Liouville measure as the mesh of the graph converges to $0$. 

The last part of this review (section \ref{sec:gener}) will be devoted to possible generalizations of Kahane's theory. We will discuss how to renormalize the vanishing measure \eqref{measintrorev} for $\gamma^2\geq 2d$ (when $\sigma$ is the Lebesgue measure). This yields new qualitative behaviours of the limiting measure that may be classified in two  categories:  {\it Critical Gaussian Multiplicative Chaos} or {\it Atomic Gaussian Multiplicative Chaos}. We discuss the associated relations with duality in $2d$-Liouville quantum gravity, the frozen phase of logarithmically correlated Gaussian potentials or the maximum of log-correlated Gaussian fields (including in particular the maximum of the GFF). Other possible generalizations are mentioned, like taking complex-valued fields $X$  or matrix-valued fields $X$ in \eqref{measintrorev}.

 Finally, we mention that a preceeding review on multiplicative chaos has already appeared \cite{reviewfan}. The reader may find in \cite{reviewfan} some further fields of applications including Dvoretzy covering, percolation on trees, random cascades and Riesz products that we do not review here for the sake of non-overlapping. Concerning multiplicative cascades, the reader may consult the recent review \cite{BFP}.

\subsection{A word on Quantum field theory and the Hoeg-Krohn model}
To our knowledge, the first mathematical occurence of measures of the form \eqref{measintrorev} appeared in Hoeg-Krohn's work \cite{Hoeg}. In the context of Quantum field theory,  Hoeg-Krohn focused on the case where $X$ is the two dimensional massive free field  and $\sigma$ the Lebesgue measure (to be precise, their point of view is not exactly that of random measures). He  showed that the measure $M_\gamma$ is non trivial for $\gamma^2<2$, thus working below the $L^2$-threshold; under the assumption $\gamma^2<2$, one can perform $L^2$-computations which considerably simplifies the study of Gaussian multiplicative chaos measures (see subsection \ref{seminal}). This work led to other works \cite{Albeverio1,Albeverio2} in dimension $2$ which generalized some of the initial results of \cite{Hoeg}. Nonetheless, it seems that none of these works focused on building a general theory applicable to a wide range of random measures in all dimensions.

\subsection{Notations}
The truncated logarithm $\ln_+$ is the function $\ln_+(x)=\max(0,\ln x)$.  The relation $f\asymp g$ means that there exists a positive constant $c>0$ such that $c^{-1}f(x)\leq g(x)\leq c f(x)$ for all $x$ under consideration.

We denote by $(D,\rho)$ a metric space $D$ equipped with its metric $\rho$. This metric space is endowed with its Borelian sigma algebra $\mathcal{B}(D)$. 

Given a domain $D$ of $\R^d$, we denote by $ H^1_0(D)$ the classical Sobolev space defined as the Hilbert space closure  with respect to the Dirichlet inner product of the set of smooth compactly supported functions on $D$.
Finally, $\gamma$ will denote the intermittency parameter. We will suppose that $\gamma \geq 0$.  

\subsection*{Acknowledgements}
We would like to thank M. Bauer, F. David, B. Eynard, J.F. Le Gall for interesting discussions on these topics. Special thanks are addressed to all those who read and suggested improvements of a prior draft of this manuscript: J. Barral, N. Curien and C. Garban.

\section{State of the art since Kahane}\label{state}
 
\subsection{The seminal work of Kahane in 1985}\label{seminal}
The theory of multiplicative chaos was first defined rigorously by Kahane in 1985 in the article \cite{cf:Kah} to which the reader is referred for further details or definitions. More specifically, 
Kahane built a theory relying on the notion of $\sigma$-positive type kernel. Consider a locally compact metric space $(D,\rho)$. A function $K:D \times D \rightarrow \R_{+}\cup \lbrace \infty \rbrace$ is of $\sigma$-positive type if there exists a sequence $(K_{k})_k$ of continuous nonnegative and positive definite kernels $K_k:D \times D \rightarrow \R_{+}$ such that:
\begin{equation}\label{eq:defsigma}
\forall x,y\in D,\quad K(x,y)=\sum_{k \geq 1}K_{k}(x,y).
\end{equation}
    
It is worth pointing out here that Kahane's theory uses nonnegativeness of the kernels $(K_k)_k$ for the only sake of an easy formulation of a uniqueness criterion (see Theorem \ref{th:uniqueness} below). If the reader is not interested in the uniqueness part of Gaussian multiplicative chaos theory, he may skip this assumption of nonnegativeness as we will discuss in section \ref{sec:extent} a more elaborate uniqueness criterion. If $K$ is a  kernel of $\sigma$-positive type with decomposition (\ref{eq:defsigma}), one can consider a sequence of independent centered Gaussian processes $(Y_{k})_{k \geq 1}$ with covariance kernels $(K_k)_k$. Then the Gaussian process 
$$X_n=\sum_{k=1}^nY_k$$ has covariance kernel  $\sum_{k=1}^{n}K_{k}$. Given a Radon measure $\sigma$ on $D$, it is proved in \cite{cf:Kah} that the sequence of random measures $(M_{n})_n$ given by:  
\begin{equation}\label{eq:deformellen}
\forall A \in \mathcal{B}(D), \quad M_{n,\gamma}(A)=\int_{A}e^{ \gamma X_{n}(x)-\frac{\gamma^2}{2}\E[X_{n}(x)^2]}\sigma(dx)
 \end{equation}  
converges almost surely in the space of Radon measures (equipped with the topology of weak convergence) towards a random measure $M$, which is called Gaussian multiplicative chaos\footnote{Private communication with J.¬. Kahane: The terminology multiplicative chaos was adopted since this theory may be seen as a multiplicative counterpart of the additive Wiener chaos theory. Actually, this is Paul L\'evy himself who suggested to J.P. Kahane in the seventies to construct a multiplicative theory of random variables, arguing that this should be as fundamental as the additive theory of random variables. It took almost ten years to Kahane to build his theory of Gaussian multiplicative chaos.} with kernel $K$ acting on $\sigma$. Basically, this convergence relies on the fact that for each compact set $A$, the sequence $(M_{n,\gamma}(A))_n$ is a nonnegative martingale. This martingale structure ensuring the almost sure convergence of \eqref{eq:deformellen} at low cost is the main motivation for considering kernels of $\sigma$-positive type.

 Then Kahane established a whole set of properties of this chaos that he derived from the following comparison principle:
 \begin{theorem}{\bf Convexity inequalities. [Kahane, 1985]. }\label{th:compar}
Let $(A_i)_{1\leq i \leq n}$ and $(B_i)_{1\leq i \leq n}$ be two centered Gaussian vectors such that:
$$\forall i, j,\quad \E[A_iA_j]\leq \E[B_iB_j]. $$
Then  for all combinations of nonnegative weights $(p_i)_{1\leq i \leq n}$ and all convex (resp. concave) functions $F:\R_+\to \R$ with at most polynomial growth at infinity
\begin{equation}\label{eq:compar}
\E\Big[F\big(\sum_{i=1}^np_ie^{A_i-\frac{1}{2}\E[A_i^2]}\big)\Big]\leq (\text{resp. } \geq)\,\E\Big[F\big(\sum_{i=1}^np_ie^{B_i-\frac{1}{2}\E[B_i^2]}\big)\Big].
\end{equation}
\end{theorem}

\begin{rem}
Note that theorem \ref{th:compar} is very general and can be useful for instance in the study of random Gibbs measures where the Hamiltonian is a Gaussian variable. This is for instance the case in the Sherrington-Kirkpatrick (SK) model of statistical physics. In an important work, the authors of \cite{GuTo} rederived (\ref{eq:compar}) with $F(x)= \ln x$ in the context of the SK model and used this inequality to compare the model of size $N$ with two independent subsystems of size $N_1$ and $N_2$ with $N=N_1+N_2$. As a consequence, they noticed that the expected finite size free energy was subadditive with respect to its size hence obtaining the existence of the limiting free energy as the size of the system goes to infinity.     
\end{rem}

This ingenious inequality sheds some light on the mechanism of Gaussian multiplicative chaos. For instance, let us stress that a kernel $K$ of $\sigma$-positive type admits infinitely many decompositions of the form \eqref{eq:defsigma}: you can obtain other decompositions by changing the order of the kernels $K_k$, by gathering them, etc... so there are possibly quite different kernels $(K_k')_k$ whose sum is $K$. The important question thus is: does the law of the limiting measure $M$ depend on the choice of the decomposition $(K_k)_k$ in \eqref{eq:defsigma}?

\begin{theorem}{\bf Uniqueness. [Kahane, 1985]. }\label{th:uniqueness}
The law of the limiting measure $M_{\gamma}$ does not depend on the sequence of nonnegative and positive definite kernels $(K_{k})_{k \geq 1}$ used in the decomposition (\ref{eq:defsigma}) of $K$.
\end{theorem}
Thus, the theory enables to give a unique and mathematically rigorous definition to a random measure $M$ in $D$ defined formally by:
\begin{equation}\label{eq:deformelle}
\forall A \in \mathcal{B}(D), \quad M_\gamma(A)=\int_{A}e^{\gamma X(x)-\frac{\gamma^2}{2}E[X(x)^2]}\,\sigma(dx).
 \end{equation}  
where $(X(x))_{x \in D}$ is a centered "Gaussian field" whose covariance $K$ is a $\sigma$-positive type kernel.  
To show the usefulness of Theorem \ref{th:compar}, it is worth giving a few words about the proof of Theorem \ref{th:uniqueness}. It works roughly as follows. Assume that you have two decompositions $(K_{k})_{k \geq 1}$ and $(K_{k}')_{k \geq 1}$ of $K$ with associated Gaussian process sequences $(X_n)_n$ and $(X_n')_n$ and associated measures  $(M_n)_n$ and $(M_n')_n$. Both sequences $(\sum_{k=1}^{n}K_{k})_n$ and $(\sum_{k=1}^{n}K_{k}')_n$ converge pointwise towards $K$ in a nondecreasing way. Therefore, if we choose a compact set $T\subset D$ then,  for each fixed $p\geq 1$ and $\epsilon >0$,  the Dini theorem entails that $$\sum_{k=1}^{p}K_{k}\leq \epsilon+\sum_{k=1}^{q}K_{k}'$$ for $q$ large enough on $T\times T$. Since $\sum_{k=1}^{q}K_{k}'$ (resp. $\sum_{k=1}^{p}K_{k}$) is the covariance kernel of $X'_q$ (resp. $X_p$), we can apply Kahane's convexity inequalities and get, for each bounded convex function $F:\R_+\to\R$:
$$\E[F(M_p(A))]\leq \E[F(e^{\sqrt{\epsilon}\gamma Z-\frac{\epsilon \gamma^2}{2}}M'_q(A))],$$ where $Z$ is a standard Gaussian random variable independent of $M'$.
By taking the limit as $q$ tends to $\infty$, and then $p\to \infty$, we obtain
$$\E[F(M (A))]\leq \E[F(e^{\sqrt{\epsilon}\gamma Z-\frac{\epsilon \gamma^2}{2}}M' (A))].$$ Since $\epsilon>0$ can be chosen arbitrarily small, we deduce 
$$\E[F(M (A))]\leq \E[F(M' (A))].$$
The converse inequality is proved in the same way, showing $\E[F(M (A))]= \E[F(M' (A))]$ for each bounded convex function $F$. By choosing $F(x)=e^{-\lambda x}$ for $\lambda>0$, we deduce that the measures $M$ and $M'$ have the same law. 

\vspace{2mm}
However, the simplicity of the convergence does not solve the question of non-degeneracy of the limiting measure $M$: it is possible that $M$ identically vanishes. A $0-1$ law argument straightforwardly shows that the event "$M$ is identically null" has probability $0$ or $1$. It seems difficult to state a general decision rule to decide whether $M$ is degenerate or not.  It depends in an intricate way on the covariance structure, i.e. the kernel $K$, and on the measure $\sigma$. So Kahane focused on the situation when the kernel $K$ and the measure $\sigma$ are intertwined via the metric structure of $D$.
More precisely, he assumed that $K$ can be written as
\begin{equation}\label{Klog}
\forall x,y\in D,\quad K(x,y)=\ln_+ \frac{T}{\rho(x,y)}+g(x,y)
\end{equation}
where $T>0$, $g:D\times D\to\R$ is a bounded continuous function and $\sigma$ is in the class $R_\alpha^+$ (denoted $M_{\alpha^+}$ in Kahane's paper):
\begin{definition}
For $\alpha>0$, a Borel measure $\sigma$ is said to be in the class $R^+_\alpha$  if for all $\epsilon>0$ there is $\delta>0$, $C<\infty$ and a compact set $A_\epsilon\subset D$ such that $\sigma(D\setminus A_\epsilon)\leq \epsilon$ and: 
\begin{equation}
\forall O \text{ open set}, \quad \sigma(O\cap A_\epsilon)\leq C_\epsilon{\rm diam}_\rho(O)^{\alpha+\delta},
\end{equation}
where ${\rm diam}_\rho(O)$ is the diameter of $O$ with respect to $\rho$.
\end{definition}
For instance, the Lebesgue measure of $\R^d$, restricted to any bounded domain of $\R^d$ equipped with the Euclidean distance $\rho$, is in the class $R_\alpha^+$ for all $\alpha<d$. 

The above definition looks like a H\"older condition for measures. It is intimately related to the notion of measure  with finite $\beta$-energy: a  Borel measure $\sigma$ is said to be of finite $\beta$-energy if
\begin{equation}\label{betaenergy}
I_\beta(\sigma)=\int\int\frac{1}{\rho(x,y)^\beta}\sigma(dx)\sigma(dy)<+\infty.
\end{equation}
Indeed, if $\sigma$ has a finite   $\beta$-energy then $\sigma\in R^+_\alpha$ for all $\alpha<\beta$.
Conversely, if $\sigma$ is in the class $R^+_\alpha$, then the measure $\sigma_{A_\epsilon}(dx)=\ind_{A_\epsilon}(x)\sigma(dx)$ has finite $\beta$-energy for all $\beta<\alpha+\epsilon$. 

To have a flavor of  the forthcoming results, let us treat the following simple situation, which we call 
 the "below $L^2$-threshold" case. It is about formulating a criterion ensuring that the martingale $(M_n(A))_n$ is bounded in $L^2$ for some given bounded set $A$, and therefore uniformly integrable.  A straightforward computation shows that:
 \begin{align*}
\E[M_{n,\gamma}(A)^2]=&\int_A\int_A\E[e^{ \gamma X_{n}(x)-\frac{\gamma^2}{2}\E[X_{n}(x)^2]}e^{ \gamma X_{n}(y)-\frac{\gamma^2}{2}\E[X_{n}(y)^2]}]\,\sigma(dx)\sigma(dy)\\
\leq &\int_A\int_A e^{\gamma^2 K(x,y)}\,\sigma(dx)\sigma(dy)\\
\leq & C\int_A\int_A \frac{1}{\rho(x,y)^{\gamma^2}}\,\sigma(dx)\sigma(dy).
\end{align*}
Therefore, if the measure $\ind_A(x)\sigma(dx)$ has finite $\gamma^2$-energy, the martingale $(M_n(A))_n$ is bounded in $L^2$ and therefore converges towards a non trivial limit. In the case when $\sigma$ is the Lebesgue measure of $\R^d$ and $\rho$ the Euclidean distance, this condition simply reads $\gamma^2<d$.

Kahane proved the following highly deeper result: 
\begin{theorem}{\bf Non-degeneracy. [Kahane, 1985]. }\label{th:degeneracy}
Assume that the kernel $K$ takes on the form \eqref{Klog} and that the measure $\sigma$ is in the class $R_\alpha^+$ for some $\alpha>0$. Then, for each compact set $A$, the sequence $(M_{n,\gamma}(A))_n$ is a uniformly integrable martingale. Hence
$$\gamma^2<2\alpha \Rightarrow M_\gamma \text{ is non degenerate}.$$
\end{theorem}
As a by-product of his proof, Kahane also shows the following result concerning the dimension of the carrier of the measure $M_\gamma$:

\begin{theorem}{\bf Structure of the carrier. [Kahane, 1985]. }\label{th:carrier}
Assume that the kernel $K$ takes on the form \eqref{Klog} and that the measure $\sigma$ is in the class $R_\alpha^+$ for some $\alpha>0$. If $\gamma^2<2\alpha$, the measure $M_\gamma$ is non degenerate and is almost surely in the class $R^+_{\alpha-\frac{\gamma^2}{2}}$.
\end{theorem}
Note that Theorem \ref{th:carrier} implies that the measure $M$ cannot give positive mass to a set of Hausdorff dimension less or equal than $\alpha-\frac{\gamma^2}{2}$. Therefore the measure $M$ cannot possess atoms if $\sigma$ is in some class $R_\alpha^+$. We will see in subsection \ref{sub:thick} that the Hausdorff dimension of the carrier is exactly $d-\frac{\gamma^2}{2}$ when $\sigma$ is the Lebesgue measure on a domain of $\R^d$.

\begin{rem}
Stronger versions of Theorem \ref{th:degeneracy} and Theorem \ref{th:carrier} were proved recently in \cite{Rnew11}.
\end{rem}

\begin{rem}
If $\sigma$ is the Lebesgue measure on some domain $D \subset \R^d$ equipped with the Euclidean metric, then $\sigma$ is in the class $R_{d-\epsilon}^+$ for all $\epsilon>0$ hence, if $\gamma^2<2 d$, the measure $M_\gamma$ is non degenerate and in the class $R^+_{d-\frac{\gamma^2}{2}-\epsilon}$ for all $\epsilon>0$.
\end{rem}

Partial converses of Theorem \ref{th:degeneracy} are more intricate. Kahane first gave a general necessary condition:
\begin{theorem}{\bf Necessary condition of non-degeneracy. [Kahane, 1985]. }\label{th:necess1}
Assume $(D,\rho)$ is a locally compact metric space and:\\
-the function $(t,s)\mapsto \rho(t,s)^2$ is of negative type,\\
-$\sigma$ has the doubling property, namely that there exists a constant $C$ such that
$$\forall x\in T,\forall r>0, \quad \sigma (B(x,2r))\leq C \sigma(B(x,r)), $$
-the kernel $K$ takes on the form \eqref{Klog}.\\
 Denote by ${\rm dim}(D)$ the Hausdorff dimension of $D$. If $\gamma^2> 2 {\rm dim}(D)$ then  $M_\gamma$  is   degenerate.
\end{theorem}
As pointed out by Kahane, for a squared distance of negative type, the assumptions of the above theorem  are satisfied when the triple $(D,\rho,\sigma)$ admits a Lipschitz immersion into a finite-dimensional space. Let us also stress that the  critical situation $\gamma^2=2{\rm dim}(D)$ is not settled by this theorem.
Nevertheless, he reinforced his assumptions to prove:

\begin{theorem}{\bf Necessary and sufficient condition of non-degeneracy. [Kahane, 1985]. }\label{th:necess2}
Assume $(D,\rho)$ is a $d$-dimensional  manifold of class $C^1$ and let $\sigma$ be its volume form (or any Radon measure absolutely continuous w.r.t  the volume form with a bounded density). Assume that the kernel $K$ takes on the form \eqref{Klog}. Then
$$ M_\gamma \text{ is non degenerate}\Rightarrow \gamma^2<2d .$$
\end{theorem}

\subsubsection*{More results in the Euclidean space}

When the metric space $(D,\rho)$ is an open subset of $\R^d$ for some $d\geq 1$ equipped with the Euclidian distance and $\sigma$ is the Lebesgue measure, the previous results can be strengthened. We first point out that Theorem \ref{th:necess2} applies and the non-degeneracy necessary and sufficient condition reads $\gamma^2<2d$. 

\begin{figure}[h]
\centering
\subfloat[$\gamma=0.2$]{\includegraphics[width=0.47\linewidth]{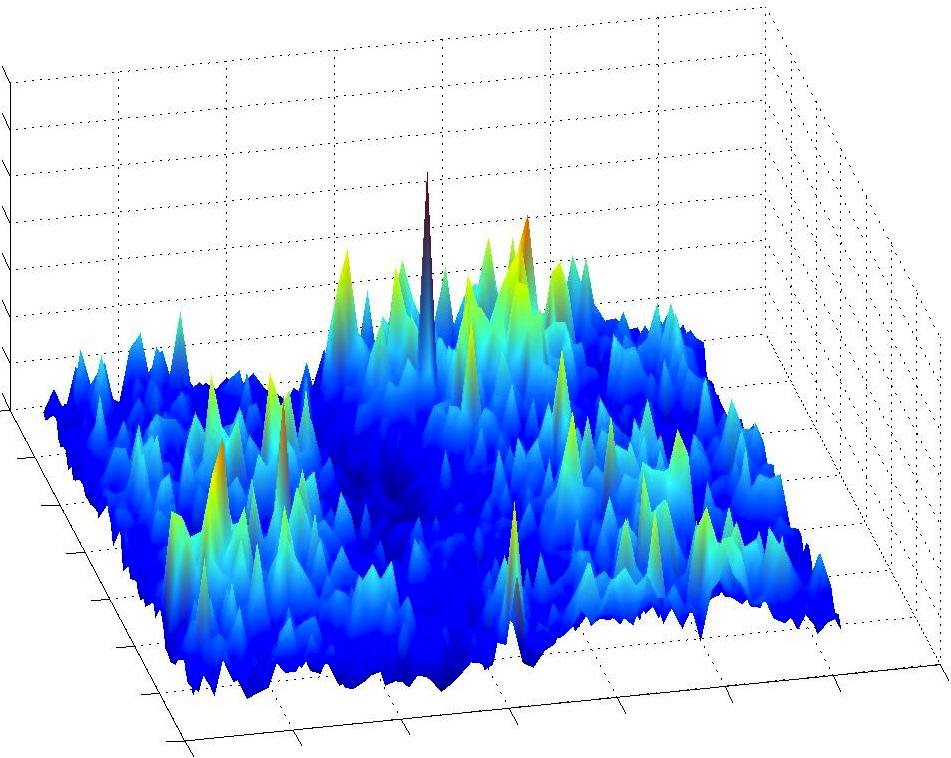}}
\,\,\subfloat[$\gamma=1$]{\includegraphics[width=0.47\linewidth]{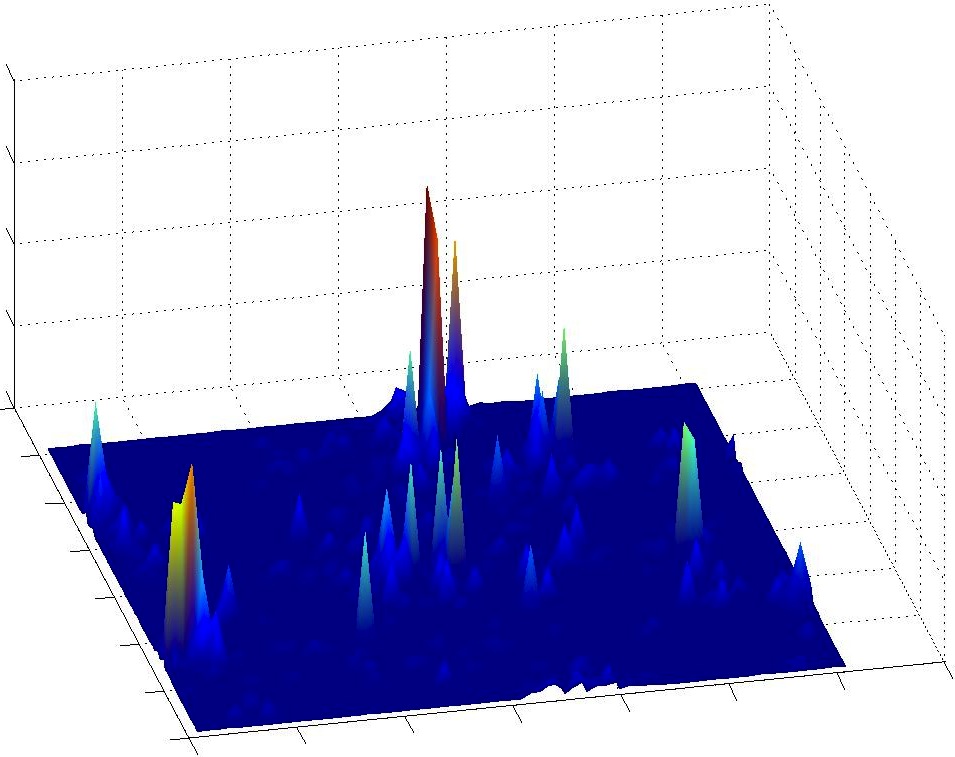}}\\
\subfloat[$\gamma=1.8$]{\includegraphics[width=0.47\linewidth]{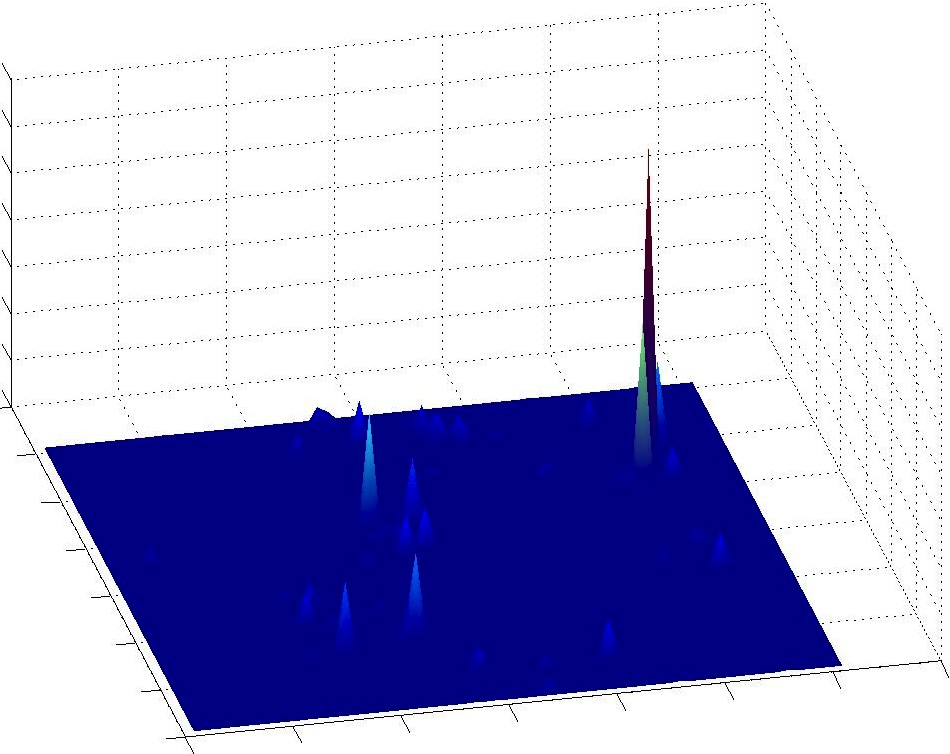}}
\caption{Influence of the intermittency parameter $\gamma$: we plot the "density profile" of a $2$-dimensional multiplicative chaos illustrating the "clustering" of the activity of the chaos when $\gamma$ grows.}
\label{interm}
\end{figure}
 
Kahane's convexity inequalities (Theorem \ref{th:compar}) allow us to give a complete description of the moments of $M_\gamma$:
\begin{theorem}{\bf Positive Moments. [Kahane, 1985]. }\label{th:moments}
If the measure $M_\gamma$ is non degenerate, that is $\gamma^2<2d$, the measure $M_\gamma$ admits finite positive moments of order $p$ for all $p\in ]0,\frac{2 d}{\gamma^2}[$. More precisely, for all compact set $A\subset D$ and $p\in ]0,\frac{2d}{\gamma^2}[$, we have $\E[M_\gamma(A)^p]<+\infty$.
\end{theorem}
 
Basically, it suffices to prove the finiteness of the moments for your favorite kernel $K$ of the type \eqref{Klog} and deduce that the conclusions remain valid for all the kernels of the type \eqref{Klog} via Theorem \ref{th:compar}. 

We know turn to the existence of negative moments which was not investigated by Kahane. We have:

\begin{theorem}{\bf Negative Moments.}\label{th:momentsneg}
If the measure $M_\gamma$ is non degenerate, that is $\gamma^2<2d$, the measure $M_\gamma$ admits finite negative moments of order $p$ for all $p\in ]-\infty,0[$. More precisely, for any compact nonempty Euclidean ball $A\subset D$ and $p\in ]-\infty,0[$, we have $\E[M_\gamma(A)^p]<+\infty$.
\end{theorem}

The finiteness of moments of negative order is proved in \cite{Molchan} in the case of discrete cascades. By adapting the argument of \cite{Molchan} and using the convexity inequalities \ref{eq:compar}, theorem \ref{th:momentsneg} is proved in \cite{cf:RoVa}.

Let us also point out an important result in \cite{BXJ} where the authors compute the tail distributions of the measure $M_\gamma$ in dimension $1$ with the kernel $K(x,y)=\ln_+\frac{1}{|x-y|}$ (which is of $\sigma$-positive type, see Proposition \ref{d=12} below):
\begin{theorem}{\bf Distribution tails. [Barral, Jin, 2012]. }
If $A$ is some nonempty segment of $\R$ then there exists a constant $c>0$ such that
$$\lim_{x\to +\infty}x^{\frac{2}{\gamma^2}}\Pb(M_\gamma(A)>x)=c.$$
\end{theorem}
We conclude this theoretical background by pointing out some further interesting properties that can be exhibited when the integrating measure $\sigma$ in \eqref{eq:deformelle} is the Lebesgue measure. The first observation that can be made is that the measure $M_\gamma$ is stationary in space as soon as we consider a stationary Gaussian distribution $X$. This is due to the translation invariance of the Lebesgue measure and the stationarity of $X$. Furthermore, when $M$ is non-degenerate, it can be shown that the support of $M_\gamma$ is almost surely the whole of $\R^d$. This results from the $0-1$ law of Kolmogorov: if you consider a ball $B$, the  $0-1$ law tells you that the event $\{M_\gamma(B)>0\}$ has probability $0$ or $1$. Indeed, we have
\begin{equation}\label{oun}
\inf_{x\in B}e^{\gamma X_n(x)-\frac{\gamma^2}{2}\E[X_n^2(x)]}M_{\gamma,n}(B)\leq  M_\gamma(B)\leq \sup_{x\in B}e^{\gamma X_n(x)-\frac{\gamma^2}{2}\E[X_n^2(x)]}M_{\gamma,n}(B),
\end{equation}
 where $M_{\gamma,n}(B)$ is the Gaussian multiplicative chaos
$$M_{\gamma,n}(dx)=\lim_{k\to\infty}e^{\gamma (X_k-X_n)(x)-\frac{\gamma^2}{2}\E[(X_k-X_n)^2(x)]}\,dx.$$
Since for $k> n$:
$$X_k-X_n=\sum_{p=n+1}^kY_p,$$ we may say that we have just removed the dependency on the first $n$ fields $Y_1,\dots,Y_n$.
Therefore \eqref{oun} entails that for any $n\geq 0$
$$\{M_{\gamma}(B)>0\}=\{M_{\gamma,n}(B)>0\}.$$
Since the event $\{M_{\gamma,n}(B)>0\}$ is independent of the fields $Y_1,\dots,Y_n$, we deduce that the event $\{M_{\gamma}(B)>0\}$ belongs to the asymptotic sigma-algebra generated by the fields $(Y_n)_n$ in such a way that it has probability $0$ or $1$. When $M_\gamma$ is non-degenerate, it has clearly probability $1$. This argument can be reproduced for every ball chosen among a countable family of balls $(B_n)_n$ generating the open sets of $\R^d$. Therefore, the event $\bigcap_n\{M_\gamma(B_n)>0\}$ has probability $1$, proving that almost surely the support of $M_\gamma$ is the whole of $\R^d$.

\subsection{Examples of kernels of $\sigma$-positive type}\label{sec:example}
In this section, we give a few important examples of $\sigma$-positive kernels.  

\subsubsection{Exact kernels}
 
We consider for  $T>0$ the kernel
\begin{equation}\label{formln}
\forall x,y\in\R^d,\quad K(x,y)=\ln_+\frac{T}{|x-y|}.
\end{equation}
It is of $\sigma$-positive type in dimension $d=1,2$ and is involved in exact scaling relations as explained in subsection \ref{sec:multi}.
 
\subsubsection{Star scale invariant kernels}
A simple way of constructing $\sigma$-positive kernels on $\R^d$ is to consider
\begin{equation}\label{star}
\forall x,y\in\R^d,\quad K(x,y)=\int_1^{\infty}\frac{k\big((x-y)u\big)}{u}\,du,
\end{equation}
 where $k$ is a continuous function of positive type with $k(0)=1$.  Such kernels are related to the notion of $\star$-scale invariance (see subsection \ref{sec:multi}). Whole plane massive Green functions are $\star$-scale invariant.

\subsubsection{Green functions}
If we consider a bounded domain $D$ of $\R^2$, the Green function $G$ of the Laplacian with $0$-boundary condition is of $\sigma$-positive type. The corresponding Gaussian distribution $X$ with covariance $G$ is the Gaussian Free Field (GFF for short). The associated Gaussian multiplicative chaos is called the Liouville measure. More details about these claims are given in subsection \ref{LQG}.

\subsection{Different notions of stochastic scale invariance}\label{sec:multi}
In this subsection, we consider the Euclidian framework. More precisely, we consider  an open set $D$ of $\R^d$ and  Gaussian multiplicative chaos of the type
\begin{equation}\label{euclid:M}
\forall A\in \mathcal{B}(D),\quad M_\gamma(A)=\int_A e^{\gamma X(x)-\frac{\gamma^2}{2}\E[X(x)^2]}\,dx
\end{equation}
where the Gaussian distribution $X$ has covariance kernel of the form:
\begin{equation}\label{euclid:K}
K(x,y)= \ln_+\frac{T}{|x-y|}+g(x,y)
\end{equation}
for some continuous and bounded function $g$ over $D^2$. Then the power-law spectrum of  such Gaussian multiplicative chaos  presents some interesting features, such as  non-linearity (in the parameter $q$ in the following theorem):
\begin{theorem}\label{th:spectrum}
Assume that the kernel $K$ takes on the form \eqref{Klog} with $\gamma^2<2d$. Choose a point $x\in D$. For each $q\in [0,\frac{2d}{\gamma^2}[$ we have
$$\E\big[M_\gamma(B(x,r))^q\big]\asymp r^{\xi(q)}\quad \text{as }r\to 0$$ where $\xi$ is the {\it structure exponent} of the measure $M$:
$$\forall q \in [0,\frac{2d}{\gamma^2}[,\quad \xi(q)=\big(d+\frac{\gamma^2}{2}\big)q-\frac{\gamma^2}{2}q^2.$$
\end{theorem}

\noindent {\it Heuristic proof.} For simplicity, assume that $T=1$ in \eqref{euclid:K}. By making a change of variables, we get:
\begin{align*}
M_\gamma(B(x,r)) =& \int_{B(x,r)}e^{\gamma X(y)-\frac{\gamma^2}{2}\E[X(y)^2]}\,dy\\
=& r^d\int_{B(0,1)}e^{\gamma X(x+ry)-\frac{\gamma^2}{2}\E[X(x+ry)^2]}\,dy.
\end{align*}
Then we observe that the field $(X(x+ry))_y$ has a covariance structure approximatively given for $r\leq 1$ by:
\begin{align*}
\E[X(x+ry)X(x+rz)]\simeq & \ln_+\frac{1}{r|y-z|}\\
=&\ln\frac{1}{r}+\ln_+\frac{1}{|y-z|}\\
\simeq & \ln\frac{1}{r}+\E[X(x+y)X(x+z)].
\end{align*}
The above relation gives us the following (good) approximation in law
$$ (X(x+ry))_{y\in B(0,1)}\simeq\Omega_r+(X(x+y))_{y\in B(0,1)}$$
where $\Omega_r$ is a centered Gaussian random variable independent of the field $(X(x+y))_{y\in B(0,1)}$ and with variance $\ln\frac{1}{r}$. Therefore
\begin{align*}
M_\gamma(B(x,r)) =&  r^d\int_{B(0,1)}e^{\gamma X(x+ry)-\frac{\gamma^2}{2}\E[X(x+ry)^2]}\,dy\\
\simeq & r^d e^{\gamma \Omega_r-\frac{\gamma^2}{2}\E[\Omega_r^2]}\int_{B(0,1)}e^{\gamma X(x+y)-\frac{\gamma^2}{2}\E[X(x+y)^2]}\,dy\\
=&  r^d e^{\gamma \Omega_r-\frac{\gamma^2}{2}\E[\Omega_r^2]}M_\gamma (B(x,1)).
\end{align*}
By taking the $q$-th power and integrating, we get as $r\to 0$
$$\E\big[M_\gamma(B(x,r))^q\big]\simeq  r^{\xi(q)}\E[M_\gamma(B(x,1)^q],$$
thus explaining the theorem. Actually, the rigorous proof of this result is very close to the heuristic developed here.\qed

Notice that the quadratic structure of the structure exponent is intimately related to the Gaussian nature of the random distribution $X$. Random measures with a non-linear power-law spectrum are often called {\it multifractal}. That is why Gaussian multiplicative chaos (and other possible extensions) are sometimes called Multifractal Random Measures  (MRM for short) in the literature. It is also natural to wonder if some specific choice of the covariance kernel $K$ may lead to replacing the symbol $\asymp$ in Theorem \ref{th:spectrum} by the symbol $=$. It turns out that this question is related to some specific scaling relation, which we describe below.

\subsubsection*{Exact stochastic scale invariance}
The notion of "exact stochastic scale invariance" relies on the additive properties of the logarithm function. Roughly speaking, in order to produce kernels with exact scaling relations, kernels $K$ of the type \eqref{Klog} with $g=0$ (and $\rho$ the Euclidian distance) must be considered. So we first focus on the $\sigma$-positive type of such kernels:

\begin{proposition}\label{d=12}
 For $d\leq 2$ and $T>0$, the function
 \begin{equation}\label{eq:exact12}
 x\in\R^d\mapsto \lambda^2\ln_+\frac{T}{|x|}
 \end{equation}
  is of $\sigma$-positive type. 
 \end{proposition} 
 
 \noindent {\it Proof.} A straightforward computation yields:
 $$\ln_+\frac{T}{|x|}=\int_0^{+\infty}(t-|x|)_+\nu_T(dt)$$
 where $\nu_T$ is the measure ($\delta_T$ is the Dirac mass at $T$):
 $$ \nu_T(dt)=\ind_{[0,T]}(t)\frac{dt}{t^2}+\frac{1}{T}\delta_T(dt).$$ 
Hence for any $\mu>0$, we have:
$$\ln_+\frac{T}{|x|}=\frac{1}{\mu}\ln_+\frac{T^\mu}{|x|^\mu}=\int_0^{+\infty}(t-|x|^\mu)_+\nu_{T^\mu}(dt).$$
By using a Chasles relation in the integral of the right-hand side, proving that this kernel is of $\sigma$-positive type thus boils down to  considering the possible values of  $\mu>0$ such that the function  $(1-|x|^\mu)_+$ is of positive type: this is the Kuttner-Golubov problem (see \cite{golubov}). 

For $d=1$, it is straightforward to see that $(1-|x|)_+$ is of positive type (compute the inverse Fourier transform). In dimension 2, Pasenchenko \cite{cf:PaYu} proved that the function $(1-|x|^{1/2})_+$ is of positive type on  $\R^2$. We can thus write
$$\ln_+\frac{T}{|x|}=\sum_{n\geq 1}K_n(x)$$ with
$$K_n(x)=\int_{\frac{1}{n}}^{\frac{1}{n-1}}(t-|x|^\mu)_+\nu_{T^\mu}(dt)$$
with $\mu=1$ in dimension $1$ and $\mu=1/2$ in dimension $2$.\qed
 
\begin{theorem}{\bf Exact stochastic scale invariance. {\bf [Bacry, Muzy 2003].}}\label{th:exact12}
Let  $K$ be the covariance kernel given by \eqref{eq:exact12} in dimension $d=1$ or $d=2$.
The associated Gaussian multiplicative chaos $M_\gamma$   is exactly stochastic scale invariant:
\begin{equation}\label{essi}
\forall  \lambda\in]0,1],\quad (M(\lambda A))_{A\subset B(0,T/2)}\stackrel{law}{=}  \lambda^d e^{\Omega_\lambda-\frac{1}{2}\E[\Omega_\lambda^2]}(M(  A))_{A\subset B(0,T/2)},
\end{equation}
 where $\Omega_\lambda$ is a Gaussian random variable,  independent of the measure $(M(  A))_{A\subset B(0,T/2)}$, with mean $0$   and variance $\gamma^2\ln \frac{1}{\lambda}$ .
 \end{theorem}  

We stress that the above equality in law is to be understood in the sense
$$\big(M(\lambda A_1),\dots,M(\lambda A_p)\big)\stackrel{law}{=}\big(\lambda^d e^{\Omega_\lambda-\frac{1}{2}\E[\Omega_\lambda^2]}M(  A_1),\dots,\lambda^d e^{\Omega_\lambda-\frac{1}{2}\E[\Omega_\lambda^2]}M(  A_p)\big)$$ for all possible choice $A_1,\dots,A_p$ of Borelian subsets of the ball $B(0,T/2)$.

Notice that such a scaling property makes obvious several computations related to the measure $M_\gamma$.
The reader may, for instance, observe that the heuristic proof of Theorem \ref{th:spectrum} becomes rigorous for such a measure. It is also not difficult to see that this scaling relation is bound to be valid only locally (over a ball) and cannot hold on the whole space: the logarithm is not of positive type over the whole of $\R^d$. Exact scaling relations were introduced in \cite{cf:BaMu} in dimension $1$ together with further generalizations in the case of log-infinitely divisible random measures. 

It is natural to wonder how to construct such measures in dimension higher than $3$. The procedure is somewhat complicated by the following observations:  for $d= 3$, it is an open question to know whether the kernel \eqref{eq:exact12} is of $\sigma$-positive type and for $d\geq 4$, it is even not of positive type. 

Another approach has therefore been suggested in \cite{Rnew9}. The main ideas are the following: what only matters to construct exactly stochastic scale invariant measures is that the covariance kernel must be the logarithm function over a ball centered at $0$. The way of "truncating" the logarithm does not matter in order to obtain the scaling relation \eqref{essi} but may be sensitive to the dimension when regarding positive definiteness: truncating the logarithm with the function $\ln_+$ does not resist   increasing the dimension. In  \cite{Rnew9}, another truncation is suggested:  we can find an isotropic function  $g:\R^d\to\R$ that is constant on a neighborhood of $0$ and such that  the kernel
 \begin{equation}\label{eq:kmmrm}
K(x)=\ln_+\frac{T}{|x|}+g(x)
\end{equation}
is of $\sigma$-positive type. Briefly, the construction is the following: let us denote by $S$ the sphere of $\R^d$ and by $ \sigma$ the unique uniform measure on the sphere such that $\sigma(S)=1$. This (probability) measure is invariant under rotations. Let us  define the function
\begin{equation}\label{def:F}
K(x)=\int_S\ln_+\frac{T}{|\ps{x,s}|}\sigma(ds),
\end{equation}
where $\ps{\cdot,\cdot}$ stands for the canonical inner product of $\R^d$.
Since  $\sigma$ is invariant under rotations, the function $K$ is isotropic. 
Fix $x\in\R^d$ such that $|x|\leq T$ and write $x=|x|e$ where $e\in S$. Then we have
$$K(x)=\int_S\ln\frac{T}{|x||\ps{e,s}|}\sigma(ds)=\gamma^2\ln\frac{T}{|x|}+\gamma^2\int_S\ln\frac{1}{|\ps{e,s}|}\sigma(ds).$$
By invariance under rotations of $\sigma$, the second term in the right-hand side does not depend on $x$ and turns out to be  finite: this can be seen by noticing that, under $\sigma$, the random variable $\ps{e,s}$ has the law of the first entry of a Haar vector.  $K$ thus coincides with the logarithm over a neighborhood of $0$, up to an additive constant. Gaussian multiplicative chaos with associated kernel \eqref{def:F} are exactly stochastic scale invariant in the sense of \eqref{essi}. It remains an open question to know to which extent the notion of exact stochastic scale invariance uniquely determines the covariance structure of the associated kernel $K$.
 
\subsubsection*{Star scale invariance}
As explained  above, the notion of exact stochastic scale invariance is a local notion (valid only over a ball). We  now present a global notion of stochastic scale invariance, called star scale invariance in \cite{Rnew1} in reference of earlier works by Mandelbrot in the case of multiplicative cascades on trees. It stems from the need  of characterizing  Gaussian multiplicative chaos with  functional equations. Indeed, tractable functional equations  may  provide efficient tools in identifying Gaussian multiplicative chaos as, for instance,  scaling limits of discrete models.

\begin{definition}{\bf Log-normal $\star$-scale invariance.} \label{def:star}
A random Radon measure $M$ on $\R^d$ is said lognormal $\star$-scale invariant if for all $0<\varepsilon\leq 1$, $M $ obeys the cascading rule
\begin{equation}\label{starc}
 \big(M (A)\big)_{A\in\mathcal{B}(\R^d)}\stackrel{law}{=} \big(\int_Ae^{ \omega_{\varepsilon}(x)}M^\varepsilon(dx)\big)_{A\in\mathcal{B}(\R^d)}
\end{equation}
where $\omega_{\varepsilon}$ is a Gaussian process, which is assumed to be stationary and continuous in probability, and $M^{ \varepsilon}$ is a random measure independent from $X_{\varepsilon}$ satisfying the scaling relation
\begin{equation}\label{star1c}
\hspace{4cm} \big(M^{\varepsilon}(A)\big)_{A\in\mathcal{B}(\R^d)}\stackrel{law}{=}  \big(M (\frac{A}{\varepsilon})\big)_{A\in\mathcal{B}(\R^d)}.\hspace{2cm}\qed
\end{equation} 
\end{definition}

Notice  that the process $ \omega_{\varepsilon}$ is unknown. Roughly speaking, we look for  random measures that scale with an independent lognormal factor on the whole space. This property is shared by a large class of Gaussian multiplicative chaos. And for those Gaussian multiplicative chaos that do not share this property, they are very close to satisfying it. If the reader is familiar with branching random walks (BRW), here is an explanation that may help intuition. If we consider a BRW the reproduction law of which does not change with time (i.e. is the same at each generation), the law of the branching random walk will be characterized by a discrete version of the above $\star$-scale invariance called "fixed point of the smoothing transform"  (in the lognormal case of course, see \cite{durrett,Biggf,Liu}). If the reproduction law evolves in time, then we have to change things a bit to adapt to this time evolution. The same argument holds for the log-normal $\star$-scale invariance: it characterizes these Gaussian multiplicative chaos that do not vary along scales.

It is proved in \cite{Rnew1}  that $\E[e^{ \omega_{\varepsilon}(r)}]=\varepsilon^d$ as soon as the measure $M$ possesses a moment of order $1+\delta$ for some $\delta>0$. Furthermore, up to weak regularity conditions on the covariance kernel of the process $ \omega_{\varepsilon}$ summarized in the definition below, all the log-normal $\star$-scale invariant random measures with enough moments can be identified.
\begin{definition}\label{goodness}
We will say that a stationary random measure $M$ satisfies the good lognormal $\star$-scale invariance  if $M$ is lognormal $\star$-scale invariant and  for each $\epsilon<1$, the covariance kernel $K_{\epsilon}$ of the 
process $\omega_{\epsilon}$ involved in (\ref{starc})  is continuous and satisfies:
\begin{align}
 & &|k_{\epsilon}(r)|&\to 0\quad \text{ as }\quad |r|\to+\infty,\label{lim}\\
\forall r,r'\in\R^d\setminus\{0\}, & & |k_{\epsilon}(r)-k_{\epsilon}(r')|&\leq  C_{\epsilon}\theta\big(\min(|r|,|r'|)\big)|r-r'|\label{lipschitz}
\end{align} 
 for some positive constant $C_{\epsilon}$ and some  decreasing function $\theta:]0,+\infty[\to \R_+$ such that 
\begin{equation}\label{cvint}
\int_1^{+\infty}\theta(u)\ln(u)\,du<+\infty.\qed
\end{equation}
\end{definition}
Lognormal $\star$-scale invariant random measures are then characterized as:
\begin{theorem}{\bf [Allez, Rhodes, Vargas, 2011]}\label{th:star}
Let $M$ be a  good lognormal $\star$-scale invariant stationary  random measure. Assume that $$\E[M([0,1]^d)^{1+\delta}]<+\infty$$ for some $\delta>0$. Then $M$ is the product of a nonnegative random variable $Y\in L^{1+\delta}$  and an independent Gaussian multiplicative chaos: 
\begin{equation}\label{chaosmulti}
\forall A\subset \mathcal{B}(\R^d),\quad M(A)=Y\int_Ae^{\gamma X(x)-\frac{\gamma^2}{2}\E[X(x)^2]}\,dx
\end{equation}
 with associated  covariance kernel given by the improper integral
\begin{equation}\label{structmulti}
\forall x\in \R^d\setminus\{0\},\quad K(x)=\int_{1}^{+\infty}\frac{k(x u)}{u}\,du
\end{equation} for some continuous covariance function $k$ such that $k(0)=1$ and   $\gamma^2 \leq \frac{2d}{1+\delta}$. 

Conversely, given some datas $k$ and $Y$ as above, the relation \eqref{chaosmulti} defines a  lognormal $\star$-scale invariant random measure $M$ with finite  moments of order $1+\beta$ for every $\beta \in [0, \delta)$.  
\end{theorem}
It is plain to see that the covariance structure \eqref{structmulti} can be rewritten as 
\begin{equation}\label{starK}
\forall x\in \R^d\setminus\{0\},\quad K(x)= \ln_+\frac{1}{|x|}+g(x)
\end{equation} for some continuous bounded function $g$, thus making the connection with Kahane's theory presented in subsection \ref{seminal}.

The first star-scale invariant kernel appearing in the literature goes back to Kahane's original paper \cite{cf:Kah} in order to approximate the kernel $\ln_+\frac{T}{|x-y|}$ in dimension $3$ but he did not study the related scaling relations. Kahane chose the kernel
$$K(x)=\int_{1}^{+\infty}\frac{e^{-u|x|}}{u}\,du.$$

Another example was exhibited in \cite{Bar}. It corresponds to the choice
$$K(x)=\int_{1}^{+\infty}\frac{k(ux)}{u}\,du\quad \text{with}\quad k(x)=(1-\frac{|x|}{T})\ind_{[0,T]}(|x|)$$ and is based on a one dimensional geometric construction. Also, the authors in \cite{Bar}  made the connection with star scale invariance.
 
To sum up, the solutions of the star-equation are now well identified in what is called the subcritical regime, which can be characterized by the fact that associated solutions possess enough moments (larger than $1$).
It remains  to investigate the characterization of solutions with only few moments (smaller than $1$): we will see  in section \ref{sec:gener} that new types of solutions are then involved, giving rise to quite new structures of the chaos.  %
\section{Extensions of the theory}\label{sec:extent}

\subsection{Limitations of Kahane's theory}
In view of natural applications, Kahane's theory appears unsufficient. Here are a few points that the theory does not address:
\begin{itemize}
\item
A kernel $K$ of $\sigma$-positive type is nonnegative and positive definite. Is the reciprocal true? 
\item
What happens if one works with convolutions of $K$ instead of nondecreasing approximating series?
\item
Kahane's theory is a theory which ensures equality in distribution. Can one build a theory that ensures almost sure equality?  
\end{itemize}
The first point raised above is still an open question. The two following points have been addressed recently as we will now describe in the next two subsections.

\subsection{Generalized Gaussian multiplicative chaos}
Kahane's theory of Gaussian multiplicative chaos relies on the notion of kernels of $\sigma$-positive type. However, on the one hand it is not always straightforward to check such a criterion and on the other hand it seems natural to think that the theory should remain valid for kernels of positive type. A way of getting rid of $\sigma$-positive typeness has been developed in \cite{cf:RoVa1,cf:RoVa}. The idea is to make a convolution product of the covariance kernel with a sequence of continuous functions approximating the Dirac delta function in order to smooth down the singularity of the kernel.

More precisely, consider  a positive definite function $K$ in $\R^d$  (or a bounded domain of $\R^d$) such that 
\begin{equation}\label{eq:defpos}
K(x)=\ln^+\frac{T}{|x|}+g(x)
\end{equation}
 and $g(x)$ is a bounded continuous function.  
Let $\theta: \R^d \rightarrow \R$  be some continuous function with the following properties: 
\begin{enumerate}
\item
$\theta$ is positive definite,
\item
$\theta$ has compact support and is $\alpha$-H\"older for some $\alpha>0$,
\item
 $\int_{\R^d}\theta(x)dx=1$.  
\end{enumerate}

Here is the main theorem of \cite{cf:RoVa}:
\begin{theorem}{\bf [Robert, Vargas, 2008]}\label{th:chaos} 
Let $\gamma^2<2d$. For all $\epsilon >0$, we consider the centered Gaussian field $(X_{\epsilon}(x))_{x \in \R^d}$ defined by the convolution:
\begin{equation*}
E[X_{\epsilon}(x)X_{\epsilon}(y)]=(\theta^{\epsilon} \ast K)(y-x),
\end{equation*}
where $\theta^{\epsilon}=\frac{1}{\epsilon^d}\theta(\frac{.}{\epsilon})$.
Then the associated random measure 
$$\forall A\in \mathcal{B}(\R^d),\quad M_{\epsilon,\gamma}(A)=\int_A e^{\gamma X_{\epsilon}(x)-\frac{\gamma^2}{2}E[X_{\epsilon}(x)^2]}dx$$ converges in law in the space of Radon measures (equipped with the topology of weak convergence) as $\epsilon$ goes to $0$ towards a random measure $M_\gamma$, independent of the choice of the regularizing function $\theta$ with the properties 1., 2., 3. above. 
\end{theorem}

This theorem is useful to define a Gaussian multiplicative chaos associated to the kernel  $K(x)=\ln^+\frac{T}{|x|}$ in dimension $3$ and hence give a rigorous meaning to the Kolmogorov-Obukhov model. Indeed, remind that this kernel is of $\sigma$-positive type in dimension $1$ and $2$. In dimension $d=3$, the function $\ln^+\frac{T}{|x|}$ is positive definite (see \cite{cf:RoVa}) but it is an open question whether it is of $\sigma$-positive type. Therefore, we are bound to apply Theorem \ref{th:chaos} instead of Kahane's theory to define the associated chaos in dimension $3$. In dimension greater than $4$, the kernel $\ln^+\frac{T}{|x|}$ is no more positive definite. Let us also mention that it is proved in \cite{cf:RoVa} that in dimension $1,2,3$ the random variable $M_\gamma(B(0,r))$ (for $r<T$) possesses a $C^{\infty}$ density with respect to the Lebesgue measure.

\begin{rem}\label{remth:chaos}
Starting with a Gaussian field $(X(x))_{x \in \R^d}$ whose covariance is given by \eqref{eq:defpos}, one could also state theorem \ref{th:chaos} in an equivalent way in terms of the fields $X_{\epsilon}=\theta^{\epsilon} \ast X $ where $\theta$ is a function of average $1$ (but not necessarily of positive type) and which satisfies a decreasing condition at infinity.  
\end{rem}

 \subsubsection{Extension to open domains and non-stationary smooth fields} 
In what follows, we explain why the generalized Gaussian multiplicative chaos theory straightforwardly extends to the situation of bounded domains and possibly non stationary covariance kernels. Let $D$ be an open set of $\R^d$. For all $\delta>0$, we set:
\begin{equation*}
D^{(\delta)}=\{x\in D;\,\,{\rm dist}(x ,\partial D)\geq\delta\}
\end{equation*} 
By convention, if $D=\R^d$, we set $D^{(\delta)}$ equal to the Euclidean ball of center $0$ and radius $\frac{1}{\delta}$. We consider a positive definite kernel $K$ (non necessarily stationary) satisfying:
\begin{equation}\label{eq:defpos2}
\forall x,y\in D,\quad K(x,y)=\ln^+\frac{T}{|x-y|}+g(x,y)
\end{equation}
where $\gamma^2 <2d$ and $g$ is a bounded continuous function over $D^{(\delta)}$ for all $\delta>0$. Let $X$ be a random centered Gaussian distribution with covariance given by \eqref{eq:defpos2}. We introduce the following notion of smooth approximation which will play a central role in the rest of the review:

\begin{definition}{\textbf{Smooth Gaussian approximations. }}\label{def:smooth}We say a that a sequence of centered Gaussian fields $(X_\epsilon )_{\epsilon>0}$ is a smooth Gaussian approximation of $K$ if:
\begin{itemize}
\item
for all $x,y \in D$, $\E[  X_\epsilon(x) X_\epsilon(y)]$ converges to $K(x,y)$ as $\epsilon$ goes to $0$. 
\item
 for all $\delta>0$, there exists some constant $C>0$ and $\alpha>0$ such that for all $\epsilon>0$:
\begin{equation*}
\forall x,y \in D^{(\delta)}, \;  \E[  (X_\epsilon(x)-X_\epsilon(y))^2] \leq C  |x-y|^\alpha\epsilon^{-\alpha}. 
\end{equation*}
\end{itemize}
 
\end{definition}

\begin{rem}
In the above definition, the second point is a technical assumption. By standard results on Gaussian processes (see \cite{cf:LeTa} for example), this assumption implies the following useful property: for all $\delta>0$ and $A>0$, there exists $C>0$ such that for all  $x \in D^{(\delta)}$:
\begin{equation*} 
\Pb(  \sup_{y ; \: |y-x| \leq A \epsilon} | X_\epsilon(y)-X_\epsilon(x) |  \geq t )  \leq C e^{-C \frac{t^2}{2}}, \;  \; t \geq 0
\end{equation*}
It is in fact this property implied by the second point which is crucial in the proofs of the results where smooth Gaussian approximations appear.  
\end{rem}

For example, the convolution constructions of the previous section are smooth Gaussian approximations. It is not very difficult to see that convolutions of $X$ with circles, balls or smooth bounded domains are smooth Gaussian approximations of $K$. In fact, the techniques of \cite{cf:RoVa} can be straightforwardly adapted to give:
\begin{theorem}\label{th:refin1}
Assume that we are given two smooth Gaussian approximations $(X_\epsilon )_{\epsilon>0}$ and $(\bar{X}_\epsilon )_{\epsilon>0}$ of $K$ such that:  

\vspace{1mm}

\noindent 1) for some $\gamma^2<2d$, the random measure $M_{\epsilon,\gamma}(dx):=e^{\gamma X_\epsilon(x)- \frac{\gamma^2}{2}E[ X_\epsilon(x)^2  ]}dx$ converges almost surely (possibly along some subsequence) to some random Radon measure $M_\gamma$ in the sense of weak convergence of measures.

\vspace{1mm}

\noindent  2) for all $\delta>0$ and $A>0$, 
\begin{equation*}
 \underset{\underset{|x-y| \leq A \epsilon}{x,y \in D^{(\delta)},}}{\sup}  \Big|  \E[  X_\epsilon(x) X_\epsilon(y)]- \E[  \bar{X}_\epsilon(x) \bar{X}_\epsilon(y)] \Big| \leq \bar{C}_A.
\end{equation*}
where $\bar{C}_A>0$ is some constant independent from $\epsilon$.

\vspace{1mm}

\noindent  3) for all $\delta>0$, 
\begin{equation*}
C_A = \underset{\epsilon \to 0} {\overline{\lim}} \underset{\underset{|x-y| \geq A \epsilon}{x,y \in D^{(\delta)},}}{\sup}\Big|  \E[  X_\epsilon(x) X_\epsilon(y)]- \E[  \bar{X}_\epsilon(x) \bar{X}_\epsilon(y)] \Big|
\end{equation*}
goes to $0$ as $A$ goes to infinity.

Under the above assumptions,  the  random measure $$\bar{M}_{\epsilon,\gamma}(dx)= e^{\gamma \bar{X}_{\epsilon}(x)-\frac{\gamma^2}{2}E[\bar{X}_{\epsilon}(x)^2]}dx$$ converges in law as $\epsilon$ goes to $0$ in the space of Radon measures on $D$ (equipped with the topology of weak convergence) towards the random measure $M_\gamma$.
 \end{theorem}

Let us mention here a simple and straightforward consequence of this theorem. Let $X$ be a random centered Gaussian distribution with covariance given by \eqref{eq:defpos2}. We  consider a function $\theta$ of class $C^1$ such that
\begin{equation}\label{int.theta}
\int_{\R^d}\theta(x)dx=1
\end{equation}
and compactly supported in the ball $B(0,1)$. Then we define the  $\epsilon$-approximation of $X$ by:
$$\bar{X}_\epsilon(x)=\int X(x)\ast \theta^{\epsilon} (x-y)\,dy.$$
Of course, the field  $\bar{X}_\epsilon$ only makes sense for $x\in D^{(\epsilon)}$.  Then the law of the limiting measure
$$M_\gamma=\lim_{\epsilon\to 0}e^{\gamma X_\epsilon(x)- \frac{1}{2}E[ X_\epsilon(x)^2  ]}dx$$ does not depend on the choice of the regularizing function $\theta$.

We will see in subsection \ref{LQG} another application of this theorem in the case of the so-called Liouville measure.

\subsection{Duplantier-Sheffield's approach in dimension $2$}\label{dupshe}
After the two aforementioned theories (Gaussian multiplicative chaos \cite{cf:Kah} and its generalized version \cite{cf:RoVa1,cf:RoVa}), the authors of \cite{cf:DuSh} came up with another contribution to Gaussian multiplicative chaos theory in the special case when the Gaussian distribution $X$ in \eqref{eq:deformelle} is the GFF in a bounded domain and $\sigma$ the Lebesgue measure. The situation can be roughly summarized as follows. Assume for instance that you are given a two-dimensional GFF $X$ on a bounded domain $D$ and you want to define the approximations $(X_n)_n$ in \eqref{eq:deformellen} almost surely as measurable functions of the whole distribution $X$. For instance, you may define $X_n$ as the projections of $X$ along the first $n$ vectors of an orthonormal basis of the Sobolev space $H_0^1(D)$  
and obtain a first way of defining almost surely the associated multiplicative chaos,  called Liouville measure,  as a measurable function of the whole GFF distribution. This falls under the scope of Kahane's theory since you are adding independent Gaussian fields. To be precise, the vectors of the basis are not necessarily nonnegative so that the construction of the measure works but you do not get the uniqueness property of Theorem \ref{th:uniqueness} (except in some special cases: for instance, the Haar basis is composed of nonnegative functions and therefore in this case you are strictly within Kahane's framework).

A second way of defining cutoff approximations of $X$ that are  measurable functions of the whole GFF distribution is to use convolution techniques: you may also define $X_\epsilon(x)$ as the average value of $X$ over the circle of radius $\epsilon$ centered at $x$ (or any other regularizing function) and plug this quantity in place of $X_n$ in \eqref{eq:deformellen}. The approximating measures
$$M_{\epsilon,\gamma}(dx)= e^{\gamma X_\epsilon(x)- \frac{\gamma^2}{2}E[ X_\epsilon(x)^2  ]}dx$$
 are then measurable functions of the  GFF. Though convergence in law follows from the techniques of the generalized Gaussian multiplicative chaos theory \cite{cf:RoVa1,cf:RoVa},   it is important to focus on {\it almost sure} convergence of these measures.
 
\begin{theorem}{\bf [Duplantier, Sheffield, 2008]. }  
 The approximating measures $(M_{\epsilon,\gamma})_\epsilon$ almost surely converge as $\epsilon\to 0$ for the topology of weak convergence of measures along suitable deterministic subsequences.
\end{theorem} 

Another important point is the following. We have now at our disposal two ways to produce the Liouville measure as the almost sure limit of suitable approximations: $H_0^1(D)$ expansions or circle average. In this formulation, it is implicitly assumed that these two approximations yield the same limiting object. Is it really the case? Let us first stress that the  two above constructions have the same law, and they also have the same law as the Liouville measure based on a white noise decomposition of the GFF introduced in \cite{Rnew10} (see Theorem \ref{uniqshef} below). There is another important contribution in \cite{cf:DuSh} to the theory of Gaussian multiplicative chaos

\begin{theorem}{\bf [Duplantier, Sheffield, 2008]. }  
The Liouville measures constructed via  an expansion along an orthonormal basis of $H_0^1(D)$ or circle average approximations are almost surely the same measures.
\end{theorem}

Though carried out in the case of GFF, this uniqueness result makes sense in a more general context: to which extent can we prove that Gaussian multiplicative chaos constructed with different approximations almost surely defined as measurable functions of the whole Gaussian distribution coincide almost surely?

\section{Multifractal analysis of the measures}
The purpose of this section is to give a brief insight into multifractal analysis. More precisely, Kahane proved that the carrier of a Gaussian multiplicative chaos has Hausdorff dimension greater or equal to $d-\frac{\gamma^2}{2}$. It is natural to wonder whether further pieces of information can be given about the structure of this carrier. We will relate this question to the Parisi-Frisch formalism \cite{ParisiFrisch}. Consider a smooth cutoff approximation $(X_\epsilon)_{\epsilon}$ of the field $X$ with covariance $K$ given by \eqref{eq:defpos2} on a domain $D$ and such that $X_{\epsilon'}-X_\epsilon$ is independent from $\sigma \lbrace  X_u; \: u \geq \epsilon \rbrace$ for all $\epsilon'<\epsilon$. We suppose that there exists a constant $C>0$ such that $ \ln \frac{1}{\epsilon}-C  \leq  \E[  X_\epsilon(x)^2  ]  \leq  \ln \frac{1}{\epsilon}+C$ for all $x \in D$. In this context, Kahane introduced the set of points
$$\big\{x\in D;\,\,\, \lim_{\epsilon}\frac{X_\epsilon (x)}{-\ln \epsilon}=\gamma   \big\}$$
and showed that it gives full measure to the chaos $M_\gamma$, thus showing that it has Hausdorff dimension greater or equal to $d-\frac{\gamma^2}{2}$. This set of points has been called {\it thick points} of the GFF in  \cite{peres} when the field $X$ is the GFF and $X_\epsilon$ corresponds to circle averages (hence, we are not exactly in the framework of Kahane as circle averages of the GFF do not correspond to adding independent fields; nonetheless, the two frameworks are very similar). We stick to this terminology as we feel that it is well-sounding. The Hausdorff dimension of thick points is derived in \cite{peres} in the context of circle averages of the GFF. Let us further mention that exact equality of Hausdorff  dimension is well known in the closely related  context of Mandelbrot's multiplicative cascades, see \cite{Barral1} and references therein. In the next subsection, we will generalize the results in \cite{peres, cf:Kah} to a large class of log-correlated Gaussian fields in all dimensions. After we will discuss general multifractal analysis. 
\subsection{The Peyri\`ere probability measure and the thick points}\label{sub:thick}

Let us first set the framework of this subsection. Consider a bounded domain $D\subset \R^d$ and a positive definite kernel $K$ (non necessarily stationary) satisfying:
\begin{equation}\label{eq:defposns}
\forall x,y\in D,\quad K(x,y)=\ln^+\frac{T}{|x-y|}+g(x,y)
\end{equation}
where $g$ is a bounded continuous function over $D$.  We work on a fixed probability space $(\Omega,\mathcal{F},\Pb)$. On this space, we consider a centered Gaussian field $(X_\epsilon(x))_{\epsilon >0, x \in D}$ such that: 
\begin{enumerate}
\item
$(\epsilon,x) \rightarrow X_\epsilon(x)$ is almost surely continuous on $]0,\infty[ \times D$,
\item
$(X_\epsilon )_{\epsilon>0}$ is a smooth Gaussian approximation of $K$ (in the sense of Definition \ref{def:smooth}),
\item
for each fixed $x \in D$, $\epsilon \rightarrow X_\epsilon(x)$ has independent increments, 
\item
There exists a constant $C>0$ such that for all $\eta<\epsilon$:
\begin{equation}\label{am.hyp}
 \ln_+\frac{1}{ |x-y|+\epsilon}-C  \leq  \E[  X_\epsilon(x) X_\eta(y)  ]  \leq  \ln_+\frac{1}{ |x-y|+\epsilon}+C, \; \; x,y \in D.
\end{equation}
\item
for all $\gamma^2<2d$, there exists a random Radon measure $M_\gamma$ such that $$M_{\epsilon,\gamma}(dx)=e^{\gamma X_\epsilon(x)- \frac{1}{2}E[ X_\epsilon(x)^2  ]}dx$$ converges almost surely (possibly along some subsequence) to $M_\gamma$ on $D$.
\end{enumerate}
For $\gamma^2<2d$ and $q \in ]0, \frac{\sqrt{2d}}{\gamma}[$, we set:
 \begin{equation*}
 K_{\gamma,q}= \left\{ x \in D; \: \underset{\epsilon \to 0}{\lim} \: \frac{\ln M_{\gamma}(B(x,\epsilon))}{\ln \epsilon} = d+(\frac{1}{2}-q)\frac{\gamma^2}{2} \right\} \cap \left\{ x \in D; \: \underset{\epsilon \to 0}{\lim} \: \frac{ X_\epsilon(x)}{-\ln \epsilon} = \gamma q \right\}.
 \end{equation*}

We can now state the following theorems   (see appendix): 

\begin{theorem}\label{analysemulti}
For all $\gamma^2<2d$, almost surely, the set $K_{\gamma,q}$ gives full mass to the measure $M_{q \gamma}$, i.e. $M_{q \gamma}(^c K_{\gamma,q})=0$.
\end{theorem}

\begin{theorem}\label{cor:analysemulti}
For all $\gamma^2<2d$, almost surely, the set  $K_{\gamma,1}$ has Hausdorff dimension $d-\frac{\gamma^2}{2}$.
\end{theorem}

A quite simple and concise proof of the above theorem is gathered in the appendix and, to our knowledge, is new in such generality (though much more is known in dimension 1: see \cite{Barral1,BS}). In dimension $2$ and in the context of circle average approximations of the GFF, another paper \cite{peres} focused on the thick points of the GFF, improving Kahane's result by giving the upper bound for the Hausdorff dimension. Let us also stress here that the place of the "almost sure" is important in Theorem \ref{analysemulti} and Corollary \ref{cor:analysemulti}. In the case of multiplicative cascades, a stronger statement is proved in \cite{barralJTP}, where the "almost sure" is valid simultaneously for all $\gamma^2<2d$.

In case the reader wishes to skip the whole proof of Theorem \ref{analysemulti}, we sketch here  Kahane's argument about the lower bound for the Hausdorff dimension of the thick points of $X$.  We focus on the case where $q=1$ and $D$ has Lebesgue measure $1$. The key point is to introduce the so-called Peyri\`ere probability measure (also called {\it rooted measure} in \cite{cf:DuSh,Rnew7,Rnew12}):
\begin{equation}\label{defPeyriere}
Q(F(\omega,x)) := \E[ \int_D F(\omega,x) M_\gamma (dx)    ] .
\end{equation}
This measure was introduced by Peyri\`ere in the context of discrete cascades (see \cite{KP}) and was used by Kahane in his work on Gaussian multiplicative chaos. It is obvious that under this measure the process $(x,t) \rightarrow X_{e^{-t}}(x)$ is a Brownian motion with drift $\gamma$, hence we have $M_{\gamma}(^c \tilde{K}_{\gamma}))=0$ where:
\begin{equation}\label{defthick}
 \tilde{K}_{\gamma}= \left\{ x \in D; \: \underset{\epsilon \to 0}{\lim} \: \frac{ X_\epsilon(x)}{-\ln \epsilon} = \gamma  \right\}.
\end{equation}
In particular, it is straightforward from Theorem \ref{th:carrier} that $ \tilde{K}_{\gamma}$ has Hausdorff dimension at least $d-\frac{\gamma^2}{2}$. Let us further stress that similar ideas are used in \cite{peres}.

\subsection{General multifractal formalism: a heuristic introduction}
In this section, we discuss a bit of general multifractal formalism for the measures $M_\gamma$. \textbf{This discussion is essentialy based on heuristics though it should not be very difficult to make it rigorous mathematically}. Let us mention that the study of multifractal formalism represents a very wide domain of mathematics and physics (in particular turbulence); therefore, being exhaustive in this field is way beyond the scope of this review where we  only mention the case of measures of the form $e^{X(x)} dx$ for some log correlated gaussian field $X$. 

Multifractal analysis is based on the calculation of the $L^q$-spectrum of the measure $M_\gamma$, defined as  
$$
q\in\mathbb{R}\mapsto \tau_{M}(q)=\liminf_{r\to 0^+}\frac{\log \sup\Big \{\sum_i  M_\gamma(B(x_i,r))^q\Big \}}{\log(r)},
$$
where the supremum is taken over all the centered packing of $[0,1]^d$ by closed balls of radius  $r$. This study is achieved in \cite{Barral1,BS} in dimension $1$ and, with a bit of effort, it should not be difficult to generalize this result to higher dimensions: 
$$
\tau_{  M_\gamma}(q)= 
\begin{cases}
 (\sqrt{d}+\frac{\gamma}{\sqrt{2}})^2 q&\text{if } q\le -\frac{\sqrt{2d}}{\gamma},\\
  \xi(q)-d&\text{if }  q\in [-\frac{\sqrt{2d}}{\gamma},\frac{\sqrt{2d}}{\gamma}],\\
(\sqrt{d}-\frac{\gamma}{\sqrt{2}})^2 q &\text{if }q\ge \frac{\sqrt{2d}}{\gamma}.
\end{cases}
$$
The so-called {\it multifractal formalism} holds for $ M_\gamma$: we can relate the $L^q$-spectrum of $M_\gamma$ to the local regularity of $M_\gamma$. More precisely, define $$E_\delta=\big\{x\in [0,1]^d;\liminf_{r\to 0^+}\frac{\ln M_\gamma(B(x,r))}{\ln (r)}=\delta\big\}
\quad (\delta\ge 0).$$
The singularity spectrum of $ M_\gamma$, i.e. the mapping $\delta\ge 0\mapsto \dim E_\delta$ ($\dim$ meaning Hausdorff dimension), is given by the celebrated Parisi-Frisch formula (see \cite{ParisiFrisch}):  
\begin{equation}\label{eq:ParFri}
\delta\ge 0\mapsto \tau_{M_\gamma}^*(\delta)=\inf\{\delta q-\tau_{ M_\gamma}(q):q\in \R\} \wedge 0.
\end{equation}
Note that the above formula is a Legendre transform since the formula comes in fact from a large deviation argument (which can be made rigorous in terms of box counting dimensions by using the Gartner-Ellis theorem). Multifractal analysis is essentially focused on studying the valididity of (\ref{eq:ParFri}) in the broad context of all (random or deterministic) measures. From this and the explicit expression on $\tau_{M_\gamma}$, one can deduce the following dimension result:
\begin{equation}\label{quadratic}   
\dim E_\delta= \begin{cases}
d-\frac{1}{2}(\frac{d}{\gamma}+\frac{\gamma}{2} -\frac{\delta}{\gamma})^2, &\text{if }   \delta \in [   (\sqrt{d}-\frac{\gamma}{\sqrt{2}})^2, (\sqrt{d}+\frac{\gamma}{\sqrt{2}})^2  ].\\
  0&\text{elsewhere}.\\
\end{cases}
\end{equation}
In essence, multifractal formalism and computing the dimension of thick points is the same thing (and hence we get the same associated dimensions) if we admit the following commonly used heuristic:
\begin{equation}\label{heuristic}
M_\gamma(B(x,r)) \underset{r \to 0}{\sim} C(x) r^d e^{\gamma X_{r}(x)-\frac{\gamma^2}{2} \ln \frac{1}{r}}
\end{equation}
where $C(x)$ is some random constant of order $1$ (nearly independent of $r$ but very dependent on $x$). Note that the notion of star scale invariance is nothing but a rigorous formulation of the heuristic (\ref{heuristic}). Swithching from multifractal formalism to thick point formalism rigorously implies handling the $C(x)$ term which is of order $1$ but fluctuates wildely as a function of $x$. As an example of application of heuristic (\ref{heuristic}), let us recover formula (\ref{quadratic}) from corollary \ref{cor:analysemulti}. By using (\ref{heuristic}), we get the following equivalence:
\begin{equation*}
M_\gamma(B(x,r)) \underset{r \to 0}{\sim} r^\delta  \Leftrightarrow X_{r}(x) \underset{r \to 0}{\sim}  (\frac{d}{\gamma}+\frac{\gamma}{2} -\frac{\delta}{\gamma})   \ln \frac{1}{r},
\end{equation*}
hence leading to formula (\ref{quadratic}) thanks to corollary \ref{cor:analysemulti}.

\section{Applications of Gaussian multiplicative chaos}
In this section, we review some applications in direct relation with Kahane's theory: some of them are well known, some of them are new. Some further applications, rather related to recent generalizations of the theory, will be given in Section \ref{sec:gener}.

\subsection{Volatility of a financial asset or  boundary Liouville measure}\label{sectiondim1}
The main application of the theory is to give a meaning to the "limit-lognormal" model introduced by Mandelbrot in \cite{cf:Man}. The "limit-lognormal" model corresponds to the choice of a stationary kernel $K$ on $\R\times \R$ given by:
\begin{equation}\label{eq:modeleog}
K(x,y)= \ln^{+}\frac{T}{|x-y|}+g(x,y) 
\end{equation}
where $T$ is a positive parameters, $g$ is a bounded continuous function and $\sigma$ is chosen to be the Lebesgue measure on $\R$. This model has many applications: part of them are discussed below.

If $(X(t))_{t \geq 0}$ is the logarithm of the price of a financial asset, the volatility $M_\gamma$ of the asset on the interval $[0,t]$ is by definition equal to the quadratic variation of $X$:
\begin{equation*}
M_\gamma[0,t]=\lim_{n \to \infty}\sum_{k=1}^{n}(X(tk/n)-X(t(k-1)/n))^2
\end{equation*}
The volatility $M$ can be viewed as a random measure on $\R$. The choice for $M_\gamma$ of multiplicative chaos associated to the kernel $K(x,y)=\ln^{+}\frac{T}{|x-y|}$ satisfies many empirical properties measured on financial markets: lognormality of the volatility, long range correlations (see \cite{cf:Cizeau} for a study of the SP500 index and components and \cite{cf:Co} for a general review). With this kernel, the measure is called the lognormal multifractal random measure (MRM) and is a particular case of the log-infinitely divisble multifractal random measures (see \cite{Bar,cf:Sch} for the log-poisson case and \cite{cf:BaMu} for the general case). 
Note that $K$ is indeed of $\sigma$-positive type   so $M_\gamma$ is well defined. In the context of finance, $\gamma^2$ is called the intermittency parameter in analogy with turbulence and $T$ is the correlation length. Volatility modeling and forecasting is an important field of finance since it is related to option pricing and risk forecasting: we refer to \cite{cf:DuRoVa} for the problem of forecasting volatility with this choice of $M_\gamma$.   

Given the volatility $M_\gamma$, the most natural way to construct a model for the (log) price $X$ is to set:
\begin{equation}\label{eq:MRW}
X(t)=B_{M_\gamma[0,t]}                      
\end{equation}
where $(B_{t})_{t \geq 0}$ is a Brownian motion independent of $M_\gamma$. Formula (\ref{eq:MRW}) defines the Multifractal Random Walk (MRW) first introduced in \cite{cf:BaDeMu} (see \cite{cf:BaKoMu} for a recent review of  financial applications of the MRW model).

\begin{figure}[h]
\centering
\subfloat[SP500 returns 2001-2009]{\includegraphics[width=0.47\linewidth]{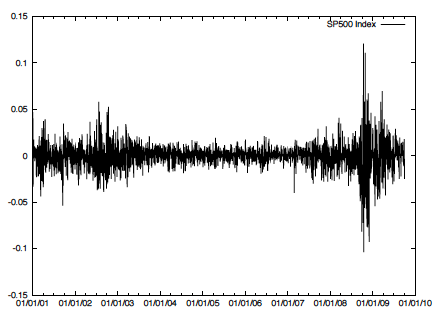}}
\,\,\subfloat[Returns simulated with Black-Scholes]{\includegraphics[width=0.47\linewidth]{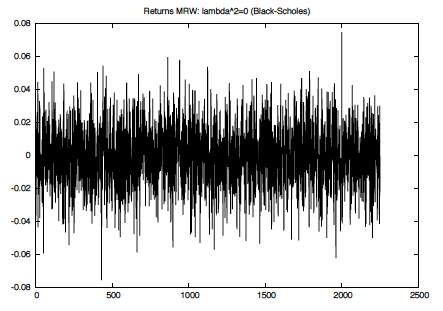}}\\
\subfloat[Returns simulated with MRW]{\includegraphics[width=0.47\linewidth]{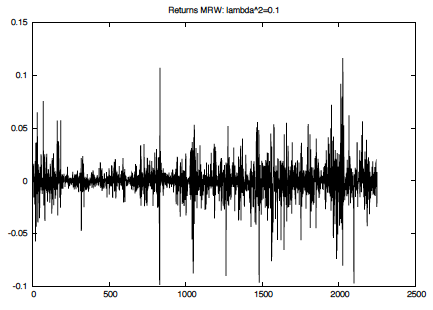}}
\,\,\subfloat[Returns simulated with MRW]{\includegraphics[width=0.47\linewidth]{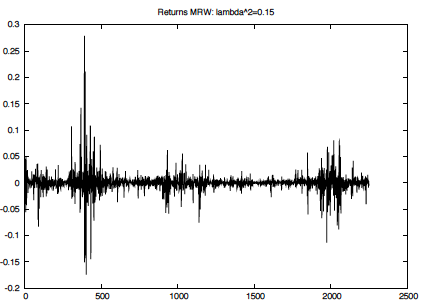}}
\caption{Intermittency in financial markets}
\label{interm}
\end{figure}
In Figure \ref{interm}, observe the burst of activity (intermittency) for the SP500 index or the MRW model. This reflects the Parisi-Frisch formalism (or multifractal formalism): the (strong) variations of regularity are related to the non-linearity of the power law spectrum.

\subsubsection{MRM with infinite correlation length or boundary Liouville measure}
Motivated by financial applications where the correlation length $T$ is too big to be measured on markets, the authors of \cite{cf:DuRoVa} adressed the issue of forecasting volatility $M_\gamma$ in the limit $T \to \infty$. More precisely, they start by forecasting the log volatility, the $1/f$ noise which lives in the quotient space of distributions defined up to some additive constant. In this context, taking the exponential, the associated random measure $M_\gamma$ is then defined up to a multiplicative constant and is the boundary Liouville measure in the upper half plane considered in \cite{cf:DuSh}.

\subsection{Liouville Quantum Gravity and KPZ}\label{LQG}
Let us first roughly explain the original motivations coming from the physics literature.
We want to define a random distribution $(g,X)$, the partition function of which formally writes
\begin{equation}\label{partition}
Z=\int \mathcal{D}g\mathcal{D}X\,e^{-S_M(X,g)-\lambda V_g(\Sigma)}
\end{equation}
where $S_M$ is some conformally invariant action for matter fields coupled to a compact simply connected two dimensional surface $\Sigma$ with metric $g$, $\lambda$ is a constant (we do not discuss its value), $V_g$ is the volume form  of $g$ and $X$ is an embedding from $\Sigma$ into a $c$-dimensional spacetime. Here we adopt standard path integral notations: the above integral just means that we sum over all possible embeddings and metrics.

\begin{example} For the free bosonic string, we consider the Polyakov action $$S_M=\frac{1}{8\pi}\int \partial_gX \cdot\partial_g X\,dV_g$$ where $X$ specifies the embedding of $\Sigma$ into flat $D$-dimensional space-time.
\end{example}
\begin{example} For the massive Ising model, we consider  $$S_M=\int \bar{X}\big(\partial_gX +m\bar{X}X \big)\,dV_g.$$ 
\end{example}
 
Notice that every metric $g$ on $\Sigma$ can be decomposed as $f^*g=e^{\varphi}g_0$, where $g_0$ is a fixed metric on $\Sigma$, $f$ is a $g_0$-diffeomorphism and $f^*g$ is the pullback metric of $g$ along $f$. So we may perform the path integral by gauge-fixing: we choose a gauge defining an equivalence class over all the metrics and perform the above integral over a slice that cuts through once each gauge equivalence class. In view of the above factorization property, a natural choice of the gauge is the conformal gauge. We choose a family $(\hat{g})$ of representatives of each equivalence class of conformally equivalent metrics and we perform   a sum over $(\hat{g})$ and over the equivalence class of $\hat{g}$ for all $\hat{g}$. The Jacobian of such a "change of variables" is  the so-called  Faddeev-Popov determinant $\triangle_{FP}(\hat{g})$.  We do not detail this here but the reader is referred to \cite{Pol,polch,cf:Da} for further details and references. Let us just say that once this determinant has been computed,  it remains to make sure that  the quantity resulting from these computations does not depend on the choice of the family of representatives $(\hat{g})$: in physics language, we have to compute the Weyl anomaly.  Performing these computations lead to considering the Liouville action  ($R_g$ is the Ricci tensor of the metric $g$)
$$S_L(\varphi, g )=\frac{1}{2\pi \gamma^2}\int (  \partial_{g}\varphi\cdot\partial_{g}\varphi+QR_{g}\varphi+\mu e^{\varphi})dV_{g}$$
and the matter action $S_M(X,\hat{g})$
in such a way that
$$\int \mathcal{D}X \mathcal{D}g\, e^{-S(X,g)}= \int  DX  Df D\varphi \,\triangle_{FP}(\hat{g}) \,e^{-S_M(X,\hat{g})-S_L(\varphi,\hat{g} )} .$$
For $c\in ]-\infty,1]$, the value of $\gamma\in [0,2]$ is related to the central charge of the matter by
$$\gamma= \frac{\sqrt{25-c}-\sqrt{1-c}}{\sqrt{6}}.$$
When the cosmological constant $\mu$ is set to $0$, the Liouville action reduces to that  of a free massless boson. Mathematically speaking the corresponding field $\varphi$ is a Gaussian Free Field. In critical Liouville quantum gravity, we are therefore led to considering random metrics of the form $e^{\gamma \varphi}\hat{g}$ and an area measure $e^{\gamma \varphi}dV_{\hat{g}}$, where $\varphi$ is a Free Field in the background metric $\hat{g}$. Therefore, the world sheet $\Sigma$ may be equipped with two metrics: the background metric $\hat{g}$ 
and the quantum metric $e^{\gamma \varphi}\hat{g}$. 

Furthermore, conditionally on a fixed  background metric $\hat{g}$ (just discarding  $\hat{g}$  from the randomness), the metric $e^{\gamma \varphi}\hat{g}$ and the matter field are independent as may be seen from the resulting partition function. Knizhnik, Polyakov and Zamolodchikov have derived in \cite{cf:KPZ} a relation between the scaling exponents of the background metric and the quantum metric 
$e^{\gamma \varphi}\hat{g}$, the so-called KPZ formula. This is a very rough description  of Liouville quantum gravity  and the reader may consult \cite{Pol,cf:KPZ,polch,cf:Da} for further details. 

Let us now explain why this KPZ formula may be of interest in the study of models of statistical physics at their critical point. Physicists understood a long time ago (see \cite{amb,amb2,cf:KPZ,cf:Da,DistKa} and certainly many others) that this continuum model of quantum gravity admits a discretized counterpart via random triangulations (or other $p$-angulations) of surfaces. The prototype of such $p$-angulations is the Brownian map studied in (see  \cite{LeGallS,legall,LGM,miermont}) and corresponds to the pure gravity case $c=0$. But we may also couple a model of statistical physics (for instance random walks, percolation, Ising model, Potts model,...) to discrete quantum gravity, i.e. by considering a model of statistical physics on the $p$-angulation in such a way that, as in the continuum case, the partition function involves both the $p$-angulation and that of the model of statistical physics. The point is that, at their critical point, these models in two dimensions should behave as a conformal field theory and may be thought of as the matter field described above. By taking the limit as the discretization step goes to $0$, these models of discrete quantum gravity should converge towards Liouville quantum gravity. Interestingly, the independence of the Liouville field and the matter field (conditionally on the background metric) suggests that the same phenomena should occur when taking the limit in the discrete model. Therefore the KPZ formula may be applied to this "discrete matter field", which are roughly independent of the fluctuating quantum metric: it  becomes particularly useful when the scaling exponents of a particular model can be more easily computed in its quantum gravity form than its background one, or vice-versa. For instance, Duplantier \cite{dupRW} have used these techniques to conjecture the exact values of the Brownian intersection exponents, which were finally rigorously derived in \cite{LSW1,LSW2,LSW3} via Schramm-Loewner-Evolution (SLE) techniques. The reader may consult \cite{disc:da1,disc:da2,disc:Froh} for "pure string models" with $c=0$ or $c=-2$, and \cite{boul,dupkos,kaz1,kaz2} for critical systems on random $p$-angulations like $Q$-Potts model, percolation or tree like polymers.

Understanding Liouville quantum gravity from a mathematical rigorous angle is a wide task, which mathematicians have tackled only recently and may take on various aspects, some of them have obviously connections with Gaussian multiplicative chaos theory.  As explained above, the mathematical formulation of the problem of constructing (critical) $2d$-Liouville quantum gravity could be roughly summarized as follows: construct a random metric on a two dimensional Riemannian manifold $D$, say a domain of $\R^2$ (or the sphere) equipped with the Euclidean metric $dz^2$, which takes on the form
\begin{equation}\label{i.metric}
e^{\gamma X(z)}dz^2
\end{equation}
where $X$ is a Gaussian Free Field (or possibly other Free Fields) on the manifold $D$   and $\gamma\in [0,2)$ is a coupling constant. 
\begin{figure}[h]
\centering
\includegraphics[width=0.8\linewidth]{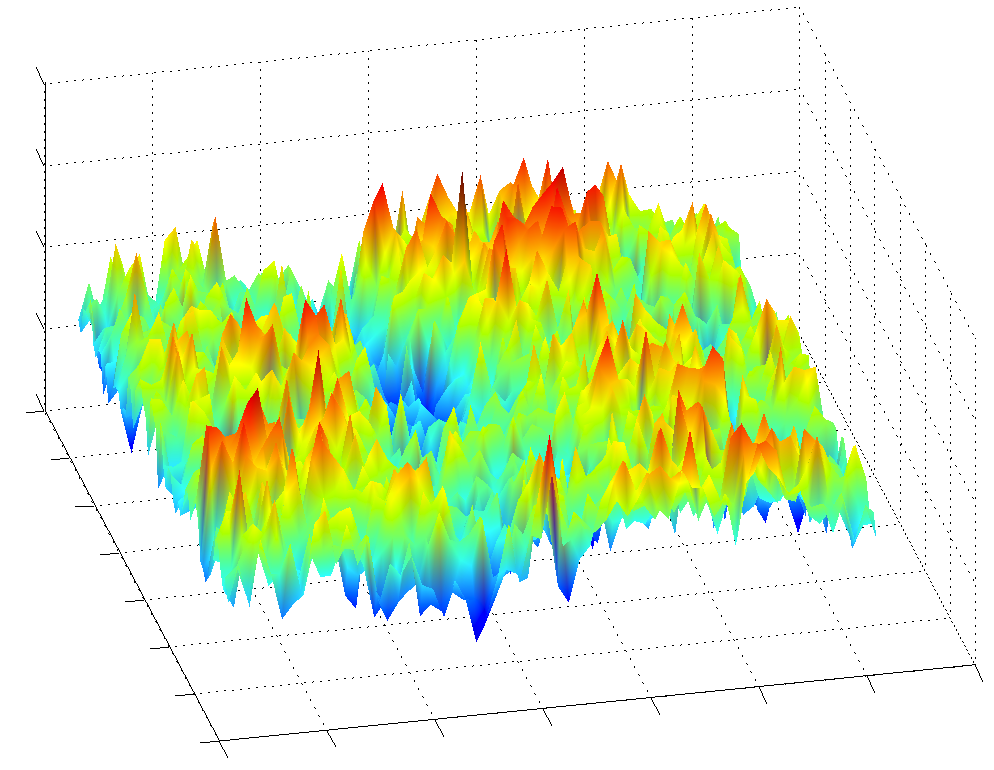}
\caption{Simulation of a GFF on the unit torus.}
\label{fig.GFF}
\end{figure}

The issue of constructing the distance associated to the metric remains unsolved.  Yet recent progress are made in \cite{LBM,LBM2,spectral} concerning the Brownian motion, Laplace-Beltrami operator or heat kernel of $2d$-Liouville quantum gravity. Nevertheless, we focus below on  the volume form, which obviously falls under the scope of Gaussian multiplicative chaos theory. This theory allows us to construct  a random measure of the type: 
\begin{equation}\label{measLQG}
M_\gamma(dx)=e^{\gamma X(x)-\frac{\gamma^2}{2}\E[X(x)^2]}\,dx,
\end{equation}
which will be called {\it Liouville measure}. The points that we will address below are the following. In order to apply Gaussian multiplicative chaos theory and  define the above measure, we have   to choose a cutoff approximation of the GFF  and there are several possible choices, which we will discuss. Furthermore we will explain why these cutoff approximations lead to the same limiting measure $M_\gamma$.

Recall that the GFF over a bounded simply connected domain $D$ with for instance Dirichlet boundary condition is a centered Gaussian distribution with covariance kernel given by the Green function $G$ of the Laplacian, i.e. $\triangle G(x\cdot)=-2\pi\delta_x$, with Dirichlet boundary condition. Actually, other types of boundary conditions may be imposed but it suffices to detail the Dirichlet boundary conditions to draw a clear picture of the techniques involved.

\subsubsection*{Decomposition of the GFF via eigenfunctions of the Laplacian}
Let us consider the eigenfunctions $(e_n)_{n\geq 1}$ of the Laplacian with Dirichlet boundary conditions. They form an orthonormal basis of $L^2(D)$ with negative associated eigenvalues $(\lambda_n)_{n\geq 1}$. A natural choice of decomposition of the GFF is to write (formally):
$$X(x)=\sum_{n\geq 1}Y_n(x)$$ where $Y_n$ is a smooth Gaussian field defined by 
$$Y_n(x)=\frac{\beta_n}{\sqrt{-\lambda_n}}\, e_n(x).$$
Here $(\beta_n)_{n\geq 1}$ is a sequence of i.i.d standard Gaussian random variables. Observe that the sequence $(\beta_n)_{n\geq 1}$ can be chosen to be measurable with respect to the whole GFF distribution: it suffices to choose   $$\beta_n=(-\lambda_n)^{\frac{1}{2}}\int_DX(x) e_n(x)\,dx.$$  
The covariance kernel of $Y_n$ matches
$$k_n(x,y)=(-\lambda_n)^{-1}e_n(x)e_n(y).$$
 It is well known that the eigenfunctions are smooth so that $Y_n$ is a smooth Gaussian field. The important point here is that the approximating sequence 
 $$\forall n\geq 1,\quad X_n(x)=\sum_{k=1}^nY_k(x)$$ is almost surely defined as a function of the whole GFF distribution $X$. Furthermore, the sequence $(Y_n)_n$ are independent Gaussian processes. Though elegant and simple,  this decomposition also possesses drawbacks because the covariance kernel of each $Y_n$ is not nonnegative and, generally speaking, it is hard to get a tractable expression of $k_n$ (or rather their partial sums), excepted maybe in terms of lattice approximations (discrete GFF, see \cite{She07}).

\begin{rem}
Actually, any orthonormal basis of $H^1_0(D)$ produces a decomposition of the GFF function \`a la Kahane, i.e. a sum of independent Gaussian processes with continuous covariance kernels. Another very important decomposition relying on an $H^1_0(D)$-basis is the projection of  the GFF onto the Haar basis. In that case, the corresponding kernels $k_n$ are continuous and positive. 
\end{rem}

\subsubsection*{White noise decomposition of the GFF}
Another possible decomposition of the Green function is based on the formula:
\begin{equation*}
G_D(x,y)=  \pi \int_{0}^{\infty}p_D(t,x,y)dt.
 \end{equation*}
where $p_D$  is the (sub-Markovian) semi-group of a Brownian motion $B$ killed upon touching the boundary of $D$, namely
$$p_D(t,x,y)=P^x(B_{t} \in dy, \; T_D > t)$$
 with $T_D=\text{inf} \{t \geq 0, \; B_t\not \in D \}$.  Note that the $\pi$ term ensures that: 
 \begin{equation*}
 G_D(x,y) \underset{|x-y| \to 0}{\sim}  \ln \frac{1}{|x-y|}. 
 \end{equation*}
 Hence we can write:
\begin{equation}\label{greendecomp}
G(x,y)=\sum_{n\geq 0}K_n(x,y)\quad \text{with }\quad K_n(x,y)=\pi\int_{\frac{1}{2^{n}}}^{\frac{1}{2^{n-1}}}p_D(t,x,y)dt, \; n \geq 1
\end{equation}
and $K_0(x,y)=\pi\int_{1}^{\infty}p_D(t,x,y)dt$.
The continuity of $p_D$ implies that $K_n$ is continuous. The symmetry of $p_D$ implies that $K_n$  is positive definite. Indeed, for each smooth function $\varphi$ with compact support in $D$,  we have for $n \geq 1$:
\begin{align*}
\int_D\int_D\varphi(x)K_n(x,y)\varphi(y)\,dx\,dy&=2\pi\int_D\int_D\int_{\frac{1}{2^{n}}}^{\frac{1}{2^{n-1}}}\varphi(x)p_D(t,x,y)\varphi(y)\,dt\,dx\,dy\\
&=\pi\int_D\int_D\int_D\int_{\frac{1}{2^{n}}}^{\frac{1}{2^{n-1}}}\varphi(x)p_D(t/2,x,z)p_D(t/2,z,y)\varphi(y)\,dt\,dz\,dx\,dy\\
&=\pi\int_{\frac{1}{2^{n}}}^{\frac{1}{2^{n-1}}}\int_D\Big(\int_D\varphi(x)p_D(t/2,x,z)\,dx\Big)^2\,dt\,dz\\
&\geq 0.
\end{align*}
Since $K_n$ is obviously positive, we can apply Kahane's theory of Gaussian multiplicative chaos to define the Liouville  measure \eqref{measLQG}.
We further stress that this argument implies a white noise decomposition of the underlying GFF: the most  direct way to construct a GFF is then to consider a white noise $W$ distributed on $D\times \R_+$ and define
 $$X(x)=\sqrt{\pi}\int_{D\times \R_+}p_D(\frac{s}{2},x,z)W(dz,ds).$$
One can check that $\E[X(x)X(y)]=\pi\int_0^{\infty} p_D(s,x,y)\,ds=G_D(x,x') $. One can even work with a continous parameter $\epsilon$ and define the Liouville measure as the almost sure limit as $\epsilon \to 0$ of $M_\epsilon(dx)=e^{\gamma X_\epsilon(x)-\frac{\gamma^2}{2}\E[X_\epsilon(x)^2]}\,dx$ where the corresponding cut-off approximations $X_\epsilon$ are given by:
$$X_\epsilon(x)=\sqrt{\pi}\int_{D\times [\epsilon^2,\infty[}p_D(\frac{s}{2},x,z)W(dz,ds).$$
 Indeed, within this framework introduced in \cite{Rnew10}, the sequence $(M_\epsilon(A))_{\epsilon >0}$ is a positive martingale for all compact set $A$. Note the following expression for the covariance of $X_\epsilon$:
 \begin{equation*}
   E[ X_\epsilon(x) X_\epsilon(y)   ]=\pi\int_{\epsilon^2}^{\infty} p_D(s,x,y)\,ds
  \end{equation*} 
Once we define the Liouville measure with this white noise construction, it is not hard to see that we fall in fact under the scope of theorem \ref{th:refin1}  (see our theorem \ref{uniqshef} below). In particular, we claim:
\begin{lemma}\label{SmoothGFF}
The sequence $X_\epsilon$ is a smooth Gaussian approximation of $G_D$. 
\end{lemma}
We will sketch a proof of this point in the appendix.

\subsubsection*{Circle average}

As explained in subsection \ref{dupshe}, the authors in \cite{cf:DuSh} have suggested a slightly different approach: instead of using the $\sigma$-positivity of the covariance kernel of the GFF to construct an approximating sequence \eqref{eq:deformellen} that is a martingale, they regularize the GFF along circles to construct their approximating sequence. More precisely, consider a GFF $X$ and define $X_\epsilon(z)$ as the mean value of $X$ along the circle centered at $z$ with radius $\epsilon$, formally understood as:
$$X_\epsilon(x)=\frac{1}{2\pi}\int_0^{2\pi}X(x+\epsilon e^{i \theta})\,d\theta.$$
The covariance kernel is given by
$$G_\epsilon(z,z')=\int \int G_D(x,y)\, \mu_\epsilon^x(du)\mu_\epsilon^y(dv)$$
where $\mu_\epsilon^x(du)$ stands for the uniform probability measure on the circle centered at $x$ with radius $\epsilon$. This expression can be given a rigorous sense \cite{cf:DuSh}. The main advantage of this construction is that it is well fitted to play with   the spatial Markov property of the GFF. Nevertheless, the increments $(X_\epsilon-X_{\epsilon'})_{\epsilon<\epsilon'}$ are not independent, getting trickier the proof of the almost sure convergence of the chaos. On the other hand, this circle average construction falls under the scope of the regularization procedures developed in \cite{cf:RoVa} in order to get convergence and uniqueness in law.

\subsubsection*{Equivalence of the  constructions}

The first question that you must have in mind is: "To which extent do the above cut-off approximations yield the same limiting multiplicative chaos?". We claim:

\begin{theorem}\label{uniqshef}
The law of the limiting chaos does not depend on the cutoff approximations listed above, namely  white noise decomposition, eigenvalues of the Laplacian, $H^1_0(D)$ expansions or circle average (or more generally convolution by $C^1$ functions or averages on smooth domains, like ball-averages...). 
\end{theorem}

Before proving this theorem, let us make some further comments. In dimension $2$, the Lebesgue measure is obviously in the class $R_{2-\epsilon}^+$ for all $\epsilon>0$ and the Green function can be rewritten as 
\begin{equation}\label{Green}
G_D(x,y)=\ln_+\frac{1}{|x-y|}+g(x,y)
\end{equation}
for some bounded continuous function $g$. Therefore, all the Kahane  machinery applies. In particular, Theorem \ref{th:degeneracy} ensures that the Liouville measure is non trivial if and only if  $\gamma^2<4$, whatever the choice of the cut-off approximation.

\vspace{3mm}
\noindent {\it Proof of Theorem \ref{uniqshef}.}

Almost sure equivalence between a given $H^1_0(D)$ expansion and circle average is already proved in \cite{cf:DuSh}, and therefore equivalence in law holds. 

\hspace{0.2 cm}

\noindent
\emph{First proof:}

Therefore, it suffices to prove equivalence in law between the white noise decomposition and the circle average construction to get theorem \ref{uniqshef}.  In view of theorem \ref{th:refin1}, one could for instance establish:
\begin{itemize}
\item
For all $\delta>0$ and $A>0$, 
\begin{equation*}
 \underset{\underset{|x-y| \leq A \epsilon}{x,y \in D^{(\delta)},}}{\sup}  |\int \int G_D(x,y)\, \mu_\epsilon^x(du)\mu_\epsilon^y(dv)  -  \pi\int_{\epsilon^2}^{\infty} p_D(s,x,y)\,ds   | \leq \bar{C}_A.
\end{equation*}
where $\bar{C}_A>0$ is some constant independent from $\epsilon$.
\item
For all $\delta>0$, 
\begin{equation*}
C_A = \underset{\epsilon \to 0} {\overline{\lim}} \underset{\underset{|x-y| \geq A \epsilon}{x,y \in D^{(\delta)},}}{\sup}  |\int \int G_D(x,y)\, \mu_\epsilon^x(du)\mu_\epsilon^y(dv) -\pi\int_{\epsilon^2}^{\infty} p_D(s,x,y)\,ds|
\end{equation*}
goes to $0$ as $A$ goes to infinity.
\end{itemize}
in order to prove equivalence between circle average and white noise decomposition. These two estimates are not very difficult to obtain but we will not detail this point since a second very direct proof is possible.

\hspace{0.2 cm}

\noindent
\emph{Second proof:}
Therefore, it suffices to prove equivalence in law between the white noise decomposition and one $H^1_0(D)$ expansion.
Just notice that the white noise decomposition and the expansion along the Haar basis correspond to two $\sigma$-finite decompositions (\ref{eq:defsigma}) of the Green function. Hence, by theorem \ref{th:uniqueness}, the two constructions are equivalent in law.

\qed

\subsubsection*{KPZ formula: almost sure Hausdorff version}\label{sec:KPZ}
As explained above, the KPZ formula has been introduced in Liouville quantum gravity and can be thought of as a bridge between the values of the scaling exponents computed with the quantum metric $e^{\gamma X(z)}\,dz^2$  and the scaling exponents  computed with the standard Euclidian metric. Dealing with metrics is convenient to have a direct definition of the scaling exponents but, as previously explained,  a rigorous construction of the quantum metric has not been achieved yet.  Nevertheless, a definition of scaling exponents via measures instead of metrics is also possible, and this is what we discuss below.

We consider the Liouville measure over a bounded domain $D\subset\R^2$:
\begin{equation}\label{chaos:KPZ}
M(dx)=\int e^{\gamma X(x)-\frac{\gamma^2}{2}\E[X(x)^2]}\,dx
\end{equation}
where $\gamma^2<4$ and $X$ is a GFF on $D$, say with Dirichlet boundary conditions. If $M$ were a random metric, we could associate a notion of (random) Hausdorff dimension to this metric. Since $M$ is only a measure, the associated notion of Hausdorff dimension is not straightforward, excepted maybe in dimension $1$ \cite{Benj,Rnew10}. We can nevertheless associate to the measure $M$ a notion of Hausdorff dimension: this just consists in replacing carefully quantities related to distances in the standard definition of Hausdorff dimension with similar quantities defined in terms of measures. This yields: given a Radon measure $\mu$ on $\R^d$ and $s\in [0,1]$, we define for a Borelian set $A$ of $\R^d$:
$$H^{s,\delta}_\mu(A)= \inf \big\{\sum_k \mu(B_k)^{s} \big\}$$ where the infimum runs over all the covering $(B_k)_k$ of $A$ with open Euclidean balls with radius $r_k\leq \delta$. Since the mapping $\delta>0\mapsto H^{s,\delta}_\mu(A)$ is decreasing, we can define the  s-dimensional $\mu$-Hausdorff metric outer measure:
$$H^{s}_\mu(A)=\lim_{\delta\to 0}H^{s,\delta}_\mu(A).$$
The limit exists but may be infinite. Since $H^{s,\delta}_\mu$ is metric, all the Borelian sets are $H^s_\mu$-measurable. The $\mu$-Hausdorff dimension of the set $A$ is then defined as the value
\begin{equation} 
{\rm dim}_\mu(A)=\inf\{s\geq 0; \,\,H^s_\mu(A)=0\}.
\end{equation}
Notice that ${\rm dim}_\mu(A)\in [0,1]$.
When $\mu$ is diffuse (without atoms), the $\mu$-Hausdorff dimension of a set $A$ can also be expressed as:
\begin{equation}\label{HSD}
{\rm dim}_\mu(A)=\sup\{s\geq 0; \,\,H^s_\mu(A)=+\infty\}.
\end{equation}
Therefore, when $\mu$ is diffuse, the above relations allow us to characterize the $\mu$-Hausdorff dimension of the set $A$ as the critical value at which 
the mapping $s\mapsto H^s_\mu(A)$ jumps from $+\infty$ to $0$.

For a given  compact set $K$ of $D$  (or a random compact set independent of $M$), the KPZ formula establishes a relation between the Hausdorff dimension of $K$ computed with $\mu=M$, call it ${\rm dim}_M(K)$, and the Hausdorff dimension of $K$ computed with $\mu$ equal to the Lebesgue measure, call it ${\rm dim}_{Leb}(K)$. We claim (see \cite{Rnew10} for a proof of this statement, or also \cite{Rnew4}):

\begin{theorem}\label{KPZ}{{\bf KPZ formula [Rhodes, Vargas, 2008]}}
Let $K$ be a compact set of $D$. Almost surely, we have the relation:
$${\rm dim}_{Leb}(K)=(1+\frac{\gamma^2}{4}){\rm dim}_M(K)-\frac{\gamma^2}{4}{\rm dim}_M(K)^2.$$
\end{theorem}

We develop below a heuristic to understand what is behind the KPZ formula, mainly the power law spectrum of the measure as explained in Theorem \ref{th:spectrum}. To begin with, we recall the definition of the $s$-dimensional $M$-Hausdorff metric outer measure
$$H^s_M(K)=\lim_{\delta\to 0}\inf\Big\{\sum_n M(B(x_n,r_n))^{s};\quad K\subset \bigcup_n B(x_n,r_n),\,\,r_n\leq \delta\Big\}.$$
Take the expectation and perform an outrageous inversion of limits:
\begin{align*}
\E[H^s_M(K)] =&\lim_{\delta\to 0}\inf\Big\{\sum_n \E[ M(B(x_n,r_n))^{s}];\quad K\subset \bigcup_n B(x_n,r_n),\,\,r_n\leq \delta\Big\}.
\end{align*}
Now compute the expectations via Theorem \ref{th:spectrum} to get:
\begin{align*}
\E[H^s_M(K)] \asymp &\,C_s \lim_{\delta\to 0}\inf\Big\{\sum_n r_n^{\xi(s)/2};\quad K\subset \bigcup_n B(x_n,r_n),\,\,r_n\leq \delta\Big\}\\
=& C_s H^{\xi(s)/2}_{Leb}(K).
\end{align*}
Because $$\xi(s)/2=\big(1+\frac{\gamma^2}{4}\big)s-\frac{\gamma^2}{4}s^2,$$ we recover at least heuristically the KPZ formula. Nevertheless, we draw attention to the fact that we do not claim that the relation
$$\E[H^s_M(K)] \asymp H^{\xi(s)/2}_{Leb}(K)$$ is true. There are possibly   logarithmic corrections in the choice of the gauge function involved in the definition of the $s$-dimensional Hausdorff measures for such a relation to be true.

\subsubsection*{KPZ formula: expected  box counting version}\label{sec:KPZDuSh}
In this subsection, we summarize the KPZ statements proved in \cite{cf:DuSh}. The KPZ theorem of \cite{cf:DuSh} relies on the notion of expected box counting dimension as a definition of scaling exponents. To state the theorem, one must introduce the following definition:

\begin{definition}{\bf isothermal quantum ball.} For any fixed measure $\mu$ on $D$, let $B^\delta(z)$ be the Euclidean ball centered at $z$
with radius given by $\mu(B^\delta(z)) = \delta$.  If
there does not exist a unique $\delta$ with this property, take the
radius to be $\sup \{ \varepsilon: \mu(B_\varepsilon(z)) \leq \delta
\}$.

When $\mu$ is the measure $M$ of \eqref{chaos:KPZ}, the ball $B^\delta(z)$ is called the isothermal quantum ball
of area $\delta$ centered at $z$. When $\mu$ is the Lebesgue measure then $B^\delta(z)$ is nothing but the Euclidean ball centered at $z$ and radius $\epsilon$ where $\delta = \pi \epsilon^2$, denoted by $B_{\epsilon}(z)$.
\end{definition}

Given a subset $K \subset D$, the $\epsilon$-neighborhood of $K$ is defined by: 
$$B_\epsilon(K) = \lbrace z : B_\varepsilon(z) \cap K \not = \emptyset \rbrace.$$
The {\it isothermal quantum $\delta$-neighborhood} of
$K$ is defined by:
$$B^\delta(K) = \lbrace z : B^\delta(z) \cap K \not = \emptyset \rbrace.$$

Finally, the authors of \cite{cf:DuSh} introduce the notion of scaling exponent. 
Fix $\gamma \in [0,2)$ and let $\lambda$ denote Lebesgue measure on
$D$. A fractal subset $K$ of
$D$ has {\it Euclidean expectation dimension} $2-2x$ and {\it
Euclidean scaling exponent $x$} if the expected area of
$B_\epsilon(K)$ decays like $\epsilon^{2x} =
(\epsilon^2)^x$, i.e.,
$$\lim_{\varepsilon \to 0} \frac{ \log \mathbb E \lambda(B_\varepsilon(X)) }{\log \varepsilon^2} = x.$$
The set $K$ has {\it quantum scaling exponent $\Delta$} if  we have
$$\lim_{\delta \to 0} \frac{ \log \E M (B^\delta(X)) }{\log \delta} = \Delta.$$

\begin{theorem}{\bf [Duplantier, Sheffield, 2008]} \label{KPZscaling}
Fix $\gamma \in [0,2)$ and a compact subset $K$ of $D$.  If
$K$ has Euclidean scaling exponent $x \geq 0$ then it
has quantum scaling exponent $\Delta$, where $\Delta$ is the
non-negative solution to
\begin{equation} \label{KPZDupSheff} 
x = \frac{\gamma^2}{4} \Delta^2 + \left( 1 - \frac{\gamma^2}{4}\right)\Delta.
\end{equation}
\end{theorem}

This theorem also extends to the case where $K$ is a random compact set independent of $M$.

In the physics litterature, the KPZ relation is usually stated under this form (\ref{KPZDupSheff}) in which case $x$ and $\Delta$ are the weights of conformal operators. To get a formulation in terms of dimensions, one must make the correspondence $2-2x \leftrightarrow {\rm dim}_{Leb}(K)$ and  $2-2\Delta \leftrightarrow {\rm dim}_M(K)$.

Let us finally mention that in \cite{cf:DuSh} is also proved  a one dimensional boundary version of KPZ. This corresponds to proving the theorem with the lognormal MRM measure of section \ref{sectiondim1}. 

\begin{rem}{\bf Further comments and references on KPZ.} 
The KPZ formula has been proved in \cite{Benj} in the case of multiplicative cascades in dimension $1$ (see also \cite{ismael} for a multidimensional version), in \cite{cf:DuSh} in the case where $X$ is a $2$-dimensional GFF, and in \cite{Rnew10} (see also \cite{Rnew4}) in the case where $X$ is a  log-correlated infinitely divisible field in any dimension. Roughly speaking, infinitely divisible fields are to the family of random distributions  what L\'evy processes are to the family of stochastic processes. Log-correlated Gaussian fields, like two dimensional Free Fields, are a subclass of log-correlated infinitely divisible fields.
Therefore, the main point here is to draw attention to the fact that   the KPZ formula is a property specific to log-correlated fields: it is neither specific to the dimension, nor to the conformal invariance of the $2d$-GFF, nor to the Gaussian nature: the only point that makes $2d$-Liouville quantum gravity (i.e. Gaussian multiplicative chaos with respect to the $2d$-GFF) satisfy a KPZ relation is the fact that the Green function of the Laplacian in dimension $2$ (and in dimension $2$ only) has a logarithmic singularity.
Also, it may be interesting to know if a "Liouville quantum gravity" picture can be drawn for log-correlated infinitely divisible fields instead of Gaussian Free Fields. 
\end{rem}

\subsubsection*{Liouville quantum gravity and KPZ on Riemannian surfaces}
One may wonder what becomes Liouville quantum gravity and the  KPZ formula on a $n$-dimensional Riemannian manifold $(S,g)$ where $g$ is the Riemannian tensor of the manifold. By Liouville quantum gravity, we mean here a Gaussian multiplicative chaos with respect to a Gaussian distribution $X$ defined on the manifold. As long as the random Gaussian distribution $X$ possesses a kernel of $\sigma$-positive type, Kahane's theory allows to define a Gaussian multiplicative chaos associated to this Gaussian distribution. If the covariance kernel of the Gaussian distribution is of the type \eqref{Klog} (where  $\rho$ is the distance   associated to the Riemannian metric) and the measure $\sigma$ in \eqref{eq:deformelle} is the volume form on $S$, then the non-degeneracy conditions of the chaos is $\gamma^2<2n$ (Theorem \ref{th:degeneracy}). Since a $n$-dimensional Riemann surface is locally isometric to the unit ball of $\R^n$, we deduce from \cite{Rnew10} (or \cite{Rnew4}) that the KPZ formula holds for the Gaussian multiplicative  chaos $M$ on this Riemann surface: it reads
\begin{theorem}\label{KPZS}{{\bf KPZ formula on Riemann manifolds [\cite{Rnew10}, 2008]}}
Let $K$ be a compact set of $S$. Almost surely, we have the relation:
$${\rm dim}_{\sigma}(K)=(1+\frac{\gamma^2}{2n}){\rm dim}_M(K)-\frac{\gamma^2}{2n}{\rm dim}_M(K)^2.$$
\end{theorem}
In particular, we see that the curvature of the surface does not affect the KPZ relation. For instance, we can consider  a $2$-dimensional Riemann surface, like  a sphere or an hyperbolic half-plane, and  the GFF on a domain of this surface with appropriate boundary conditions in order to define the associated Liouville measure. As explained above, in dimensions different from $2$, the GFF does not possess logarithmic correlations so that it does not make sense to look for KPZ relations based on the GFF.  Nevertheless, in dimensions different from $2$, it is plain to construct other log-correlated Gaussian distributions $X$: various examples of log-correlated Gaussian fields are described in the present manuscript but also in \cite{DSRV3}).

Another situation of interest is to consider massive or generalized Free Fields. On a domain $D\subset \R^2$, the Massive Free Field (MFF) with Dirichlet boundary conditions is defined as a standard Gaussian in the Hilbert space defined as the closure of Schwartz functions over $D$ with respect to the inner product
$$(f,g)_h=m^2(f,g)_{L^2(D)}-(f,\triangle g)_{L^2(D)}.$$
The real $m>0$ is called the mass. Its action on $L^2(D)$  can   be seen as a Gaussian distribution with covariance kernel given by  the Green function $G_m^D$ of the operator $m^2-\triangle$, i.e.:
$$(m^2-\triangle)G_m^D(x,\cdot)=2\pi\delta_x $$
with Dirichlet boundary conditions. When $D$ is the whole plane, the massive Green kernel is a star-scale invariant kernel in the sense of \cite{Rnew1}. We may also consider {\it Generalized Free Fields} as defined in \cite{glimm}.
\begin{definition}
A {\it Generalized Free Field} with Dirichlet boundary conditions over a domain $D\subset \R^2$ is defined as a random centered  Gaussian distribution  (say on the space of Schwartz functions on $D$) the covariance kernel of which is given by
$$G^D_\varrho(x,y)=\int_0^{+\infty}G^D_m(x,y)\,\varrho(dm).$$
where $G_m^D$ is the massive Green function on $D$ with mass $m$ and $\varrho$ is a  Radon measure on $\R_+$, called the K\"allen-Lehmann weight, satisfying 
$$\int_0^{1}-\ln m\,\varrho(dm)<+\infty$$ and 
$$\forall k\in\N,\quad \int_0^{+\infty}m^k\,\varrho(dm)<+\infty.$$
\end{definition}
The construction of such a field can be straightforwardly adapted from the previous white noise decomposition of the free field. Denote by $p_D(t,x,y)$ the transition densities of a Brownian motion killed upon touching the boundary of $D$. Consider a Gaussian noise, white in space and time and $\varrho$-colored in mass, i.e. a Gaussian random measure $W(dx,ds,dm)$ distributed on $D\times R_+\times\R_+$ such that for all Borel sets $A,A'\subset D$, $B,B',C,C'\subset \R_+$ 
$$\E[W(A,B,C)W(A',B',C')]=|A\cap A'| \,|B\cap B'|\, \varrho(C\cap C'),$$
 where $|A|$ stands for the Lebesgue measure of $A$. The generalized free field can then be defined 
 $$X(x)=\sqrt{\pi}\int_{D\times \R_+\times\R_+} e^{-\frac{m^2}{4}s}p_D(\frac{s}{2},x,y)\,W(dy,ds,dm),$$ and cutoff approximations:
 $$\forall \epsilon\geq 0,\quad X_\epsilon(x)=\sqrt{\pi}\int_{D\times [\epsilon,+\infty[\times\R_+} e^{-\frac{m^2}{4}s}p_D(\frac{s}{2},x,y)\,W(dy,ds,dm).$$
 For all these fields, theory of Gaussian multiplicative chaos and Theorem \ref{KPZS} apply. Let us just stress that the result of the KPZ formula remains unchanged for Massive Free Fields whereas $\gamma^2$ must be replaced with $\gamma^2\varrho(\R_+)$ for Generalized Free Fields. 
  
A last example of interest is the {\it boundary Liouville measure} introduced in \cite{cf:DuSh}. The authors suggests to consider a smooth domain $D$ of $\R^2$ together with  a GFF $X$ on $D$ with free boundary conditions  and to define the boundary Liouville measure on  $\partial D$ as
 $$\nu(dx)=\lim_{\epsilon\to 0}\int_{ \cdot}e^{\frac{\gamma}{2}X_\epsilon(x)-\frac{\gamma^2}{4}\ln\frac{1}{\epsilon}}\,dx$$
 where $dx$ stands for the length measure on $\partial D$ and $X_\epsilon(x)$ is the mean value of $X$ over $\partial B(x,\epsilon)\cap D$. Actually, the boundary Liouville measure is nothing but a Gaussian multiplicative chaos over a $1$-dimensional Riemannian manifold (in fact $C^1$ is enough). By noticing that the correlations are given by 
\begin{equation}\label{boundcorrel}
\E[\frac{\gamma}{2}X_\epsilon(x)\frac{\gamma}{2}X_\epsilon(y)]=\frac{\gamma^2}{2}\ln_+\frac{1}{|x-y|}+g(x,y)
\end{equation}
 for some continuous bounded function $g$, Theorem \ref{th:degeneracy} ensures that the boundary Liouville measure is non-degenerate provided that $\gamma^2<4$. Furthermore, the boundary KPZ holds: just take care of replacing $\gamma^2$ in Theorem \ref{KPZS} by $\gamma^2/2$ because of the unusual normalization  in \ref{boundcorrel} (the term in front of the log is $\frac{\gamma^2}{2}$, averaging along semi-circles yields an extra factor $2$). The expected box counting KPZ formula for the boundary Liouville measure is proved in \cite{cf:DuSh} whereas the Hausdorff dimension version of the KPZ formula for boundary Liouville measure is proved in \cite{Rnew10}.

\subsection{Convergence of discrete Liouville measures on isoradial graphs}\label{cvdlm}

Here we consider a  planar (isoradial) graphs. We weight the vertices of this graph by the exponential of the Discrete Gaussian Free Field (DGFF for short) to obtain a discrete Liouville measure: the measure having as density exponential of the DGFF with respect to the discrete canonical volume measure on the graph. We will prove that this measure weakly converges in law towards the Liouvile measure as  the mesh size of the graph converges to $0$.

Before stating a clear theorem, we need to explain the framework in further details. The reader may find in \cite{chelkak} all the basic tools (and much more) about isoradial graphs described below. We stick to the notations used in \cite{chelkak}. A planar graph $\Gamma$ embedded in $\C$ is called isoradial iff each face is inscribed into a circle of common radius $\epsilon$. If all circle centers are inside the corresponding faces, then one can naturally embed the dual graph $\Gamma^*$ in $\C$ isoradially with the same $\epsilon$, taking the circle centers as vertices of $\Gamma^*$. The name rhombic lattice is sometimes used for these graphs because all the quadrilateral faces of the corresponding bipartite graph $\Lambda$ (having vertices $\Gamma\cup\Gamma^*$) are rhombi with sides of length $\epsilon$. We will make the following assumption (see \cite{chelkak})
\begin{equation}\label{angle}
\text{the rhombi half-angles are uniformly bounded away from } 0 \text{ and }\frac{\pi}{2}.
\end{equation}
Roughly speaking, $\Lambda$ does not possess too flat rhombi. This entails that the Euclidean distance between the vertices of $\Gamma$ is comparable to the graph distance. For such graphs, the canonical volume measure $\mu_\Gamma(A)$ of a subset $A\subset\Gamma$ is given by
$$ \mu_\Gamma(A)=\sum_{z\in A}W_z$$ 
where the weight $W_z$ of $z\in \Gamma$ is the Lebesgue measure of the face of the dual graph containing $z$.

We consider a sequence $(\Gamma_n)_n$ of isoradial graphs as indicated above with radius $(\epsilon_n)_n$. For all the quantities defined above, the subscript $n$ means that it is related to the graph $\Gamma_n$. For instance, $ \mu_n$ stands for the volume measure of the graph $\Gamma_n$. We also assume that the radius $\epsilon_n$ of $\Gamma_n$ goes to $0$ as $n\to\infty$. 

Let us now consider a bounded simply connected domain $D$ of $\C$. Let us denote by $D_n$ the graph $D\cap \Gamma_n$, i.e. we keep the vertices and edges that entirely lie in $D$. We consider a DGFF $X_n$ on the vertices of $D\cap \Gamma_n$. Recall that the DGFF is a collection $(X_n(z))_{z \in D_n}$ of centered Gaussian random variables with covariance kernel given by $2\pi G_{D_n}$, where $G_{D_n}$ denotes the Green function  on $D_n$ with $0$-boundary condition (see \cite[Definition 2.6]{chelkak}). 

Then, for $\gamma\in [0,2[$, we define the discrete Liouville measure on $D_n$:
\begin{equation}
M_{n,\gamma}(dz)=e^{\gamma X_n(z)-\frac{\gamma^2}{2}\E[X_n(z)^2]}\,\mu_n(dz).
\end{equation}

\begin{theorem}\label{DLQG}
Let $D$ be a bounded simply connected open domain in $\C$. For $\gamma\in [0,2[$, the discrete Liouville measures $ (M_{n,\gamma}(dz))_n$ on $D_n$ weakly converge in law towards the Liouville measure on $D$ (see subsection \ref{LQG}), that is the Gaussian multiplicative chaos:
$$M_\gamma(dx)=e^{\gamma X(x)-\frac{\gamma^2}{2}\E[X(x)^2]}\,dx$$ where $X$ is a GFF on $D$ with $0$-boundary condition. 
\end{theorem}

The proof can be found in Appendix \ref{discLQG}. In fact, we will prove the result for the discrete GFF on the square lattice with mesh $\epsilon$ going to $0$. It turns out that the proof of Theorem \ref{DLQG}
also works at criticality (see subsection \ref{sectheory:deriv} for further details and references), yielding the following:

\begin{theorem}\label{DLQG2}
For $\gamma=2$, the discrete critical Liouville measures $ (\sqrt{\ln \frac{1}{\epsilon_n}} M_{n,2}(dz))_n$ on $D_n$ weakly converge in law towards the critical Liouville measure on $D$ (see \cite{Rnew7,Rnew12} or subsection \ref{sectheory:deriv} below), that is the Gaussian multiplicative chaos:
$$M'(dx)=\sqrt{\frac{2}{\pi}}(2\E[X(x)^2]-X(x))e^{2 X(x)-2\E[X(x)^2]}\,dx$$ where $X$ is the GFF on $D$ with Dirichlet boundary condition. 
\end{theorem}

\subsection{Kolmogorov-Obhukov model in turbulence}
We refer to \cite{cf:Fr} for an introduction to the statistical theory of 3 dimensional turbulence. Consider a stationary flow at high Reynolds number, that is when the velocity of the fluid is large in comparison with the viscosity forces. It is believed that at small scales the velocity field of the flow is homogeneous and isotropic in space . By small scales we mean scales much smaller than the integral scale $R$ characteristic of the time stationary force driving the flow. In the works \cite{cf:Kol} and \cite{cf:Obu}, Kolmogorov and Obukhov proposed to model the mean energy dissipation per unit mass in a ball $B(x,l)$ of center $x$ and radius $l<<R$ by a random variable $\epsilon_{l}$ such that $\ln(\epsilon_{l})$ is normal  with variance $\sigma_{l}^2$ given by:
\begin{equation*}
\sigma_{l}^2=\lambda^2\ln(\frac{R}{l})+A
\end{equation*}      
where $A$ is a constant and $\lambda^2$ is the intermittency parameter. As noted by Mandelbrot (\cite{cf:Man}), the only way to define such a model is to construct a random measure $\epsilon$ by a limit procedure. Then, one can define $\epsilon_{l}$ by the formula:
\begin{equation*}
 \epsilon_{l}=\frac{3<\epsilon>}{4 \pi l^3}\epsilon(B(x,l))
\end{equation*}
where $<\epsilon>$ is the average mean energy disspation per unit mass. 
Formally, one is looking for a random measure $\epsilon$ such that:
\begin{equation}\label{eq:epsilon}
\forall A \in \mathcal{B}(\R^d), \quad \epsilon(A)=\int_{A}e^{\gamma X(x)-\frac{\gamma^2}{2}E[X(x)^2]}dx
 \end{equation}  
where $(X(x))_{x \in \R^d}$ is a "Gaussian field" whose covariance kernel $K$ is given by \eqref{euclid:K}. Therefore, one can give a rigorous meaning to energy dissipation (\ref{eq:epsilon}) by using Gaussian multiplicative chaos theory. Let us mention here that the objective of describing a stochastic representation of the velocity field is a much more ambitious task (see subsection \ref{matrix} for more comments and perspectives on this topic).

\subsection{Decaying Burgers turbulence}
Consider the Burgers equation
\begin{equation}
\partial_t v+v\nabla v=\nu \nabla v+f(x,t)\quad \text{with initial condition }v(x,0)=v_0(x).
\end{equation}
$\nu$ is the {\it viscosity } parameter. When $f\not =0$, this equation is called {\it (randomly) forced Burgers equation}. When $f=0$ and $v_0(x)\not =0$, this equation is called {\it decaying Burgers turbulence}. 

We consider here the case $f=0$ with random initial data $v_0$. The solution of {\it decaying Burgers turbulence}, via Hopf-Cole transform, is given by $v(t,x)=\nabla (-2 \nu\ln Z(t,x))$, where
\begin{equation}\label{hopf}
Z(t,x)=\int_\R e^{-\frac{1}{2\nu}\frac{|y-x|^2}{2t}-\frac{1}{2\nu}V(y)}\frac{dy}{\sqrt{4\pi \nu t}}.
\end{equation}
where $v_0=\nabla V$. If we choose $V$ as a log-correlated Gaussian random potential, i.e. with covariance kernel of the type \eqref{euclid:K} then $Z$ appears as a Gaussian multiplicative chaos of the type \eqref{eq:deformelle} and integrating measure $\sigma$ given by the standard heat kernel on $\R$. This situation corresponds to a power law correlated random profile of Gaussian distributed initial velocities: 
$$\E[v_0(x)v_0(y)]\sim |x-y|^{-2}.$$ Gaussian distributed initial velocities with power law correlations is a subject of interest in Burgers turbulence. The reader may consult \cite{FLR} for an account of physical motivations, further references and a study of the present situation.

It is shown in \cite{FLR} that this equation exhibits a phase transition at a critical viscosity parameter $\nu_c$, which is related to the phase transition of Gaussian multiplicative chaos (see Theorem \ref{th:degeneracy}). For $\nu<\nu_c$, freezing phenomena occur, highlighting a glassy phase conjecturally as that explained in subsection \ref{sectheory:atomic} (frozen phase).
   
\section{Generalizations of the theory}\label{sec:gener}
In this section, we review some generalizations of Kahane's theory while making connections with possible applications.

\subsection{Gaussian multiplicative chaos at criticality}\label{sectheory:deriv}
Kahane's construction of Gaussian multiplicative chaos  makes sense when the factor $\gamma$ appearing in \eqref{euclid:M} satisfies $\gamma^2<2d$. This gives rise to the issue of constructing random measures in the same spirit for $\gamma^2\geq 2d$. Necessarily these measures will present a different structure. The case $\gamma^2=2d$ is of special interest since it corresponds to a phase transition. We  call this situation the {\it critical case}. Kahane's theory ensures that the associated martingale $(M_n)_n$ appearing in \eqref{eq:deformellen} almost surely converges towards $0$. Several approaches are possible in order to make sense of a suitable random measure corresponding to criticality $\gamma^2=2d$. It turns out that all these approaches are conjecturally the same.

The first approach is based on the convergence  of  a fundamental object that is called {\it derivative martingale}.  Such an object has been intensively studied in the case of multiplicative cascades, branching random walks  \cite{BiKi,Kyp} or branching Brownian motions \cite{neveu}. The convergence is achieved in \cite{Rnew7} in the context of Gaussian multiplicative chaos associated to star scale invariant kernels \eqref{structmulti} and may be achieved for other type kernels provided that one can use a white-noise decomposition cutoff (see \cite{Rnew12}), for instance including the GFF in a bounded domain (see previous subsection \textbf{White noise decomposition of the GFF}). The limiting measure can formally be written as
\begin{align}\label{derivintro}
M'(A) =&  \int_A\bigl(\gamma\,\E[X^2(x)]-X(x)\bigr)e^{\gamma X(x)-\frac{\gamma^2}{2}\E[X^2(x)]}\,dx \quad \text{with }\gamma=\sqrt{2d} . 
\end{align}
The reader may object that this construction should be possible for every value $\gamma^2\leq 2d$ and therefore argue that the interest for the only value $\gamma^2= 2d$ is not natural.  It turns out that there is an abrupt change in the behaviour of this object at the critical value: the exponential term penalizes  those points $x$ where the "value" of the process $X(x)$ lies above the expectation term $\gamma \E[X(x)^2]$ with a strength that depends on $\gamma$: for $\gamma^2<2d$ the measure may be well defined but the penalization term is not strong enough so that the sign of the term $\bigl(\gamma\,\E[X^2(x)]-X(x)\bigr)$ alternates. For $\gamma^2=2d$, the exponential term is strongly penalizing, forcing the measure $M'$ to be nonnegative, which is not straightforward at first sight. Furthermore, in the situation when the random distribution $X$ possesses a star scale invariant kernel, the limit of the derivative martingale $M'$ yields a solution to the star equation for the only value $\gamma^2= 2d $. Indeed, when trying to derive this scaling relation for $M'$, you are left with an unusual extra term proportional to the standard chaos $M$. For $\gamma^2<2d$, this term does not vanish whereas for $\gamma^2=2d$ Kahane's theory ensures that this extra term disappear, making the measure $M'$ star scale invariant. It is also proved in \cite{Rnew7} that the random measure $M'$ has almost surely full support and no atom. Further improvements are made in \cite{basic}: the authors  determine  the exact asymptotics of the right tail of the distribution of the total mass of the measure, and an almost sure upper bound for the modulus of continuity of the cumulative distribution function of the measure. A lower bound for the increments of the measure is also investigated,  showing that the measure is supported on a set of Hausdorff dimension $0$.  

\begin{figure}[h]
\begin{center}
\includegraphics[angle=0,width=.93290\linewidth]{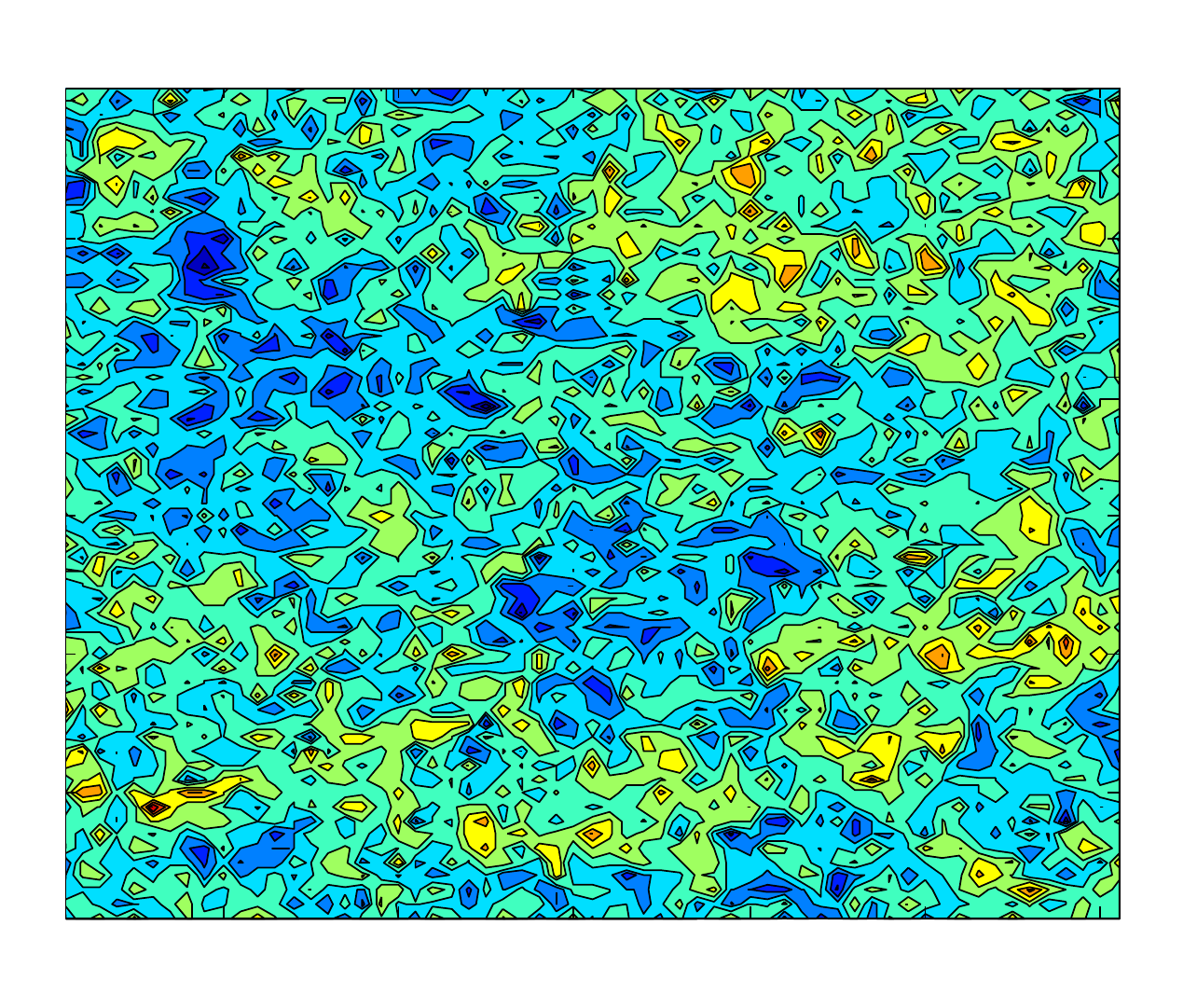}
\caption{Height landscape of the derivative martingale measure plotted with a logarithmic scale color-bar, showing that the measure is very ``peaked'' (for $t=12$, a multiplicative factor of about $10^{8}$ stands between extreme values, i.e., between warm and cold colors).}
\label{fig1}
\end{center}
\end{figure}

Another possible approach is to find a suitable renormalization of \eqref{eq:deformellen} in order to get a non trivial limit. This approach is carried out in \cite{Rnew12}, based on the works \cite{lalley} for branching Brownian motions or \cite{AidShi,HuShi} for branching random walks. The renormalization turns out to be the square root of the variance of the Gaussian field $X_n$, the convergence holds in probability  and the limit measure is  the derivative measure $M'$,  up to a perfectly explicit multiplicative factor. 

This approach is very convenient in order  to use Kahane's convexity inequalities with the measure $M'$, which is not obvious at first sight when looking at the expression \eqref{derivintro}. In particular, a complete description of the moments of $M'$ is achieved in \cite{Rnew7,Rnew12} as well as the computation of the power-law spectrum of $M'$. 

Let us point out that another approach of criticality via a KPP equation has been investigated in the case of multiplicative cascades  \cite{BKNSW} or branching Brownian motion \cite{lalley}. No  Gaussian multiplicative chaos counterpart has ever been derived rigorously. The main  reason is the structure of correlations, which are more intricate than in the discrete framework.  But this is clearly a perspective of interest.

In Liouville quantum gravity, the Liouville measure at criticality $\gamma=2$ corresponds to a central charge $c=1$. It is expected, for instance,  to be the scaling limit of the $O(n=2)$ loop model or the $Q=4$- Potts model  (see \cite{Nienhuis} for a description of these models). Let us mention that, in the standard physics literature about $c=1$ Liouville field theory  \cite{BKZ,GM,GZ,GrossKleban,KKK,Kleb2,Parisi,Polch}, the Tachyon field $\varphi e^{2 \varphi}$ presents   an unusual dependence on the Liouville field (or GFF) $\varphi$ (compare with the usual form $ e^{\gamma \varphi}$ for $c<1$). This is clearly in relation with the mathematical formulation \eqref{derivintro}. The KPZ formula (Theorem \ref{KPZ}) has been proved in \cite{Rnew12} at criticality, i.e. for $\gamma=2$.
\subsection{Atomic Gaussian multiplicative chaos}\label{sectheory:atomic}

We know review some recent progress motivated by the super-critical regime $\gamma^2>2d$. They have been mainly inspired by the seminal paper \cite{durrett} in the context of multiplicative cascades. Typically, a new class of multiplicative chaos emerges, called {\it atomic multiplicative chaos}, which can  be split in two parts. 

\subsubsection*{Dual phase}

Before coming to mathematical considerations, let us first briefly outline the physics motivations for considering the dual phase. The name {\it dual phase} comes from physics literature about Liouville quantum gravity. Remind that the coupling constant $\gamma$ appearing in \eqref{i.metric} or \eqref{measLQG} is determined by the underlying model coupled to gravity. It is related to the so-called {\it central charge} of the model  by the relation \cite{cf:KPZ}
$$c=1-\frac{6}{4}\Big(\gamma-\frac{4}{\gamma}\Big)^2.$$
The central charge belongs to $]-\infty,1]$. The special case $c=1$ (and $\gamma=2$) corresponds to criticality. Otherwise, the central charge belongs to $]-\infty,1[$ yielding two possible $\gamma> 0$. 
The first solution $\gamma$ belongs to $]0,2[$ and is given
\begin{equation}\label{gam1}
\gamma=\frac{\sqrt{25-c}-\sqrt{1-c}}{\sqrt{6}} . 
\end{equation}
 It corresponds to the {\it standard branches of Liouville gravity} detailed in subsection \ref{LQG}. The second solution, call it $\bar{\gamma}$, belongs to $]2,+\infty[$:
\begin{equation}\label{gam2}
\overline{\gamma}=\frac{\sqrt{25-c}+\sqrt{1-c}}{\sqrt{6}} . 
\end{equation}
Observe that $\gamma \overline{\gamma}=4$.

In a series of papers \cite{alvarez,Das,dupdual,Kleb1,Kleb2,Kleb3,korchemsky},  physicists have  investigated what they called {\it the other/non-standard/unconventional/dual branches of gravity} initially via modified random matrix models generating random surfaces  in order to interpret these other possible values of the coupling constant $\overline{\gamma}\in]2,+\infty[$.  They have noticed several intriguing relations between the standard and dual branches of the Liouville action, laying the foundations of what they called {\it Duality of Liouville quantum gravity}.  More recently, duality of Liouville quantum gravity has been digged up in \cite{PRL,Dup:houches} at an heuristic yet interesting level. 

The purpose of what follows is to present mathematical results and conjectures on the mechanisms involved in this dual phase, which turns out to be a very rich area from the mathematical or physics angles.  To understand mathematically how to handle the dual phase, let us continue this pedagogical introduction in the context of Liouville quantum gravity. So we consider a bounded domain $D$ of $\R^2$ and a GFF $X$ on $D$ with Dirichlet boundary conditions. The first observation is that Theorem \ref{th:necess2} tells us that the standard chaos
$$M_{\bar{\gamma}}(dx)=e^{\bar{\gamma}X(x)-\frac{\bar{\gamma}^2}{2}\E[X(x)^2]}\,dx$$ reduces to $0$ because $\bar{\gamma}^2>4$. Therefore another construction has to be found. To understand which construction is involved, let us assume for a while that the chaos $M_{\bar{\gamma}}$ does not reduce to $0$ (though it does). Then, as explained in Theorem \ref{th:spectrum} (and in its proof), the chaos $M_{\bar{\gamma}}$ would then satisfy the scaling relation
\begin{equation}
M_{\bar{\gamma}}(B(x,r))=r^2e^{\bar{\gamma} \Omega_r-\frac{\bar{\gamma}^2}{2}\E[\Omega_r^2]}M_{\bar{\gamma}}(B(x,1))
\end{equation}
for some Gaussian random variable $\Omega_r$ independent of $M_{\bar{\gamma}}$ and variance $\ln \frac{1}{r}$. A straightforward computation shows that this relation may be rewritten as:
\begin{equation}
M_{\bar{\gamma}}(B(x,r))=\Big(r^2e^{\gamma \Omega_r-\frac{ \gamma ^2}{2}\E[\Omega_r^2]}\Big)^{1/\alpha}M_{\bar{\gamma}}(B(x,1))
\end{equation}
where $\alpha=\frac{4}{\bar{\gamma}^2}$ (recall that $\gamma \bar{\gamma}=4$). This relation suggests that if $M_{\bar{\gamma}}$ were non degenerate, it should scale  for small $r$ like
\begin{equation}\label{dualsmoothing}
M_{\bar{\gamma}}(B(x,r))\simeq \Big(M_{\gamma}(B(x,r))\Big)^{1/\alpha}.
\end{equation}
The convenient fact of the above relation is that  the Gaussian multiplicative chaos $M_\gamma$ in the right-hand side  is non trivial since $\gamma^2<4$. Therefore, we are heuristically looking for a measure that may be interpreted as the $\frac{1}{\alpha}$-th root of the standard chaos $M_\gamma$. From the mathematical angle, this $\frac{1}{\alpha}$-th root perfectly makes sense in terms of independently scattered random measures with prescribed control measure: roughly speaking, this consists in throwing a stable point process over a landscape described by the standard Gaussian multiplicative chaos $M_\gamma$. More precisely, the law of this measure can be generated as follows:
\begin{itemize}
\item sample the standard Gaussian multiplicative chaos  
$$M_\gamma (dx)=  e^{  \gamma X(x)-\frac{\gamma^2}{2}\E[X(x)^2]}\,dx.$$
\item sample a random measure $M_{\bar{\gamma}}$ whose law, conditionally to $M_\gamma$, is that of  an independently scattered random measure  characterized by $$\forall q\geq 0,\quad \E[e^{-q M_{\bar{\gamma}}(A)}|M_\gamma]=e^{-q^\alpha M_\gamma(A)}. $$
\end{itemize}
This construction is called subordination procedure. Interestingly,  it is a purely atomic random measure as  suggested in physics literature.
 
Until now, we have identified, at least in law, what kind of object we are looking for to model the dual branch of Liouville quantum gravity. Let us also mention here that we recover the suggestion made in \cite{Dup:houches} concerning the dual measure. It seems that the arguments used in \cite{Dup:houches} to guess the exact form of the dual Liouville measure are based on the knowledge of the scaling exponents involved in the dual branch of gravity.

From the theoretical angle, this is the beginning of interesting questions. First, we mention that the scaling relation \eqref{dualsmoothing} has an exact equivalent in the context of branching random walks \cite{durrett,Biggf,Liu}, the solutions of which are called "the fixed points of the smoothing transform". Interestingly, it is proved that the corresponding scaling relation has a unique solution in law, which is the exact equivalent for branching random walks of the measure $M_{\bar{\gamma}}$ described above. This is a strong argument to validate the law of the above measure. The question of uniqueness in the context of Gaussian multiplicative chaos theory can be formulated in terms of star scale invariance and will be discussed in subsection \ref{conj:star}. Second, the above  construction of $M_{\bar{\gamma}}$ is not quite satisfactory. Indeed, the measure $M_{\bar{\gamma}}$ is not a measurable function of the only field $X$: the measure $M_{\bar{\gamma}}$ also depends on the randomness originating from the sampling of the point process. It is not a measure intrinsically generated by the only field $X$. Third, the construction of $M_{\bar{\gamma}}$ by subordination  does not  appear as an almost sure limit of a suitably renormalized sequence in the spirit of \eqref{eq:deformellen} for instance. We will see that this does not allow us to understand the mechanisms involved in duality and that this is   also related to the latter objection. Finally, a KPZ formula for the dual measure has to be proved rigorously.

\vspace{2mm} 
A different way of constructing the measure $M_{\bar{\gamma}}$   has been carried out in \cite{Rnew4}.  It is important to understand here that the main differences with the subordination procedure take place at the level of almost sure properties of the  measure $M_{\bar{\gamma}}$. Basically, the construction consists in  replacing the $dx$ measure in \eqref{euclid:M} by an independently scattered infinitely divisible random measure, call it $n_\alpha$ for some $\alpha\in]0,1[$,  whose law  is characterized by
$$\forall q\geq 0,\quad \E[e^{-q n_\alpha(A)}]=e^{-  q^\alpha |A|}$$
for all Borelian set $A$ of $D$. The random measure $n_\alpha$ is assumed to be independent of the free field $X$. For  $\bar{\gamma}^2>4$, it then makes sense to define an approximating sequence in the spirit of \eqref{eq:deformellen} (with the same notations):
\begin{equation}\label{theory:atomic}
M_{n,\bar{\gamma}}(dx)=  e^{\bar{\gamma} X_n(x)-\frac{\alpha\bar{\gamma}^2}{2 }\E[ X_n(x)^2]}\,n_\alpha(dx).
\end{equation}
\begin{theorem}{\bf Convergence of the atomic Gaussian multiplicative chaos [\cite{Rnew4}, 2012]}
Almost surely in $n_\alpha$, the sequence of random measures $(M_{n,\bar{\gamma}})_n$ converges in $\Pb^X$-probability  towards a non trivial limiting measure $$M_{\bar{\gamma}}(dx)=e^{\bar{\gamma} X(x)-\frac{\alpha\bar{\gamma}^2}{2 }\E[ X(x)^2]}\,n_\alpha(dx),$$ called  atomic Gaussian multiplicative chaos, which is purely atomic.  
\end{theorem}

We point out that the above construction is not a Gaussian multiplicative chaos in the usual sense. Indeed the lognormal weight is not normalized to have expectation $1$. The expectation blows up, giving rise to a situation that qualitatively deeply differs from standard Gaussian multiplicative chaos theory. Therefore, the simplicity of Kahane's construction is lost in this construction: \\
1) we do not deal with martingales because the weights are not normalized,\\
2) these measures are not integrable because the measure $n_\alpha$ is not.

\begin{figure}[h]
\begin{center}
\includegraphics[width=15cm,height=6cm]{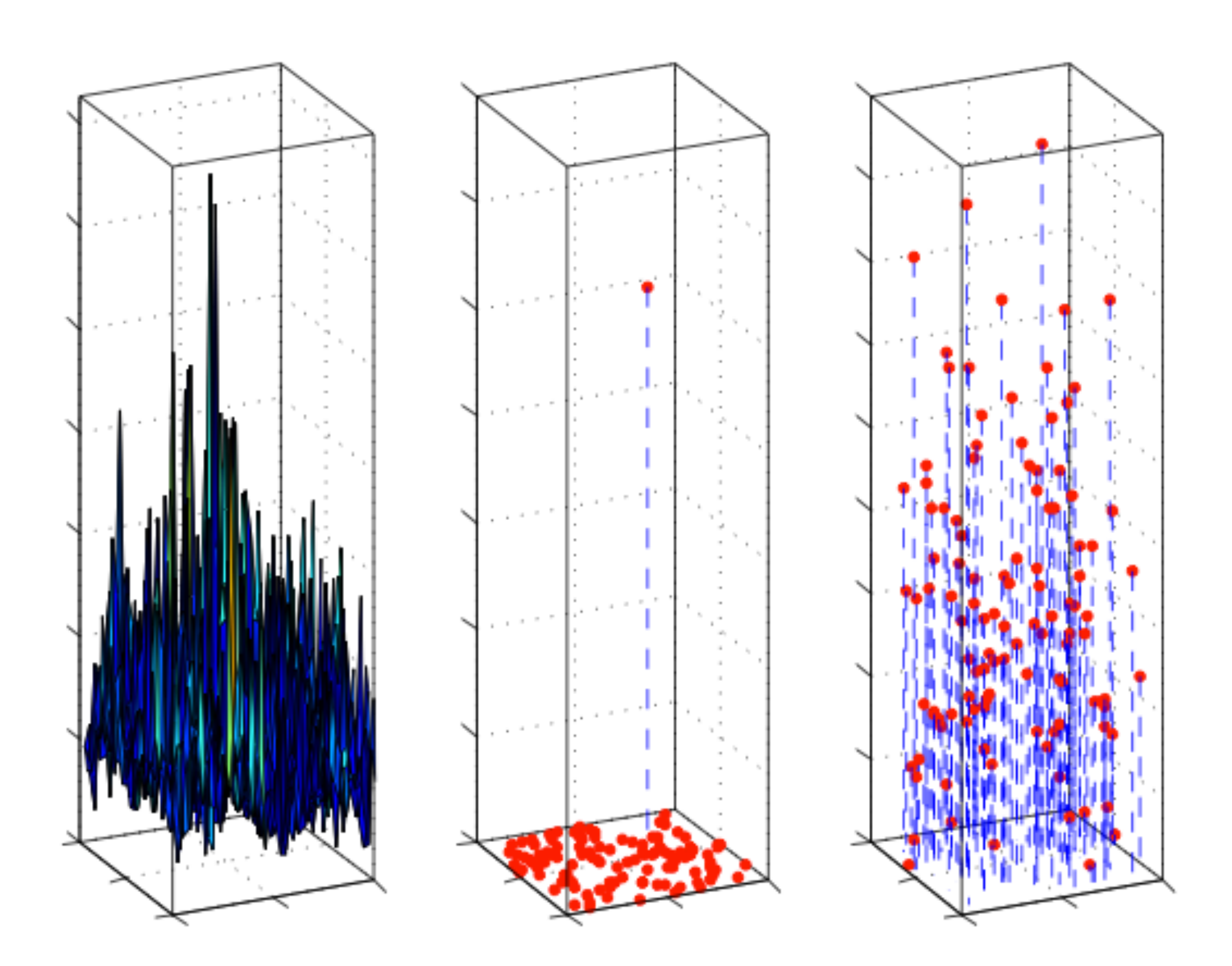}
\caption{Chaos and dual chaos for the value $\gamma^2=0,25$ (and then $\alpha=0.051$). Left: ``density" of the standard chaos. Middle: position and weights of the atoms of the dual measure. Right: position and weights of the atoms of the dual measure  with a logarithmic ordinate scale. 
}
\label{figdual}
\end{center}
\end{figure}

Surprisingly, the atomic Gaussian multiplicative chaos is independent of the measure $n_\alpha$ used in the construction:
\begin{theorem}
The random measure $M_{\bar{\gamma}}$ obtained as the limit in probability of the sequence $(M_{n,\bar{\gamma}})_n$ defined by \eqref{theory:atomic} is independent of the stable measure $n_\alpha$. It is therefore a measurable function of the only field $X$.
\end{theorem}

\noindent {\it Proof.} It is standard that the stable measure $n_\alpha$ can be decomposed as 
$$n_\alpha(A)=\int_A\int_{\R_+}z \,m_\alpha(dx,dz)$$ where $m_\alpha$ is a Poisson random measure on $\R^d\times \R_+$ with intensity given by $dx\frac{dz}{z^{1+\alpha}}$. Let us first prove  that  $\int_B\int_c^{+\infty}e^{\frac{\gamma}{\alpha }X(x)-\frac{\gamma^2}{2 \alpha}\E[X(x)^2]}z\,m_\alpha(dx,dz) =0$ for all $c>0$ and all balls $B$.
To this purpose, observe that this term is the limit in probability as $n\to\infty$ of
$$\int_B\int_c^{+\infty}e^{\frac{\gamma}{\alpha }X_n(x)-\frac{\gamma^2}{2 \alpha}\E[X_n(x)^2]}z\,m_\alpha(dx,dz),$$ which is nothing but  a finite sum of terms like
$$Z_i e^{\frac{\gamma}{\alpha }X_n(x_i)-\frac{\gamma^2}{2 \alpha}\E[X_n(x_i)^2]},$$ where $Z_i$ are random variables independent of the field $X$. By the strong law of large numbers, for all $x\in\R$, almost surely in $X$, we have
$$\lim_{n\to\infty}\frac{\gamma}{\alpha }X_n(x)-\frac{\gamma^2}{2 \alpha}\E[X_n(x)^2]=-\infty.$$ Our claim follows. The Kolmogorov $0-1$ law then proves the result. \qed

\vspace{2mm} 
Furthermore, the measure $M_{\bar{\gamma}}$ may be seen as the volume form associated to a Riemannian tensor of the form
$$e^{\bar{\gamma} X(x)-\frac{\alpha\bar{\gamma}^2}{2 }\E[ X(x)^2]}\,(n_\alpha(dx))^2$$
This can be seen by regularizing the measure $n_\alpha$ with a sequence of mollifiers $(\rho_n)_n$: $M_{\bar{\gamma}}$ appears as the limit as $n\to\infty$ of the volume forms associated to the smooth metric tensors obtained by mollifying the measure $n_\alpha$ with the sequence $(\rho_n)_n$.  For instance, consider the setup drawn in subsection \ref{cvdlm}. Instead of considering the canonical volume form $\mu_n$ on $D_n$, meaning each vertice $z$ of $D_n$ has a weight corresponding to the Lebesgue measure of the face of the dual graph containing $z$, we assign to each vertic $z$ of $D_n$ the weight $W_z^\alpha$ of the $n_\alpha$-measure of the dual face of $D_n$ containing $z$. We define the measure
$$\mu_{n,\alpha}(A)=\sum_{z\in A}W_z^\alpha.$$
This is some kind of Bouchaud trap model on isoradial graphs. We now consider this kind of "Bouchaud trap model" in the gravitational dressing, i.e. we consider the discrete GFF $X_n$ on $D_n$ and define the discrete dual Liouville measure for $\bar{\gamma}\in ]2,+\infty[$
$$ M_{n,\bar{\gamma}}(dz)=  e^{\bar{\gamma} X_n(z)-\frac{\alpha\bar{\gamma}^2}{2 }\E[ X_n(z)^2]}\,\mu_{n,\alpha}(dz).$$
It is straightforward to see that the sequence of measures $(M_{n,\bar{\gamma}})_n$ converges in law towards the measure  $M_{\bar{\gamma}}$. The same conclusion holds if, instead of assigning a stable law on each face of the dual graph, we assign a Random Energy Model at low temperature on the collection of faces of the dual graph of $D_n$.
 
 \begin{figure}[h]
\centering
\subfloat[$\gamma^2=0.01$]{\includegraphics[width=0.44\linewidth]{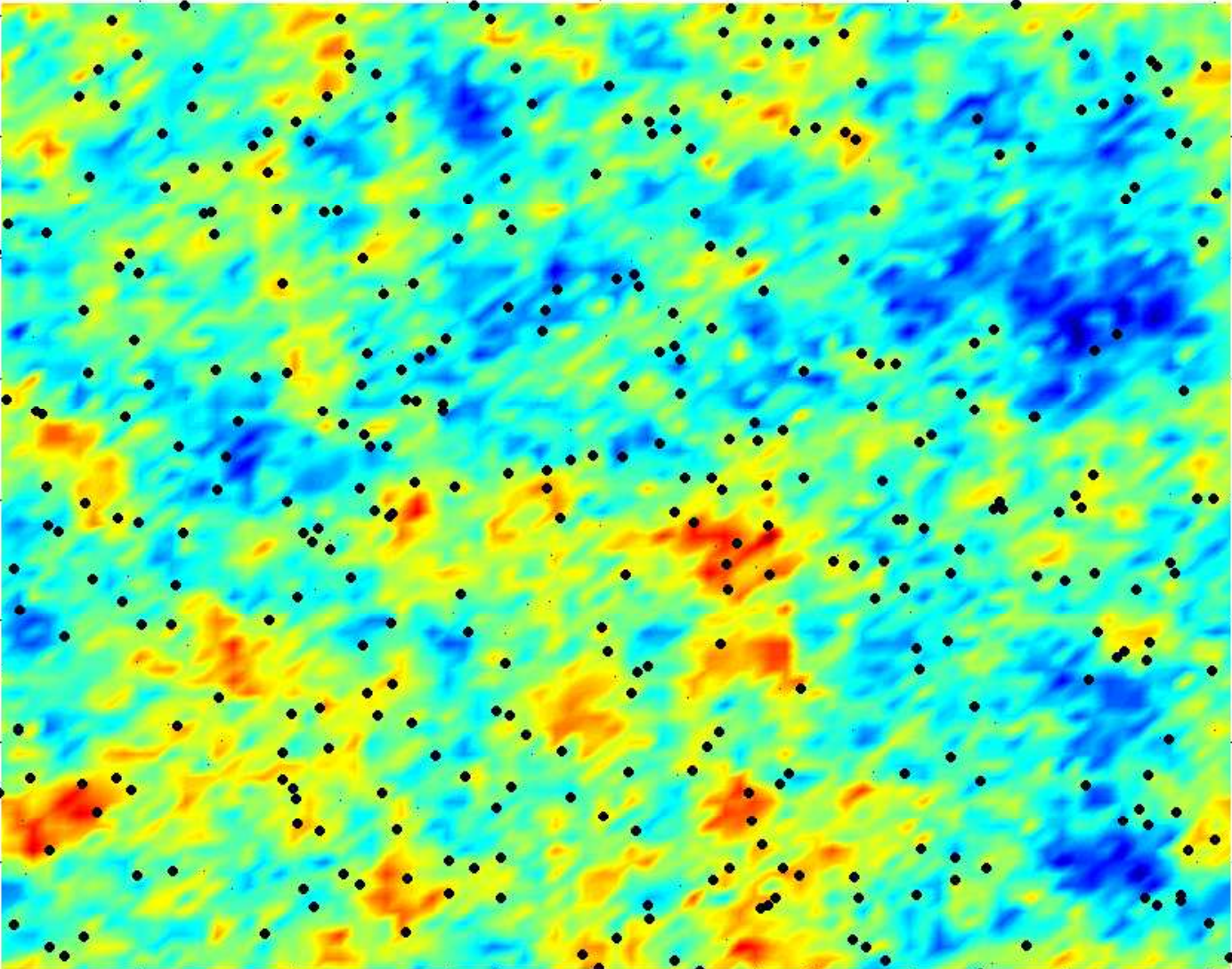}}
\,\,\subfloat[$\gamma^2=1$]{\includegraphics[width=0.44\linewidth]{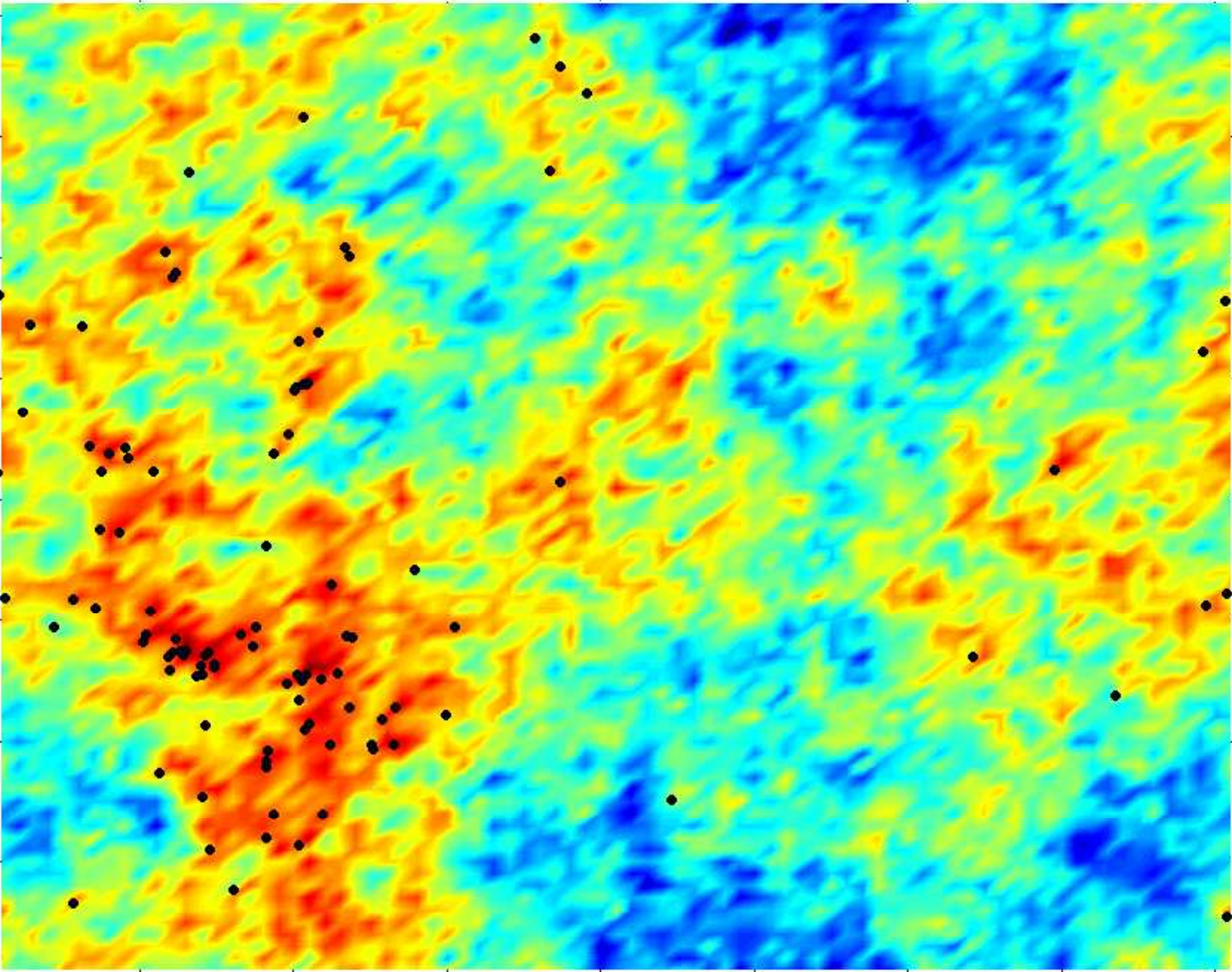}}\,\,
\subfloat[$\gamma^2=3.6$]{\includegraphics[width=0.44\linewidth]{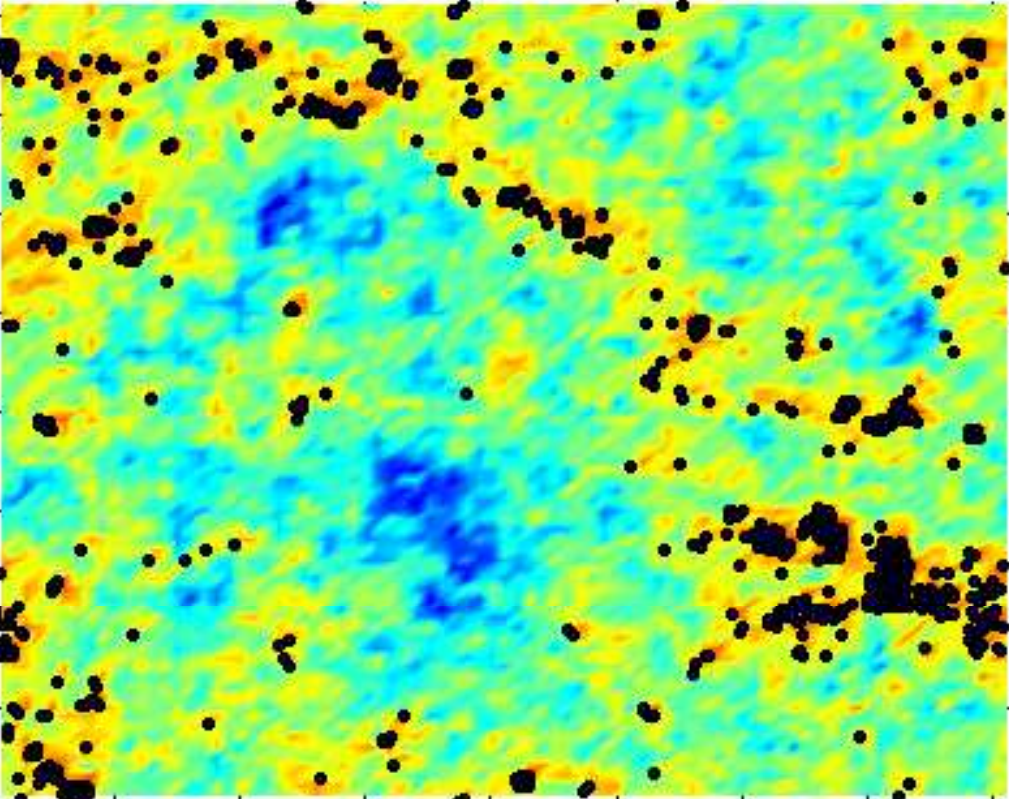}}
\caption{Influence of $\gamma$ on the spatial localization of the atoms of the dual measure.}
\label{figdual2}
\end{figure}

Let us further mention that the construction of atomic Gaussian multiplicative chaos carried out in \cite{Rnew4} is more general than the situation explained here: it is valid in any dimension $d$ for log-correlated Gaussian fields $X$, and also  for all possible values of $\alpha\in]0,1[$ and $ \bar{\gamma}^2>2d$, meaning that we do not impose $\alpha=\frac{2d}{\bar{\gamma}^2}$ ($=\frac{4}{\bar{\gamma}^2}$ in dimension $d=2$). Also, a rigorous proof of the KPZ formula for dual measures appears in \cite{Rnew4}. Special care must be taken here to handle purely atomic measures.

In Figure  \ref{figdual2}, we simulate for different value of $\gamma$  a few atoms of the dual measure (the biggest ones). The colored background stands for the height profile of the associated sub-critical measure $M$ plotted with a logarithmic intensity scale: red for areas with large mass and blue for areas with small mass. Localization of atoms is plotted in black. The larger $\gamma$ is, the more localized on  areas with large potential the atoms are.

\subsubsection*{Frozen phase}

The second part of the super-critical regime consists in throwing a stable point process over a landscape described by   the critical measure $M'$ of subsection \ref{sectheory:deriv}:
\begin{enumerate}
\item sample  the critical measure $M'$ 
$$M'(dx)= \int_A\bigl(\gamma\,\E[X^2(x)]-X(x)\bigr)e^{\gamma X(x)-\frac{\gamma^2}{2}\E[X^2(x)]}\,dx \quad \text{with }\gamma=\sqrt{2d},$$
\item sample a point process $N'_\alpha$ whose law, conditionally to $M'$, is that of  an independently scattered random measure   characterized by $$\forall q\geq 0,\quad \E[e^{-qN'_\alpha(A)}|M]=e^{-q^\alpha M'(A)}. $$\end{enumerate}
This phase is called {\it frozen} due to a linearization of the free energy involved in models converging towards these measures. Nevertheless, the term {\it frozen} may also be understood by the fact that the landscape on which points are thrown is "frozen" and matches the critical measure $M'$: only the height of the atoms varies with the parameter $\alpha$. Furthermore, these family of measures is conjectured to be involved in the glassy phase and freezing phenomena observed in log-correlated random potentials. The reader may consult \cite{CarDou,Fyo,rosso} for an account of physics motivations and results, \cite{BKNSW,Rnew3,madaule,Webb} for rigorous results in the case of discrete models and \cite{Rnew4,Rnew7} for precise conjectures in the context of Gaussian multiplicative chaos theory. More precisely, the glassy phase of log-correlated Gaussian potentials is concerned with the renormalization of measures beyond the critical value $\gamma^2>2d$. For $\gamma^2>2d$, consider the
measure:
$$M_{n,\gamma}(dx)=e^{\gamma X_n(x)-\frac{\gamma^2}{2}\E[X_n(x)^2]}\,dx,$$ where $(X_n)_n$ is your favorite cutoff approximation of the Gaussian distribution $X$ with covariance kernel of the type \eqref{Kintrorev}.
The limiting measure, as $n\to \infty$, vanishes as shown by Theorem \ref{th:necess2}. Therefore, it is natural to wonder how to renormalize the sequence of measures $(M_{n,\gamma})_n$ in order to get a non trivial limiting measure. As pointed out in \cite{Rnew7}, the sequence
\begin{align*}
\big( c_n^{\frac{3\gamma}{2\sqrt{2d}}}e^{c_n\big(\frac{\gamma }{\sqrt{2}}-\sqrt{d}\big)^2}M_{n,\gamma}(dx) \big ) _{t \geq 0}
\end{align*}
is tight and every converging subsequence is non trivial, where $c_n={\rm Var}(X_n)$. This argument is based on the results in \cite{madaule,Webb}. Let us stress that the definition of $c_n$ is clear when $X_n$ is stationary. If $X_n$ is not stationary, there is usually a clear  candidate for the definition of $c_n$, since the behaviour of the variance usually does not depend too much on the spatial localization, like in the case of the GFF. Of course, if one chooses the cutoff approximations in a very bad way, it may happen that the renormalization constant $c_n$ does not straightforwardly make sense. We do not discuss here in further details these technical considerations. 
 Based on heuristics on star scale invariance,  it is conjectured in \cite{Rnew4,Rnew7} that this renormalized sequence actually admits only one possible limit:
 \begin{conjecture}\label{conjfroz}
 Prove that:
\begin{equation}\label{renormsurcrit}
c_n^{\frac{3\gamma}{2\sqrt{2d}}}e^{c_n\big(\frac{\gamma }{\sqrt{2}}-\sqrt{d}\big)^2}M_{n,\gamma}(dx)
\stackrel{law}{\to} d_\gamma N_\alpha(dx),\quad \text{as }t\to \infty
\end{equation}
where $d_\gamma$ is a positive constant depending on $\gamma$ and the law of the random measure $N_\alpha$ is that described above with $\alpha=\frac{\sqrt{2 d}}{\gamma}$.
\end{conjecture}
In particular, physicists are interested in the behaviour of the Gibbs measure associated to $M_{n,\gamma}(dx)$ on a ball $B$. It is the measure renormalized by its total mass:
$$
G^\gamma_n(dx)=\frac{M_{n,\gamma}(dx)}{M_{n,\gamma}(B)}.$$
From \eqref{renormsurcrit}, we deduce
\begin{equation}
G^\gamma_t(dx)\stackrel{law}{\to} \frac{N_\alpha(dx)}{N_\alpha(B)},\quad \text{as }n\to \infty.
\end{equation}
The size reordered atoms of the latter object form a Poisson-Dirichlet process as conjectured by physicists \cite{CarDou}. This Poisson-Dirichlet approach  has recently made important progress  \cite{arguin}, which is presently the more accurate mathematical result to describe the glassy phase of log-correlated random potentials.  A further step is to prove Conjecture \ref{conjfroz}  as it offers a more complete picture of the underlying phenomena than the Poisson-Dirichlet approach. The main reason is that it makes precise the spatial localization of the atoms together with their heights whereas the Poisson-Dirichlet approach only focuses on the heights of the atoms. We stress that this conjecture has been proved in the context of lognormal Mandelbrot's cascades in \cite{BKNSW,Webb} and in the context of branching random walks in \cite{Rnew3}, based on results appearing in \cite{madaule}.

\begin{figure}[h]
\centering
\subfloat[$\alpha=0.2$]{\includegraphics[width=0.44\linewidth]{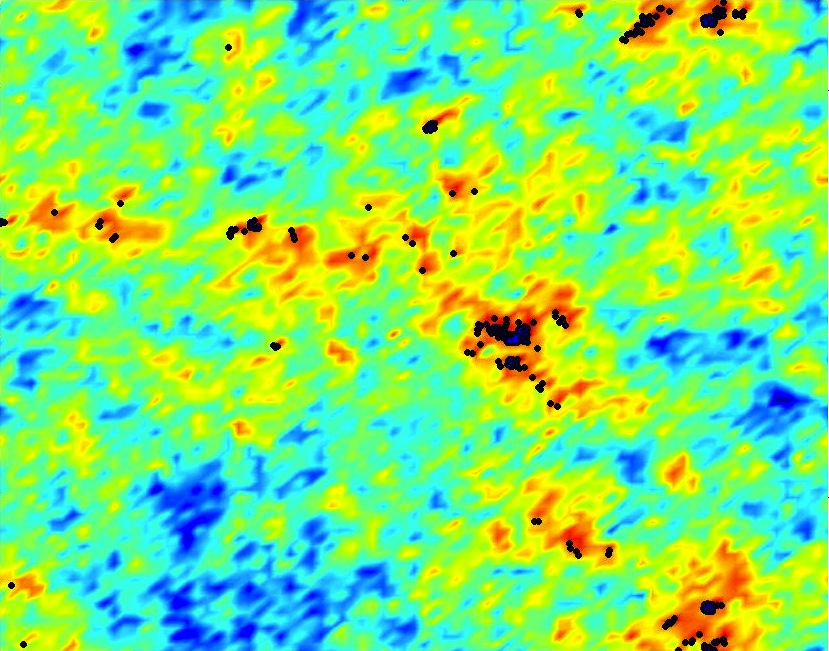}}
\,\,\subfloat[$\alpha=0.5$]{\includegraphics[width=0.44\linewidth]{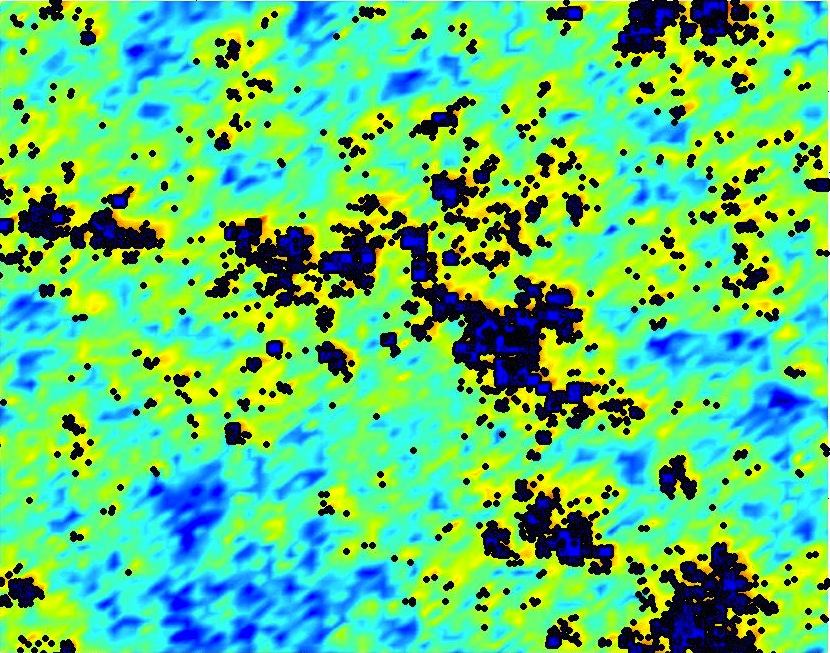}}\,\,
\caption{Localization of atoms in the frozen phase for different values of $\alpha$.}
\label{figfrozen}
\end{figure}

We further stress that the same remarks as in the dual phase remain valid here concerning an almost sure construction of the frozen phase: though the techniques developed in \cite{Rnew4} should apply, this remains to be done rigorously.

In Figure  \ref{figfrozen}, we simulate for different value of $\gamma$  a few atoms of the frozen measure $N'_\alpha$ (the biggest ones). The colored background stands for the height profile of the derivative measure $M'$ plotted with a logarithmic intensity scale: red for areas with large mass and blue for areas with small mass. Localization of atoms is plotted in black. The reader may observe the strong clustering appearing here in comparison with Figure \ref{figdual2}.
\subsection{Conjectures in  connection with lognormal  star scale invariance}\label{sub.conj:star}
Several questions remain open about star scale invariance as mentioned in \cite{Rnew4,Rnew7}, inspired by the discrete multiplicative cascades case  \cite{durrett}. It is conjectured that all the non-trivial ergodic lognormal star scale invariant random measures (actually, we need to impose a stronger dependence decay than ergodicity) belong to one of the families listed below, up to a constant multiplicative  factor. At least, we ask for ergodicity to get rid of the irrelevant random factor $Y$ of Theorem \ref{th:star}. First there must be  an $\alpha\in]0,1]$ such that
$$\E[e^{\alpha \omega_{\varepsilon}(r)}]=\varepsilon^{ d}.$$
Assuming  this, it is proved in (see \cite{Rnew1,Rnew8}) that the Gaussian process $\alpha \omega_{e^{-t}}$ can be rewritten as
$$\omega_{e^{-t}}(x)-\E[\omega_{e^{-t}}(x)]=\frac{\gamma}{\alpha}X_t(x) $$ where  $\gamma^2\leq 2d$ and  $X_t$ is a centered stationary  Gaussian field with covariance structure  given by:   
\begin{equation}\label{theory:cov}
{\rm Cov}\big(X_t(0),X_t(x)\big)= \int_1^{e^t}\frac{k(ux)}{u}\,du
\end{equation}
for some continuous covariance kernel $k$ with $k(0)=1$. Four situations may then occur, each attached with structurally different types of solutions (situations 2,3,4 are conjectures):

\begin{conjecture}\label{conj:star}{\bf Solutions to the $\star$-equation}
\begin{enumerate}
\item If $\alpha=1$ and $\gamma^2<2d$ then the law of the solution $M$ is the standard Gaussian multiplicative chaos, as stated in Theorem \ref{th:star}.
\item If  $\alpha=1$ and $\gamma^2=2d$, then the law of the solution $M$ is that of the derivative martingale $M'$ described in section \ref{sectheory:deriv}.
\item If $\alpha<1$ and $\gamma^2<2d$, then $M$ is an atomic Gaussian multiplicative chaos of the dual phase as described in subsection \ref{sectheory:atomic}.
\item If $\alpha<1$ and $\gamma^2=2d$, then $M$ is an atomic Gaussian multiplicative chaos of the frozen phase.
\end{enumerate}
\end{conjecture}

\subsection{About the maximum of log-correlated Gaussian fields and discrete GFF}
 Another active field of research concerning log-correlated Gaussian fields is the study of their maximum, say over a unit square. More precisely, given a cutoff approximation $(X_n)_n$ of a centered Gaussian
distribution $X$ with logarithmic covariance kernel of the form \eqref{euclid:K}, the question is to know how to renormalize the quantity
\begin{equation}\label{def:max}
\sup_{x \in [0,1]^d} X_n(x)
\end{equation}
in order to get a non trivial limit. Let us set $$c_n={\rm Var}(X_n)$$
with the same discussion about the definition of $c_n$ as in the "frozen part" of section \ref{sectheory:atomic}.
 It is readily seen that the quantity \eqref{def:max} goes to $0$ as $n$ goes to $\infty$, and actually, Kahane's result about the criticality of the value $\gamma^2=2d$ for Gaussian multiplicative chaos already tells you the first order term of this quantity, i.e.
 $$\sup_{x \in [0,1]^d} X_n(x)-\sqrt{2d}c_n=o(c_n), \quad \text{as }n\to\infty.$$
 Technically much more involved, you may even find the second term in the asymptotic expansion of this quantity: more precisely, the family 
  $$\big(\sup_{x \in [0,1]^d} X_n(x)- \sqrt{2d} c_n +\frac{3}{2\sqrt{2d}} \ln c_n \big)_{n \geq 0}$$ 
 converges in law as $n\to\infty$ and the limiting law in non trivial. This is proved  in \cite{bramson} (to be precise, the case of the discrete GFF is treated in dimension  $d=2$ in \cite{bramson} but this result is very likely to hold in generality). By analogy with the branching random walk case (\cite{aidekon}), it is conjectured in \cite{Rnew7} :
\begin{conjecture}\label{conj:max}
 \begin{equation*}
 \sup_{x \in [0,1]^d} X_n(x)- \sqrt{2d} c_n +\frac{3}{2\sqrt{2d}} \ln c_n \to G_d,\quad \text{ in law as }t\to \infty
 \end{equation*}
  where the distribution of $G_d$ is given in terms of the distribution  of the derivative chaos $M'([0,1]^d)$ of subsection \ref{sectheory:deriv}. More precisely, there exists some constant $c>0$ such that:
  \begin{equation}\label{eqonmax}
  \E[e^{-qG_d}]=  \frac{1}{c^q}  \Gamma(1+\frac{q}{\sqrt{2d}})  \E\big[{\big( M'([0,1]^d) \big) ^{-\frac{q}{\sqrt{2d}}}  }\big]
  \end{equation}
  where $\Gamma$ is the standard $\Gamma$-function.
 \end{conjecture}
 
The discussion leading to this conjecture in \cite{Rnew7} is related to star scale invariance as the limiting law $G_d$ can be related to the limit at zero temperature of the frozen phase. It is natural to wonder if similar results should hold for the discrete GFF in a domain $D$ on the vertices of a graph with mesh size going to $0$, thus giving a precise candidate for the limiting law appearing in \cite{bramson}. In fact, in view of the convergence of the discrete Liouville measure at criticality on isoradial graphs (see subsection \ref{cvdlm}), it is clear that the limiting distribution in \cite{bramson} can be nothing but that described in conjecture \ref{conj:max}, where the derivative chaos $M'$ is that related to the (continuous) GFF in $D$ constructed in \cite{Rnew7,Rnew12}. Of course, we do not take into account issues that may arise due to boundary fluctuations: meaning that we assume that the unit square is located far from the boundary of the domain under consideration.

\subsection{Matrix-valued Gaussian multiplicative chaos}\label{matrix}
The main motivation of the Kolmogorov K41 theory (\cite{Kol41})in fully developed turbulence and it's extensions (\cite{cf:Kol}) is to define a realistic statistical theory of an incompressible, homogeneous, isotropic and fully developed turbulent flow (see for example \cite{FoiMan01,cf:Fr}). This ambitious program consists in defining a probabilistic model for the velocity field  which satisfies the main statistical signatures observed experimentally, such as the mean energy transfer towards the small scales and the intermittency (or multifractal) phenomenon (\cite{Kol41}). Ideally, one looks for a field as close as possible to an invariant measure of the equations of motion. In  \cite{CheRob10}, the authors propose a probabilistic construction of such a velocity field. Their construction, which requires a limiting procedure, is mathematically non rigorous and is based on the short time dynamics of the Euler flow, as well as further multifractal considerations. One of the key step of this construction is the introduction of the exponential of an isotropic trace-free matrix whose entries are Gaussian variables with logarithmic correlations. So, this gives rise to the issue of constructing a theory of multiplicative chaos for  Gaussian symmetric   isotropic matrices, which has been studied in \cite{Rnew5}.

The setup is the following. Consider a log-correlated field $(X_x)_x$ of centered Gaussian symmetric matrices and apply a cut-off to the correlation kernel to regularize the logarithmic singularity (in the spirit of \eqref{eq:deformellen}) to obtain a field $(X^n_x)_x$ of centered Gaussian symmetric matrices. The main purpose is to define a matrix-valued measure of the type
\begin{equation}\label{matrix}
M(dx)=\lim_{n\to\infty}\int e^{X^n_x-c_n}\,dx
\end{equation}
 where $c_n$ is a suitable renormalization sequence. This procedure is in essence similar to the scalar case. Nevertheless, the context is highly non-commutative so that finding a proper renormalizing sequence $c_n$ is, in general, fairly tricky, to put it mildly. Of course, things can be much simplified by assuming that the field $X^n$ takes values in a commutative sub-algebra of matrices but there is a big loss of generality in that case and furthermore, this does not fit to the fields observed in $3d$-turbulence: there are some experimental evidences that  the field must be symmetric and isotropic.  Isotropy turns out to be a crucial advantage from the theoretical angle. Indeed, isotropic Gaussian matrices are invariant under the action of the orthogonal group so that the expectation of the matrix exponential $e^{X^n_x}$ must  be proportional to the identity matrix, making expression \eqref{matrix} tractable: the exact computation of $c_n$ is possible thanks to well established formula for eigenvalues of isotropic Gaussian matrices. 
 
 The convergence of the matrix-valued integrals $\int e^{X^n_x-c_n}\,dx$ is not straightforward  since  the martingale property is lost as soon as the field is not scalar. So $L^2$ computations are used in \cite{Rnew5} to prove the convergence. It results that the method is not optimal: it does not give necessary and sufficient condition for convergence. Interestingly, the computations made in \cite{Rnew5} suggest  a logarithmic correction to the power-law spectrum, which does not appear in the scalar case. This seems to be related to non-commutativity. To sum up, there are several intriguing novelties appearing with matrix-valued multiplicative chaos, which deserved to be explored. This research field is widely open.

\subsection{Beyond Gaussianity}
The question of $\star$-scale invariance may be raised in quite a more general framework than the lognormal one:
\begin{definition}\label{levystar}
A stationary random measure $M$ on $\R^d$ is said to be $\star$-scale invariant if  for all $0<\epsilon \leq 1$, $M$ obeys the cascading rule
\begin{equation}\label{eq:levystar}
 \big(M(A)\big)_{A\in\mathcal{B}(\R^d)}\stackrel{law}{=}\big(\int_Ae^{\omega_{\epsilon}(r)}M^{\epsilon}(dr)\big)_{A\in\mathcal{B}(\R^d)}
 \end{equation}
where $\omega_{\epsilon}$ is a stochastically continuous stationary process and $M^{\epsilon}$ is a random measure independent from $\omega_{\epsilon}$ satisfying the relation
$$\big(M^{\epsilon}( A)\big)_{A\in\mathcal{B}(\R^d)}\stackrel{law}{=} \big(M(\frac{A}{\epsilon})\big)_{A\in\mathcal{B}(\R^d)}.$$
\end{definition} 

This is  a quite general statement since nothing is required about the nature of the randomness of the process $\omega_{\epsilon}$. In \cite{Rnew8}, it is proved  that such a scaling relation entails infinite divisibility of the process $\omega_{\epsilon}$ provided that it possesses enough exponential moments.
This suggests the connection with random measures defined as
\begin{equation}\label{theory:levychaos}
Q(dx)=\int e^{X_x}\,dx
\end{equation}
 where $X$ is an infinitely divisible process, suitably normalized, with "logarithmic dependence". 
One is therefore led to considering multiplicative chaos with respect to infinitely divisible random processes. There has been several works in this direction \cite{cf:BaMu,Bar,Fan,Rnew9}, each of which focuses on a specific situation  that is of interest for some given property. Multifractal random measures, not necessarily lognormal, were first mathematically introduced in \cite{Bar}: the construction is geometric and produces a specific example of star scale invariant random measure in the sense of \eqref{eq:levystar}. Many tools that can be used to study measures of the type \eqref{theory:levychaos} are also developed. In \cite{cf:BaMu}, the authors adapt Barral and Mandelbrot's construction \cite{Bar} in order to construct $1$-dimensional exact stochastic scale invariant random measures (non necessarily lognormal). This construction is generalized in \cite{Rnew9} to overcome the problem of dimension and to construct  exact stochastic scale invariant random measures in any dimension. The case where $X$ is a stable infinitely divisible random field is investigated in \cite{Fan}.   In \cite{Rnew8}, the authors construct a unified theory containing all the possible solutions to \eqref{eq:levystar}. Based on a general decomposition for stochastically continuous infinitely divisible processes \cite{maru}, they define a generalized multiplicative chaos that gives sense to random measures of the form \eqref{theory:levychaos} where $X$ is a   general infinitely divisible process with "logarithmic dependence". Their main purpose is  to characterize all the solutions to \eqref{eq:levystar} possessing a moment of order $1+\delta$ for some $\delta>0$. Once again, the continuous nature of the problem imposes strong constraints on the structure of the field $X$ in the spirit of the lognormal case \eqref{structmulti}.

The infinitely divisible case offers several perspectives too: as soon as a critical point exists (the existence depends on the nature of the process $\omega_\epsilon$) the same picture as Conjecture \ref{conj:star} can be drawn, giving rise to corresponding notions/issues of duality, freezing, statistics of extreme values...  But the situation is technically much more involved and most of the picture remains conjectural.



\appendix

\section{Proof of lemma \ref{SmoothGFF}}
Fix $\delta>0$ and $A>0$. We have the following for $x,y \in D^{(\delta)}$:
\begin{align*}
 & E[ (X_\epsilon(x) -X_\epsilon(y))^2   ]  \\
 & =\pi\int_{\epsilon^2}^{\infty} p_D(t,x,x)\,dt+\pi\int_{\epsilon^2}^{\infty} p_D(t,y,y)\,dt-2\pi\int_{\epsilon^2}^{\infty} p_D(t,x,y)\,dt \\
&  = \pi\int_{\epsilon^2}^{\infty} p_D(t,x,x)-  p_D(t,x,y)  \,dt+ \pi\int_{\epsilon^2}^{\infty}  p_D(t,y,y)-  p_D(t,x,y)  \,dt   \\
&  \leq \pi\int_{\epsilon^2}^{1} p_D(t,x,x)-  p_D(t,x,y)  \,dt+ \pi\int_{\epsilon^2}^{1}  p_D(t,y,y)-  p_D(t,x,y)  \,dt  +C |y-x| \\
\end{align*}
where $C$ is some absolute constant (depending on $\delta$) . Recall the following expression for $p_D$: if $(B_s)_{s \geq 0}$ is a Brownian motion starting from $0$, we get:
\begin{equation*}
p_D(t,x,y)= P(  \forall s \leq t, \; B_s-\frac{s}{t}B_t+x+\frac{s}{t}(y-x) \in D    )  \frac{e^{-\frac{|y-x|^2}{2 t}}}{2 \pi t}
\end{equation*}
For the sake of simplicity, we will take $D=B(0,1)$ the Euclidean ball of radius $1$ and $x=0$. For $|y| \leq A \epsilon$, we get:
\begin{align*}
 &  \pi\int_{\epsilon^2}^{1} p_D(t,x,x)-  p_D(t,x,y)  \,dt \\
 & \leq \pi\int_{\epsilon^2}^{1}\frac{1}{2 \pi t} \left ( P(  \underset{s \leq t}{\sup} \; |B_s-\frac{s}{t}B_t| \leq 1    )-   P(   \underset{s \leq t}{\sup} \; |B_s-\frac{s}{t}B_t+\frac{s}{t}y| \leq 1    ) \right )   \,dt \\
 & +  \pi\int_{\epsilon^2}^{1} P(  \underset{s \leq t}{\sup} \; |B_s-\frac{s}{t}B_t+ \frac{s}{t}y| \leq 1    )  | \frac{1}{2 \pi t} -\frac{e^{-\frac{|y|^2}{2 t}}}{2 \pi t}     |  \,dt \\
& \leq \pi\int_{\epsilon^2}^{1}\frac{1}{2 \pi t}   \left ( P(  \underset{s \leq t}{\sup} \; |B_s-\frac{s}{t}B_t| \leq 1    )-   P(   \underset{s \leq t}{\sup} \; |B_s-\frac{s}{t}B_t| \leq 1-|y|    ) \right )     \,dt \\
 & +  \pi\int_{\epsilon^2}^{1}  | \frac{1}{2 \pi t} -\frac{e^{-\frac{|y|^2}{2 t}}}{2 \pi t}     |  \,dt \\
\end{align*}
Now we have $\int_{\epsilon^2}^{1}  | \frac{1}{2 \pi t} -\frac{e^{-\frac{|y|^2}{2 t}}}{2 \pi t}     |  \,dt \leq C \frac{|y|^2}{\epsilon^2}  $. We know turn to the second term. We consider the random variable $X= \underset{s \leq 1}{\sup} \; |B_s- s B_1|$. We have by scaling of Brownian motion (recall also that $\frac{|y|}{\sqrt{t}} \leq A$ for all $t \geq \epsilon^2$):
\begin{align*}
 & \leq \pi\int_{\epsilon^2}^{1}\frac{1}{2 \pi t}   \left ( P(  \underset{s \leq t}{\sup} \; |B_s-\frac{s}{t}B_t| \leq 1    )-   P(   \underset{s \leq t}{\sup} \; |B_s-\frac{s}{t}B_t| \leq 1-|y|    ) \right )     \,dt \\
 & =  \pi\int_{\epsilon^2}^{1}\frac{1}{2 \pi t}    P( \frac{1}{\sqrt{t}} -\frac{|y|}{\sqrt{t}}  \leq  X \leq \frac{1}{\sqrt{t}}    )     \,dt   \\
& \leq  C_A  |y|  \int_{\epsilon^2}^{1}\frac{1}{t^{3/2}}     \,dt  \\
& \leq C_A  \frac{|y|}{\epsilon}.
\end{align*}
In conclusion, we can find a $C_A>0$ such that for all $x,y \in D^{(\delta)}$ with $|y-x| \leq A \epsilon $:
\begin{equation*}
E[ (X_\epsilon(x) -X_\epsilon(y))^2   ] \leq C_A  \frac{|y-x|}{\epsilon}.
\end{equation*}


\section{Proof of convergence of Theorem \ref{DLQG}}\label{discLQG}
Let us carry out the proof in the case of the square lattice for simplicity. The reader may then adapt the proof to the general case using subsubsection \ref{ch-smir} below.

Let $D$ be the standard square $[0,1]^2$. We let $D_{\epsilon}$ stand for the $\epsilon^{-1}$-dilation of $D$, i.e. $D_{\epsilon}=\frac{D}{\epsilon}$ and $L_\epsilon$ for the lattice $\epsilon \Z^2\cap D$. For simplicity, let us suppose that $\epsilon=\frac{1}{2^N}$ where $N$ is an integer. We partition $D$ into the squares 
\begin{equation*}
\Delta_{\epsilon,i,j}:= [(i-\frac{1}{2}) \epsilon, (i+\frac{1}{2})\epsilon]  \times [(j-\frac{1}{2}) \epsilon, (j+\frac{1}{2})\epsilon].
\end{equation*}
We fix $\delta \in ]0,1[$ and we will work in $D^{(\delta)}= \lbrace x \in D; \; d(x,\partial D) \geq \delta \rbrace$ the points in $D$ which are at distance $\delta$ from the boundary. 

We introduce, for all $t \in [0,1]$, the Gaussian process $Z_{\epsilon}(t,x)=\sqrt{t}X_{\epsilon}(x)+\sqrt{1-t}\bar{X}_{\epsilon}(x)$  where $X_{\epsilon}(x)$ is the discrete GFF on $L_\epsilon$ (that we extend to the square $D$ by setting it constant in each $\Delta_{\epsilon,i,j}$) and $\bar{X}_{\epsilon}(x)$ is an independent Gaussian process with covariance $\pi \int_{\epsilon}^\infty p_D(t,x,y)\,dt$. In view of the techniques of generalized Gaussian multiplicative chaos \cite{cf:RoVa}, one must show the following three points:
\begin{enumerate} 
\item
For all $A>0$,
\begin{equation*}
 \underset{\underset{|x-y| \leq A \epsilon}{x,y \in D^{(\delta)},}}{\sup}  |  E[X_{\epsilon}(x)X_{\epsilon}(y)]- \pi \int_{\epsilon^2}^{\infty}p_D(t,x,y)dt    | \leq \bar{C}_A.
\end{equation*}
where $\bar{C}_A>0$ is some constant independent from $\epsilon$.

\item
If we set:
\begin{equation*}
C_A = \underset{\epsilon \to 0} {\overline{\lim}}     \underset{\underset{|x-y| \geq A \epsilon}{x,y \in D^{(\delta)},}}{\sup}  |  E[X_{\epsilon}(x)X_{\epsilon}(y)]  - \pi \int_{\epsilon^2}^{\infty}p_D(t,x,y)dt |
\end{equation*}
then $C_A$ goes to $0$ as $A$ goes to infinity.
\item 

For all $\alpha<1$, we have the following convergence uniformly in $t \in [0,1]$:
\begin{equation}\label{modulus1}
\underset{\epsilon \to 0}{\overline{\lim}}  \: E[  \left(    \sup_{{\delta} {\epsilon} \leq i,j \leq \frac{1-\delta}{\epsilon}}  \int_{\Delta_{\epsilon,i,j}}   e^{ \gamma Z_{\epsilon}(t,x)-\frac{\gamma^2}{2} E[ Z_{\epsilon}(t,x)  ]  }    dx   \right)^{\alpha}   ]=0
\end{equation}
in the subcritical case $\gamma<2$ or 
\begin{equation}\label{modulus2}
\underset{\epsilon \to 0}{\overline{\lim}}  \: E[  \left(  \sqrt{\ln \frac{1}{\epsilon}}   \sup_{{\delta} {\epsilon} \leq i,j \leq \frac{1-\delta}{\epsilon}}  \int_{\Delta_{\epsilon,i,j}}   e^{ 2 Z_{\epsilon}(t,x)-2 E[ Z_{\epsilon}(t,x)  ]  }    dx   \right)^{\alpha}   ]=0
\end{equation}
in the critical case $\gamma=2$.
\end{enumerate}

\subsubsection*{Proof of item 1,2}
Proof of item 1 is straightforward.

The covariance of $X_{\epsilon}$ is given by:
\begin{equation*}
\forall x,y \in L_\epsilon,\quad E[X_{\epsilon}(x)X_{\epsilon}(y)]= E^{\lfloor \frac{x}{\epsilon} \rfloor }  [  \sum_{i=1}^{\tau_{\partial D_{\epsilon}}} 1_{S_i= \lfloor \frac{y}{\epsilon} \rfloor}   ]
\end{equation*}
 where $(S_{i})_{i \geq 1}$ is the simple random walk in $\Z^2$ starting from $\lfloor \frac{x}{\epsilon} \rfloor$ and $\tau_{\partial D_{\epsilon}}$ stands for the first exit time of the random walk out of $D_\epsilon$. By \cite{Lawler}, we have the following expression:
 \begin{equation*}
 E[X_{\epsilon}(x)X_{\epsilon}(y)]  = \sum_{z \in \partial D_{\epsilon}} P^{\lfloor \frac{x}{\epsilon} \rfloor } ( S_{ \tau_{\partial D_{\epsilon}  }}= z  ) a(z-\lfloor \frac{y}{\epsilon} \rfloor)
 - a(\lfloor \frac{x}{\epsilon} \rfloor-\lfloor \frac{y}{\epsilon} \rfloor) 
 \end{equation*}
 where $x,y \in L_\epsilon$ and  $a(x)$ is the potential kernel for the simple random walk:
 \begin{equation*}
 a(x)=\frac{2}{\pi} \ln |x| +\frac{2 \bar{\gamma}+ \ln 8}{\pi}+ O(\frac{1}{|x|^2})
 \end{equation*}
 and $\bar{\gamma}$ is the Euler constant. Now, recall the following expression for the Green function of Brownian motion in $D$ killed upon touching the boundary:
 
 \begin{equation*}
 G_D(x,y)= \int_{\partial D} p_{D}(x,z)  \ln \frac{|z-y|}{|x-y|}  dz
 \end{equation*}
 where $p_D$ is the Poisson kernel. 
 Let $A>0$. We set $p_\epsilon(x,z)= P^{\lfloor \frac{x}{\epsilon} \rfloor } ( S_{ \tau_{\partial D_{\epsilon}  }}= \lfloor \frac{z}{\epsilon} \rfloor )/ \epsilon$, then we get for $|x-y| \geq \epsilon A$ (see \cite{Lawler}) and $x,y \in [\delta,1-\delta]^2$: 
 \begin{align*}
 E[X_{\epsilon}(x)X_{\epsilon}(y)]  & = \sum_{z \in \partial D_{\epsilon}} P^{\lfloor \frac{x}{\epsilon} \rfloor } ( S_{ \tau_{\partial D_{\epsilon}  }}= z  ) \ln | z-\lfloor \frac{y}{\epsilon} \rfloor |- \ln | \lfloor \frac{x}{\epsilon} \rfloor  -\lfloor \frac{y}{\epsilon} \rfloor |  +O(\frac{1}{A^2})  \\
& = \sum_{z \in \partial D_{\epsilon}} P^{\lfloor \frac{x}{\epsilon} \rfloor } ( S_{ \tau_{\partial D_{\epsilon}  }}= z  ) \ln | \epsilon( z-\lfloor \frac{y}{\epsilon} \rfloor) |- \ln | \epsilon( \lfloor \frac{x}{\epsilon} \rfloor  -\lfloor \frac{y}{\epsilon} \rfloor )|  +O(\frac{1}{A^2})  \\
& = \int_{\partial D} p_\epsilon(x,z) \ln | \epsilon( \lfloor \frac{z}{\epsilon} \rfloor-\lfloor \frac{y}{\epsilon} \rfloor) |dz - \ln | \epsilon( \lfloor \frac{x}{\epsilon} \rfloor  -\lfloor \frac{y}{\epsilon} \rfloor )|  +O(\frac{1}{A^2})  \\ 
  \end{align*}
By lemma B.1 in \cite{Smir}, we have that  $\sum_{z \in \partial D_{\epsilon}} P^{\lfloor \frac{x}{\epsilon} \rfloor } ( S_{ \tau_{\partial D_{\epsilon}  }}= z  ) \ln | \epsilon( z-\lfloor \frac{y}{\epsilon} \rfloor) |$ converges uniformly to $ \int_{\partial D} p_{D}(x,z)  \ln |z-y| dz$. 
By the uniform convergence above, we get that $C_A$ converges to $0$ as $A \to \infty$.

\subsubsection*{Proof of item 3}
We will work in the critical case $\gamma=2$ as the subcritical case is a straithforward adaptation of the proof in \cite{cf:RoVa}.

We have the following inequality:
\begin{align*}
& E[  \left(  \sqrt{\ln \frac{1}{\epsilon}}   \sup_{{\delta} {\epsilon} \leq i,j \leq \frac{1-\delta}{\epsilon}}  \int_{[i \epsilon, (i+1)\epsilon]  \times [j \epsilon, (j+1)\epsilon]}   e^{2 Z_{\epsilon}(t,x)-2 E[ Z_{\epsilon}(t,x)  ]  }    dx   \right)^{\alpha}   ]   \\
& \leq  E[  \left(  \sqrt{\ln \frac{1}{\epsilon}}  \epsilon^2 \sup_{{\delta} {\epsilon} \leq i,j \leq \frac{1-\delta}{\epsilon}} \sup_{x \in [i \epsilon, (i+1)\epsilon]  \times [j \epsilon, (j+1)\epsilon]}    e^{ 2 Z_{\epsilon}(t,x)-2 E[ Z_{\epsilon}(t,x)  ]  }      \right)^{\alpha}   ]   \\
\end{align*}

Let us fix $t$. If $x \in \Delta_{\epsilon,i,j}$, we decompose the process $Z_{\epsilon}(t,x)= Z_{\epsilon}(t,i,j)+\bar{Z}_{\epsilon}(t,x)$ 
where:
\begin{equation*}
Z_{\epsilon}(t,i,j) = \sqrt{t}X_{\epsilon}(x)+\sqrt{1-t}  \sqrt{\pi} \int_{D\times [\epsilon^2 ,\infty[} \underset{x \in \Delta_{\epsilon,i,j} }{\inf} p_D(\frac{s}{2},x,y)W(dy,ds).
\end{equation*}
and $(\bar{Z}_{\epsilon}(t,x))_{x \in  \Delta_{\epsilon,i,j}}$ is the remaining smooth Gaussian process. Note that the process $\bar{Z}_{\epsilon}(t,x)$ has the following properties:
\begin{itemize}
 \item
 It is continuous on $\Delta_{\epsilon,i,j}$
 \item
There exists some constant $C>0$ independent from $\epsilon,i,j,t$ and such that:
\begin{equation*}
\forall x \in \Delta_{\epsilon,i,j}, \; E[ \bar{Z}_{\epsilon}(t,x)^2 ] \leq C
\end{equation*}
\item
There exists some constant $C>0$ independent from $\epsilon,i,j,t$ and such that:
\begin{equation*}
\forall x \in \Delta_{\epsilon,i,j}, \; 0 \leq E[ Z_{\epsilon}(t,i,j)  \bar{Z}_{\epsilon}(t,x)] \leq C
\end{equation*}
\end{itemize}

We can then introduce the standard $2d$ discrete cascade  $(Z_{i,j})_{1 \leq i, j \leq \frac{1}{\epsilon}}$ independent from $\bar{Z}_{\epsilon}$ (see \cite{Rnew7}). Recall that there exists a constant $C>0$ independent from $\epsilon,i,j,t$ such that for all $(i,j),(i',j') \in [| \frac{\delta}{\epsilon}, \frac{1-\delta}{\epsilon}   |]$:
\begin{equation*}
 E[ Z_{i,j} Z_{i',j'}  ] \leq E [Z_{\epsilon}(t,i,j) Z_{\epsilon}(t,i',j') ]+C
\end{equation*}

From above, we deduce that one can find a fixed random variable $Y$ independent from $Z_{\epsilon}$ of bounded variance (with respect to all variables $\epsilon,i,j,t$) and independent Gaussian processes $(Y_{\epsilon,i,j}(t,x))_{i,j}$ (also independent from $(Z_{i,j})_{1 \leq i, j \leq \frac{1}{\epsilon}}$) such that the following holds:
\begin{itemize}
\item
For each $i,j$, the process $Y_{\epsilon,i,j}(t,x)$ is continuous in $ \Delta_{\epsilon,i,j}$ and there exists some constant $C>0$ independent from $\epsilon,i,j,t$ and such that:
\begin{equation*}
\forall x \in \Delta_{\epsilon,i,j}, \; 0 \leq E[ Y_{\epsilon,i,j}(t,x)^2] \leq C
\end{equation*}

\item
The processes $Z_{\epsilon}+Y$ and $\tilde{Z}_{\epsilon}$ have same variance, where $\tilde{Z}_{\epsilon}=Z_{i,j}+Y_{\epsilon, i,j}$. Thus, we have:
\begin{equation*}
E[(Z_{\epsilon}(t,x)+Y)^2]=E[\tilde{Z}_{\epsilon}(t,x)^2]
\end{equation*}
\item
The processe $Z_{\epsilon}+Y$ is more correlated than the process $\tilde{Z}_{\epsilon}$:
\begin{equation*}
E[(Z_{\epsilon}(t,x)+Y)(Z_{\epsilon}(t,x)+Y)] \geq E[\tilde{Z}_{\epsilon}(t,x)\tilde{Z}_{\epsilon}(t,y)]
\end{equation*}
\end{itemize}
Therefore the process $Z_{\epsilon}(t,x)+Y$ is stochastically dominated by the process  $\tilde{Z}_{\epsilon}(t,x)$ by the standard Slepian lemma (see \cite{Rnew7} for example). In conclusion, we get that:
\begin{align*}
& \leq  E[  \left(  \sqrt{\ln \frac{1}{\epsilon}}  \epsilon^2 \sup_{{\delta} {\epsilon} \leq i,j \leq \frac{1-\delta}{\epsilon}} \sup_{x \in [i \epsilon, (i+1)\epsilon]  \times [j \epsilon, (j+1)\epsilon]}    e^{ 2 Z_{\epsilon}(t,x)-2 E[ Z_{\epsilon}(t,x)  ]  }       \right)^{\alpha}   ]   \\
& \leq  C E[  \left(  \sqrt{\ln \frac{1}{\epsilon}}  \epsilon^2 \sup_{{\delta} {\epsilon} \leq i,j \leq \frac{1-\delta}{\epsilon}}   \bar{Z}_{i,j} e^{ 2 Z_{i,j}-2 E[ Z_{i,j} ]  }   \right)^{\alpha}   ]   \\
\end{align*}
where $\bar{Z}_{i,j}=\sup_{x \in \Delta_{\epsilon,i,j}}    e^{ 2 Y_{\epsilon,i,j}(x) -2 E[ Y_{\epsilon,i,j}(x)^2 ]  } $. Now, it is plain to see that the random measure which gives a mass $ \sqrt{\ln \frac{1}{\epsilon}}  \epsilon^2  \bar{Z}_{i,j} e^{ 2 Z_{i,j}-2 E[ Z_{i,j} ]  } $ to each square $ \Delta_{\epsilon,i,j}$ converges (possibly along subsequences) in law to a measure which is absolutely continuous with respect to the standard canonical cascade measure (which correponds to taking $\bar{Z}_{i,j}=1$). Since this measure has no atoms, we deduce that:
\begin{equation*} 
     E[  \left(  \sqrt{\ln \frac{1}{\epsilon}}  \epsilon^2 \sup_{{\delta} {\epsilon} \leq i,j \leq \frac{1-\delta}{\epsilon}}   \bar{Z}_{i,j} e^{ 2 Z_{i,j}-2 E[ Z_{i,j} ]  }   \right)^{\alpha}   ]  \underset{\epsilon \to 0} {\rightarrow} 0.
\end{equation*}

\subsubsection*{General case}\label{ch-smir}
To adapt the above proof to the general class of isoradial graphs, one just need the following properties of the Green function, which can be found in \cite{chelkak} with the minor difference that we define the Green function as $-2\pi$ times the Green function of \cite{chelkak}. The minus serves to have a positive definite function and the $2\pi$ is the standard normalization in Liouville quantum gravity to have the Green function that asymptotically behaves like a normalized log. In what follows, $G_{\Gamma_n}$ stands for the discrete Green function on $\Gamma_n$ and $G_{\Omega_n}$ for the discrete Green function on $\Omega_n$ with $0$ boundary condition. In particular, from Definition 2.3+Theorem 2.5 in \cite{chelkak}, we have:
\begin{enumerate}
\item $G_{\Gamma_n}(x,x)=\ln\frac{1}{\epsilon_n}+c$ for some explicit constant $c$, which does not depend on relevant quantities. 
\item $\forall x\not =y\in \Gamma_n$, $$G_{\Gamma_n}(x,y)=\ln\frac{1}{|x-y|}+O\Big(\frac{\epsilon_n^2}{|x-y|^2}\Big)$$ uniformly w.r.t. the shape of $\Gamma_n$ and $x\not =y\in \Gamma_n$.
\item $G_{\Omega_n}=G_{\Gamma_n}-G^*_{\Omega_n}$ where $(G^*_{\Omega_n})_n$ is a sequence of functions  that uniformly converges on the closed balls included in $\Omega$ towards the function $G^*$ defined by
$$\triangle G^*(\cdot ,y)=0\text{ on }  \Omega,\quad G^*(\cdot,y)=\ln\frac{1}{|\cdot-y|}\text{ on }\partial \Omega.$$
\item for each closed ball $B\subset \Omega$, there is a constant $C_B>0$ such that $\forall x,y\in B$ and $n\geq 0$
$$\ln\frac{1}{\epsilon_n}-C\leq G_{\Omega_n}\leq \ln\frac{1}{\epsilon_n}+C.$$
\end{enumerate}

\section{Multifractal formalism}
\subsection*{Proof of theorem \ref{analysemulti}}

We will prove theorem \ref{analysemulti} in dimension $d=2$ though the proof can easily be adapted to all dimensions (we will indicate when appropriate how to adapt the proof to all dimensions). The proof that we will provide is robust, i.e. it is enough to show the result for one special Gaussian sequence satisfying the above conditions to deduce the result for all Gaussian sequences by using Kahane's convexity inequalities. Therefore, let us choose the exact scale invariant field $\bar{X}$ with covariance kernel given by:
\begin{equation*}
\forall x,y\in U,\quad \E[\bar{X}(x)\bar{X}(y)]=\ln_+\frac{1}{ |x-y|}.
\end{equation*}
Let us also consider  a white noise decomposition $(\bar{X}_\epsilon)_{\epsilon\in]0,1]}$ of $\bar{X}$ as constructed in \cite{cf:RoVa} (in all dimensions, one can work with the exact stochastic scale invariant kernels introduced in \cite{Rnew9}). In particular, the process $\epsilon \rightarrow \bar{X}_\epsilon$ has independent increments and $\bar{X}_{\epsilon,\epsilon'}:=\bar{X}_{\epsilon}-\bar{X}_{\epsilon'}$ has a correlation cutoff of length $\epsilon'$ (i.e. if the Euclidean distance between two sets $A,B$ is greater than $\epsilon'$ then $(\bar{X}_{\epsilon,\epsilon'}(x))_{x \in A}$ and $(\bar{X}_{\epsilon,\epsilon'}(x))_{x \in B}$ are independent). The correlation structure of $(\bar{X}_\epsilon)_{\epsilon\in]0,1]}$ is given for $\epsilon \in]0,1]$ by:
$$\E[\bar{X}_\epsilon(x)\bar{X}_{\epsilon}(y)]=\left\{
\begin{array}{ll}
0 & \text{ if }|x-y|>1\\
\ln\frac{1}{ |x-y|}& \text{ if } \epsilon \leq |x-y|\leq 1\\
\ln \frac{1}{\epsilon}+2( 1- \frac{|x-y|^{1/2}}{\epsilon^{1/2}}  ) &  \text{ if } |y-x| \leq \epsilon
\end{array}\right..$$
In particular, we have the following exact stochastic scale invariant relation for all $\lambda<1$:
\begin{equation*}
(\bar{X}_{\lambda \epsilon}(\lambda x))_{|x| \leq 1} \overset{(Law)}{=} (\bar{X}_{\epsilon}(x))_{|x| \leq 1}+\Omega_{\frac{1}{\lambda}}
\end{equation*}
where $\Omega_{\frac{1}{\lambda}}$ is a centered Gaussian random variable of variance $\ln \frac{1}{\lambda}$ independent of $\bar{X}_{\epsilon}$.

\subsection*{Preliminary lemma}
Now we begin the proof. Recall that we work in dimension $d=2$.
 \begin{lemma}\label{am.prel}
Let $\gamma\in [0,2[$. If $a\in]0,1]$ satisfies $(2-\frac{\gamma^2}{2})a-\frac{\gamma^2 a^2}{2}>0$ then there exists a constant $C>0$ such that:
\begin{equation} 
\sup_{\epsilon\in]0,1]} \E \Big [ \Big( \int_{[0,\frac{1}{2}]^2}  \frac{e^{  \gamma  \bar{X}_{\epsilon}(v)-\frac{ \gamma^2}{2}\E[ \bar{X}_{\epsilon}(v)^2]}}{ (|v|+ \epsilon)^{\gamma^2}  } dv \Big)^a \Big]   \leq C.
\end{equation}
\end{lemma}

\proof By stochastic scale invariance, we have the following property:
\begin{align*}
  \E  \Big[  \Big(& \int_{[0,\frac{1}{2}]^2}  \frac{e^{  \gamma \bar{X}_{\epsilon}(v)-\frac{ \gamma^2}{2}\E[\bar{X}_{\epsilon}(v)^2]}}{ (|v|+ \epsilon)^{\gamma^2}  } dv  \Big)^a  \Big]  \\
  \leq & \E \Big[  \Big( \int_{[0,\frac{1}{4}]^2}  \frac{e^{  \gamma  \bar{X}_{\epsilon}(v)-\frac{ \gamma^2}{2}\E[ \bar{X}_{\epsilon}(v)^2]}}{ (|v|+ \epsilon)^{\gamma^2}  } dv  \Big)^a  \Big] + \E  \Big[  \Big( \int_{[0,\frac{1}{2}]^2 \setminus [0,\frac{1}{4}]^2}  \frac{e^{  \gamma  \bar{X}_{\epsilon}(v)-\frac{ \gamma^2}{2}\E[ \bar{X}_{\epsilon}(v)^2]}}{ (|v|+ \epsilon)^{\gamma^2}  } dv  \Big)^a  \Big]   \\
  \leq &\frac{1}{2^{a (2-\frac{\gamma^2}{2} )-\frac{\gamma^2 a^2}{2}}}  \E \Big[ \Big( \int_{[0,\frac{1}{2}]^2}  \frac{e^{  \gamma  \bar{X}_{2 \epsilon}(v)-\frac{ \gamma^2}{2}\E[ \bar{X}_{2 \epsilon}(v)^2]}}{ (|v|+ 2 \epsilon)^{\gamma^2}  } dv \Big)^a \Big]   \\
 &+ 4^{a\gamma^2}\sup_{\epsilon} \E  \Big[  \Big( \int_{[0,\frac{1}{2}]^2}  e^{  \gamma  \bar{X}_{\epsilon}(v)-\frac{ \gamma^2}{2}\E[ \bar{X}_{\epsilon}(v)^2]} dv  \Big)^a  \Big]  ,
\end{align*} 
where in the last inequality we have used stochastic scale invariance. By assumption, we have  $\rho= \frac{1}{2^{a (2-\frac{\gamma^2}{2} )-\frac{\gamma^2 a^2}{2}}} <1$. Since $\gamma\in [0,2[$, we also have:
$$C\stackrel{\text{def}}{=}4^{a\gamma^2}\sup_{\epsilon} \E  \Big[  \Big( \int_{[0,\frac{1}{2}]^2}  e^{  \gamma  \bar{X}_{\epsilon}(v)-\frac{ \gamma^2}{2}\E[ \bar{X}_{\epsilon}(v)^2]} dv  \Big)^a  \Big] <+\infty.$$
Set: 
\begin{equation*}
u_n= \E \Big [ \Big( \int_{[0,\frac{1}{2}]^2}  \frac{e^{  \gamma  \bar{X}_{\epsilon_n}(v)-\frac{ \gamma^2}{2}E[ \bar{X}_{\epsilon_n}(v)^2]}}{ (|v|+ \epsilon_n)^{\gamma^2}  } dv \Big )^a \Big] 
\end{equation*}
where $\epsilon_n=\frac{1}{2^n}$. We have just proved that $u_{n+1} \leq \rho u_n+C$ hence the sequence $(u_n)_n$ is bounded, thus giving the result.\qed

\subsection{Proof of Theorem \ref{analysemulti}}
In what follows, the value of $\alpha$ is fixed and given by $\alpha= 2+(\frac{1}{2}-q)\gamma^2$.

\subsubsection*{Lower bound}
\begin{lemma}\label{lemlow}
Fix $\beta>0$. There exist two constants $D>0$ and $\eta>0$ such that for all $r>0$, we have:
\begin{equation*}
\E\Big[   \int_{[0,1]^2}  \ind_{\{M_\gamma (B(x,r))   \geq r^{\alpha- \beta}\}}  M_{q \gamma}(dx)  \Big]  \leq D r^{\eta}.
\end{equation*}
\end{lemma}

\proof We fix $\epsilon, \epsilon'>0$ such that $\epsilon'< r \epsilon$. Let $\Omega_{1/r}$ denote a centered Gaussian random variable of variance $\ln \frac{1}{r}$. We write $\Omega_{1/r}=\sqrt{\ln \frac{1}{r}  }N$ whrere $N$ is a standard Gaussian. By using the Girsanov transform, we get:
\begin{align*}
& \E\Big[  \ind_{\{ \int_{B(x,r)}  e^{  \gamma \bar{X}_{r \epsilon}(u)-\frac{ \gamma^2}{2}\E[\bar{X}_{r \epsilon}(u)^2]} du   \geq r^{\alpha- \beta}\}} e^{ q \gamma \bar{X}_{\epsilon'}(x)-\frac{q^2 \gamma^2}{2}\E[\bar{X}_{\epsilon'}(x)^2]} \Big ]    \\
& =\Pb\Big[e^C \int_{B(x,r)}  \frac{e^{  \bar{\gamma X}_{r \epsilon}(u)-\frac{ \gamma^2}{2}\E[\bar{X}_{r \epsilon}(u)^2]}}{ (|u-x|+r \epsilon)^{q \gamma^2}  } du   \geq r^{\alpha- \beta}\Big ]  \\
& \leq \Pb \Big[ e^Ce^{\gamma \Omega_{1/r}} \int_{B(0,1)}  \frac{e^{  \gamma \bar{X}_{\epsilon}(v)-\frac{ \gamma^2}{2}\E[\bar{X}_{\epsilon}(v)^2]}}{ (|v|+ \epsilon)^{q \gamma^2}  } dv   \geq r^{- \beta}\Big ]  \\
& \leq  \Pb\big[ e^C e^{\gamma \Omega_{1/r}}  \geq r^{- \frac{\beta}{2}}\big] +  \Pb\Big[ \int_{B(0,1)}  \frac{e^{  \gamma \bar{X}_{\epsilon}(v)-\frac{ \gamma^2}{2}\E[\bar{X}_{\epsilon}(v)^2]}}{ (|v|+ \epsilon)^{q \gamma^2}  } dv   \geq r^{- \frac{\beta}{2}} \Big]  \\ 
& \leq  \Pb\big[ N  \geq \frac{\beta}{2\gamma} \sqrt{\ln \frac{1}{r}}-\frac{C}{ \sqrt{\ln \frac{1}{r}}} \big] +  \Pb\Big[ \int_{B(0,1)}  \frac{e^{  \gamma \bar{X}_{\epsilon}(v)-\frac{ \gamma^2}{2}\E[\bar{X}_{\epsilon}(v)^2]}}{ (|v|+ \epsilon)^{q \gamma^2}  } dv   \geq r^{- \frac{\beta}{2}} \Big]  \\ 
& \leq  e^{C\beta\gamma^{-2}}r^{\frac{\beta^2}{8 \gamma^2}}  +  r ^{\frac{a \beta}{2}}  \E \Big [ \Big ( \int_{B(0,1)}  \frac{e^{  \gamma  \bar{X}_{\epsilon}(v)-\frac{ \gamma^2}{2}\E[ \bar{X}_{\epsilon}(v)^2]}}{ (|v|+ \epsilon)^{\gamma^2}  } dv \Big)^a \Big]  \\ 
& \leq  e^{C\beta\gamma^{-2}}r^{\frac{\beta^2}{8 \gamma^2}} + C' r ^{\frac{a \beta}{2}} ,
\end{align*}
where $C$ is the constant appearing in \eqref{am.hyp} and $C'$ comes from Lemma \ref{am.prel}.
Therefore, by integrating this identity with respect to $dx$, we get that:
\begin{equation*}
 \E\Big[ \int_{[0,1]^2} \ind_{ \{\int_{B(x,r)}  e^{  \gamma \bar{X}_{r \epsilon}(u)-\frac{ \gamma^2}{2}E[\bar{X}_{r \epsilon}(u)^2]} du   \geq r^{\alpha- \beta}\}} e^{ q \gamma \bar{X}_{\epsilon'}(x)-\frac{q^2 \gamma^2}{2}\E[\bar{X}_{\epsilon'}(x)^2]} dx \Big ]    \leq  D r^\eta  
\end{equation*}
where $D$ and $\eta$ are independent of $\epsilon, \epsilon'$.
One can take the limit as $\epsilon' \to 0$ and conclude that:
\begin{equation*}
 \E\Big[ \int_{[0,1]^2} 1_{\{ \int_{B(x,r)}  e^{  \gamma \bar{X}_{r \epsilon}(u)-\frac{ \gamma^2}{2}\E[\bar{X}_{r \epsilon}(u)^2]} du   \geq r^{\alpha- \beta}\}} M_{q \gamma}(dx)  \Big]    \leq D  r^\eta. 
\end{equation*}
Now one can use Fatou twice to get that:
\begin{align*}
&  \E[ \int_{[0,1]^2} \ind_{\{M_{\gamma}(B(x,r)) \geq 2 r^{\alpha- \beta}\}}   M_{q \gamma}(dx) \Big ]    \\
& \leq \E\Big[ \int_{[0,1]^2} \liminf_{\epsilon \to 0} \ind_{\{ \int_{B(x,r)}  e^{  \gamma \bar{X}_{r \epsilon}(u)-\frac{ \gamma^2}{2}\E[\bar{X}_{r \epsilon}(u)^2]} du   \geq r^{\alpha- \beta}  \}}   M_{q \gamma}(dx)  \Big]    \\
& \leq  \liminf_{\epsilon \to 0}\E\Big[ \int_{[0,1]^2} \ind_{\{ \int_{B(x,r)}  e^{  \gamma \bar{X}_{r \epsilon}(u)-\frac{ \gamma^2}{2}\E[\bar{X}_{r \epsilon}(u)^2]} du   \geq r^{\alpha- \beta} \} }   M_{q \gamma}(dx)\Big  ] \\
& \leq D r^\eta .
\end{align*}
\qed


\subsubsection*{Upper bound}
Here, we prove the following result:

\begin{lemma}\label{lemup}
Fix $\beta>0$. There exist two constants $D>0$ and $\eta>0$ such that for all $r>0$, we have:
\begin{equation*}
\E[   \int_{[0,1]^2}  \ind_{\{M_{\gamma}(B(x,r)) \leq r^{\alpha+ \beta}\}}  M_{q \gamma}(dx) \Big ]  \leq D r^{\eta}.
\end{equation*}
\end{lemma}

\proof Along the same lines as the proof of the lower bound, we get that:
\begin{align*}
& \E\Big[  \ind_{\{ \int_{B(x,r)}  e^{  \gamma \bar{X}_{r \epsilon}(u)-\frac{ \gamma^2}{2}E[\bar{X}_{r \epsilon}(u)^2]} du   \leq r^{\alpha+ \beta}\}} e^{ q \gamma \bar{X}_{\epsilon'}(x)-\frac{q^2 \gamma^2}{2}\E[\bar{X}_{\epsilon'}(x)^2]} \Big ]    \\
& \leq   \Pb\Big[ e^{-C} e^{\gamma \Omega_{1/r}}  \leq r^{\frac{\beta}{2}}] + \Pb\Big[ \int_{B(0,1)}  \frac{e^{  \gamma \bar{X}_{\epsilon}(v)-\frac{ \gamma^2}{2}E[\bar{X}_{\epsilon}(v)^2]}}{ (v+ \epsilon)^{q \gamma^2}  } dv   \leq r^{ \frac{\beta}{2}} \Big]  \\ 
& \leq   e^{C\beta\gamma^{-2}} r^{\frac{\beta^2}{8 \gamma^2}} +  r ^{\frac{\beta}{2}} \E\Big[ \Big(\int_{B(0,1)}  \frac{e^{  \gamma \bar{X}_{\epsilon}(v)-\frac{ \gamma^2}{2}E[\bar{X}_{\epsilon}(v)^2]}}{ (|v|+ \epsilon)^{q \gamma^2}  } dv\Big)^{-1}\Big ]  \\ 
& \leq   e^{C\beta\gamma^{-2}}r^{\frac{\beta^2}{8 \gamma^2}} + r ^{\frac{\beta}{2}} \Big[ \Big(\int_{B(0,1)}   e^{  \gamma \bar{X}_{\epsilon}(v)-\frac{ \gamma^2}{2}E[\bar{X}_{\epsilon}(v)^2]} dv\Big)^{-1}\Big ]  \\ 
& \leq e^{C\beta\gamma^{-2}} r^{\frac{\beta^2}{8 \gamma^2}} +  r ^{\frac{\beta}{2}}    
\end{align*}
where, in the last inequality, we have used finiteness of the negative moments (see \cite{cf:RoVa}).\qed

\vspace{2mm}
We can now prove theorem \ref{analysemulti}. Fix $\beta>0$. By using lemmas \ref{lemlow}, \ref{lemup} and the Borel-Cantelli lemma, we get that $M_{q\gamma}(K_{(\beta)}^c)=0$ where:
\begin{equation*}
K_{(\beta)}=  \bigcup_{N \geq 1} \bigcap_{n \geq N}  \Big\{ x\in U;\frac{1}{2^{n(\alpha+\beta)}} \leq M_{\gamma}(B(x,\frac{1}{2^n})) \leq \frac{1}{2^{n(\alpha- \beta)}}  \Big\}.
\end{equation*}
We conclude by setting $$K= \bigcap_{j \geq 1}  K_{(1/j)}  \cap \Big\{ x \in U; \: \underset{\epsilon \to 0}{\lim} \: \frac{  X_\epsilon(x)}{-\ln \epsilon} = \gamma q \Big\}$$ since it is obvious that $M_{q\gamma} ( \big\{ x \in U; \: \underset{\epsilon \to 0}{\lim} \: \frac{  X_\epsilon(x)}{-\ln \epsilon} = \gamma q  \big\} ^c)=0$. By taking $q=1$, we have proved that
$$M_\gamma\Big(\Big\{ x \in U; \: \underset{\epsilon \to 0}{\lim} \: \frac{\ln M_{ \gamma}(B(x,\epsilon))}{\ln \epsilon} = d-\frac{  \gamma^2}{2} \Big\}^c\Big)=0,$$ and therefore
$$M_{q\gamma}\Big(\Big\{ x \in U; \: \underset{\epsilon \to 0}{\lim} \: \frac{\ln M_{ q\gamma}(B(x,\epsilon))}{\ln \epsilon} = d-\frac{ q^2 \gamma^2}{2} \Big\}^c\Big)=0.$$
Therefore $M_{q\gamma} (  \bar{K}^c)=0$ where we have set
\begin{equation*}
\bar{K}  = K_{(\beta)} \cap \Big\{ x \in U; \: \underset{\epsilon \to 0}{\lim} \: \frac{\ln M_{q\gamma}(B(x,\epsilon))}{\ln \epsilon} = d-\frac{q^2 \gamma^2}{2} \Big\}.
\end{equation*}\qed

\subsubsection*{Proof of Theorem \ref{cor:analysemulti}}
Finally, we adress the issue of the Hausdorff dimension of $K_{\gamma,1}$. By proposition 4.9 (a) in \cite{Falk}, we have $\text{dim}_{\text{H}}(  K_{\gamma,1} ) \geq d-\frac{\gamma^2}{2}$. 
Fix $\beta>0$. We cut $[0,1]^2$ into the standard $4^n$ dyadic intervals $I_j^n=x_j^n+ [0,\frac{1}{2^n}[^2$ of size length $\frac{1}{2^n}$. We introduce the following random variables:
\begin{equation*}
 \mu_j^n= \underset{\eta \to 0}{\lim}   \: 4^n \int_{I_j^n}e^{\gamma (X_{\eta}(x)-X_{1/2^n}(x_j^n)) -\frac{\gamma^2}{2} \ln (2^n / \eta) } dx, \; \; 
\end{equation*}
and
\begin{equation*}
 \alpha_j^n= \underset{x \in I_j^n}{\sup}  |X_{1/2^n}(x)-  X_{1/2^n}(x_j^n)|.
\end{equation*}
By using Kahane's convexity inequalities and the fact that the $M_\gamma(O)$ has negative moments when $O$ is an open set (see \cite{cf:RoVa}), we get the existence for all $q>0$ of a constant $C_q$ (independent of n) such that:
\begin{equation*}
\E[ \frac{1}{ (\mu_j^n)^q}] \leq C_q, \; \; n \geq 1, \; 1 \leq j \leq 4^n.
\end{equation*}
By using the fact that $(X_\epsilon )_{\epsilon>0}$ is a smooth Gaussian approximation, we also get the existence of some constant $C_q$ (independent of n) such that (see \cite{cf:LeTa} for example):
\begin{equation*}
 \E[   e^{q\alpha_j^n} ] \leq C_q, \; \; n \geq 1, \; 1 \leq j \leq 4^n.
\end{equation*}
  Therefore, one easily gets that:
\begin{equation*}
\Pb( \cup_{N \geq 1} \cap_{n \geq N}  \lbrace   \underset{1 \leq j \leq 4^n}{\inf} \mu_j^n e^{-3 \gamma \alpha_j^n }    \geq \frac{1}{2^{n \beta}} \rbrace   )=1.
\end{equation*}

For all $x$ in $I_j^n$, we get that $M_{\gamma}( I_j^n  ) \geq  \mu_j^n e^{-3 \gamma \alpha_j^n } 4^n e^{\gamma X_{1/2^n}(x)-\frac{\gamma^2}{2} \ln 2^n }$; hence, we deduce that on $K_{\gamma,1}$ we have:
\begin{equation*}
\underset{r \to 0}{\overline{\lim}} \frac{M_{\gamma}(  B(x,r) )}{r^{d-\gamma^2/2+\beta}} = \infty 
\end{equation*}
Since this is valid for all $\beta>0$, proposition 4.9 (b) in \cite{Falk} gives $\text{dim}_{\text{H}}(  K_{\gamma,1} ) \leq d-\frac{\gamma^2}{2}$. \qed


\bigskip


\begin{thebibliography}{20}


\bibitem{AidShi}
A\"{\i}d\'ekon E., Shi Z.: The Seneta-Heyde scaling for the branching random walk, arXiv:1102.0217v2.

\bibitem{Albeverio1}
Albeverio, S., Hoegh-Krohn, R.: The Wightman axioms and the mass gap for strong interactions
of exponential type in two dimensional space-time. J. Funct. Anal. 16, 39-82 (1974)


\bibitem{Albeverio2}
S. Albeverio, G. Gallavotti and R. Hoegh-Krohn: 
Some Results for the Exponential Interaction in Two or More Dimensions Commun.Math.Phys.70,187-192 (1979).

\bibitem{aidekon}
A\"{\i}d\'ekon E.: Convergence in law of the minimum of a branching random walk, 	 arXiv:1101.1810v3, to appear in Annals of Probability.

 


\bibitem{Rnew1}
Allez R., Rhodes R., Vargas V.: Lognormal $\star$-scale invariant random measures, to appear in {\it Probability Theory and Related Fields}, arXiv:1102.1895v1.

\bibitem{amb}
Ambj\o rn J., Boulatov D., Nielsen J.L., Rolf J., Watabiki Y.: The spectral dimension of $2D$ Quantum Gravity. JHEP 9802 (1998) 010, arXiv:hep-lat/9808027v1.

\bibitem{amb2}
Ambj\o rn J., Anagnostopoulos K.N., Ichihara T., Jensen L., Watabiki Y.: Quantum Geometry and Diffusions. JHEP11(1998)022.

\bibitem{arguin}
Arguin L-P., Zindy O.: Poisson-Dirichlet statistics for the extremes of a log-correlated Gaussian field, arXiv:1203.4216v1.
\bibitem{alvarez}
Alvarez-Gaum\'e L., Barbon J.L., Crnkovic C., Nucl. Phys., B394 (1993), 383.

\bibitem{cf:BaDeMu} Bacry, E., Delour, J., and Muzy, J.F.: Multifractal random walks, 
 \emph{Phys. Rev. E}, \textbf{64} (2001), 026103-026106.
 
\bibitem{cf:BaKoMu} Bacry E., Kozhemyak, A., Muzy J.-F.: Continuous
  cascade models for asset returns, available at
  www.cmap.polytechnique.fr/~bacry/biblio.html, to appear in \emph{Journal of Economic Dynamics and Control}. 
  
\bibitem{cf:BaMu} Bacry E., Muzy J.F.: Log-infinitely divisible multifractal process, 
 \emph{Communications in Mathematical Physics}, \textbf{236} (2003), 449-475.

\bibitem{ismael}
Bailleul I.: KPZ in a multidimensional random geometry of multiplicative cascades, available at "http://perso.univ-rennes1.fr/ismael.bailleul/".


\bibitem{Barral1}
Barral J.: Moments, continuit\'e, et analyse multifractale des martingales de Mandelbrot, \emph{Probab. Theory Relat. Fields}, \textbf{113} (1999), 535--569.

\bibitem{barralJTP}
Barral J., Continuity of the multifractal spectrum of a random statistically self-similar measure, J. Theoretic. Probab., 13 (2000), 1027--1060.

\bibitem{BFP}
Barral J.,  Fan A.-H., Peyrire J. , Mesures engendr\'ees par multiplications, Panoramas et Synth\`eses 32: Quelques interactions entre analyse, probabilit\'es et fractals, 2010.


\bibitem{Rnew4}
Barral J., Jin X., Rhodes R., Vargas V.: Gaussian multiplicative chaos and KPZ duality, to appear in Commun.Math.Phys., arXiv:1202.5296v2.

\bibitem{BXJ}
Barral J., Jin X.: On exact scaling log-infinitely divisible cascades, arXiv:1208.2221 [math.PR].

\bibitem{BKNSW} Barral J., Kupiainen A., Nikula M., Saksman E., Webb C.: Critical Mandelbrot cascades, arXiv:1206.5444v1.

\bibitem{basic} 
Barral J., Kupiainen A., Nikula M., Saksman E., Webb C.: Basic properties of critical lognormal multiplicative chaos, arXiv:1303.4548.

\bibitem{Bar}
Barral, J., Mandelbrot, B.B.: Multifractal products of cylindrical pulses, \emph{Probab. Theory
Relat. Fields} \textbf{124} (2002), 409-430.

\bibitem{BMI}
Barral, J., Mandelbrot, B.B.: Random Multiplicative Multifractal Measures I, II, III, in M. Lapidus, M. van Frankenhuijsen, eds., Fractal geometry and applications: a jubilee of Beno\^{\i}t Mandelbrot, Proceedings of Symposia in Pure Mathematics, vol 72., n. 2, 3-90, 2004.

\bibitem{Rnew3}
Barral J., Rhodes R., Vargas V.: Limiting laws of supercritical branching random walks,	Comptes rendus - Mathematique, 350,   pp. 535-538,  2012.

\bibitem{BS} Barral, J., Seuret, S., The singularity spectrum of
  L\'evy processes in multifractal time. {\it Adv.\ Math}.\ {\bf 214}
  (2007), 437--468.

\bibitem{bauer}
Bauer M., Bernard D., Cantini L.:  Off-Critical SLE(2) and SLE(4): a Field Theory Approach. J. Stat. Mech. (2009) P07037.


\bibitem{BenjCur}
Benjamini, I., Curien, N.: Simple random walk on the uniform infinite planar quadrangulation: Subdiffusivity via pioneer points, arXiv:1202.5454.

\bibitem{Benj}
Benjamini, I., Schramm, O.: KPZ in one dimensional random geometry of multiplicative cascades, Communications in Mathematical Physics, vol. 289, no 2, 653-662, 2009.

\bibitem{BiKi}
Biggins, J.D. and Kyprianou, A.E.: Measure change in multitype branching, \emph{Adv. Appl. Probab.} \textbf{36} (2004), 544Ð581.

\bibitem{Biggf}
Biggins, J.D. and Kyprianou, A.E.: Fixed points of the smoothing transform; the boundary case, \emph{Electronic Journal of Probability}, vol. 10 (2005).

\bibitem{boul}
D.V. Boulatov and V.A. Kazakov, Phys. Lett. 184B (1987) 247.

\bibitem{bramson}
Bramson M., Ding J., Zeitouni O.: Convergence in law of the  maximum of the two-dimensional discrete Gaussian Free Field, arXiv:1301.6669v1.

\bibitem{BKZ}
Br\'ezin E., Kazakov V.A. ,  Zamolodchikov Al.B.: Scaling violation in a field theory of closed strings in one physical dimension, \emph{Nuclear Physics} \textbf{B338},  673-688 (1990).

 \bibitem{CarDou}
Carpentier D., Le Doussal P.: Glass transition of a particle in a random potential, front selection in nonlinear RG and entropic phenomena in Liouville and Sinh-Gordon models, \emph{Phys. Rev. E} \textbf{63}:026110 (2001).

\bibitem{cf:Castaing}
Castaing B., Gagne Y., Hopfinger E.J.:
Velocity probability density-functions of high Reynolds-number turbulence,
\emph{Physica D} \textbf{46} (1990) 2, 177-200.

\bibitem{cf:Cas}
Castaing B., Gagne Y., Marchand M.:
Conditional velocity pdf in 3-D turbulence,
\emph{J. Phys. II France} \textbf{4} (1994), 1-8.


\bibitem{cf:Cha} Chainais, P.: Multidimensional infinitely divisible cascades. Application to the modelling of intermittency in turbulence, \emph{European Physical Journal B}, \textbf{51} no. 2 (2006), pp. 229-243.

\bibitem{chelkak}
Chelkak D., Smirnov S.: Discrete complex analysis on isoradial graphs, Advances in Mathematics 228 (2011) 1590Ð1630. 

\bibitem{chen}
Chen L., Jakobson D.: Gaussian Free Fields and KPZ Relation in $\R^4$,  arXiv:1210.8051.

\bibitem{Rnew5} 
 Chevillard L., Rhodes R., Vargas V.: Gaussian multiplicative chaos for symmetric isotropic matrices, arxiv.

\bibitem{CheRob10}
Chevillard L., Robert R., Vargas V. (2010), A Stochastic Representation of the Local Structure of Turbulence, \emph{EPL}, \textbf{89}, 54002. 

\bibitem{cf:Cizeau} Cizeau, P., Gopikrishnan, P., Liu, Y., Meyer, M., Peng, C.K., Stanley, E.: Statistical properties of the volatility of price fluctuations, \emph{Physical Review E}, \textbf{60} no.2 (1999), 1390-1400.


\bibitem{cf:Co} Cont, R.: Empirical properties of asset returns: stylized facts and statistical issues, 
 \emph{Quantitative Finance}, \textbf{1} no.2 (2001), 223-236.

\bibitem{cf:DaVe} Daley D.J., Vere-Jones D.: An introduction to the
  theory of point processes, Springer-Verlag, (1988).

\bibitem{Das}
Das S.R., Dhar, A., Sengupta A.M., Wadia S.R., Mod. Phys. Lett., A5 (1990), 1041.

\bibitem{disc:da1}
David F., Nucl. Phys. B257 (1985) 45.

\bibitem{disc:da2}
David F., Nucl. Phys. B257 (1985) 543.

\bibitem{cf:Da} David F.: Conformal Field Theories Coupled to 2-D Gravity in the Conformal Gauge, \emph{Mod. Phys. Lett. A} \textbf{3} 1651-1656 (1988).

\bibitem{DistKa} Distler J.,   Kawai H.: Conformal Field Theory and 2-D Quantum Gravity or Who's Afraid of Joseph Liouville?, \emph{Nucl. Phys.} \textbf{B321} 509-517 (1989).

\bibitem{cf:DuRoVa} Duchon, J., Robert, R., Vargas, V.: Forecasting volatility with the multifractal random walk model, submitted to \emph{Mathematical Finance}, available at http://arxiv.org/abs/0801.4220.    

\bibitem{Dup:houches}
Duplantier B.: A rigorous perspective on Liouville quantum gravity and KPZ, in \emph{Exact Methods in Low-dimensional Statistical Physics and Quantum Computing},
J. Jacobsen, S. Ouvry, V. Pasquier, D. Serban, and L.F. Cugliandolo, eds., Lecture Notes of the Les Houches Summer School: Volume 89, July 2008, Oxford University Press (Clarendon, Oxford) (2010).

\bibitem{dupRW}
Duplantier B.: Random walks and quantum gravity in two dimensions. Phys. Rev. Lett., 81 (25):5489-5492, 1998.

\bibitem{dupdual}
Duplantier B., in Fractal geometry and applications: a jubilee of Beno"t Mandelbrot, Part 2 (Amer. Math. Soc., Providence, RI, 2004), vol 72 of Proc. Sympos. Pure Math. 365-482, arXiv:math-ph/0303034.

\bibitem{dupkos}
Duplantier B., Kostov I.: Conformal spectra on a random surface, Phys. Rev. Lett., vol 61 no 13 (1988), 1433. 

\bibitem{DupKwo}
Duplantier D. and Kwon K.-H.: Conformal invariance and intersection of random
walks. Phys. Rev. Lett., 61:2514Ð2517, 1988.

\bibitem{PRL} Duplantier, B., Sheffield, S.: Duality and KPZ in Liouville Quantum Gravity, \emph{Physical Review Letters} \textbf{102}, 150603 (2009).

\bibitem{cf:DuSh} Duplantier, B., Sheffield, S.: Liouville Quantum Gravity and KPZ, \emph{Inventiones Mathematicae} \textbf{185} (2) (2011) 333-393.

\bibitem{Rnew7}
Duplantier  B., Rhodes R., Sheffield S., Vargas V.: Critical Gaussian multiplicative chaos: convergence of the derivative martingale, arXiv.

\bibitem{Rnew12}
Duplantier  B., Rhodes R., Sheffield S., Vargas V.: Renormalization of Critical Gaussian Multiplicative Chaos and KPZ, arXiv.

\bibitem{DSRV3}
Duplantier  B., Rhodes R., Sheffield S., Vargas V.: Log-correlated Gaussian fields: a primer, in preparation.

\bibitem{durrett} Durrett R., Liggett T.M.: Fixed points of the smoothing transformation, {\it Probability Theory and Related Fields}  \textbf{64} (3) (1983) 275-301.

\bibitem{Falk}
Falconer K.J.: The geometry of fractal sets, Cambridge University Press, 1985.

\bibitem{reviewfan}
Fan A.H.: Some topics in the theory of multiplicative chaos, in Progress in Probability, 57, Ch. Bandt, M. Zaehle and U. Mosco (Ed.), Birkhauser, 2004, 119-134.

\bibitem{Fan}
Fan A.H.: Sur le chaos de L\'evy d'indice $0<\alpha<1$, \emph{Ann. Sciences Math. Qu\'ebec}, vol 21 no. 1, 1997, p. 53-66.

\bibitem{FoiMan01}
Foias C., Manley O., Rosa R., Temam R. (2001) Navier-Stokes equations and turbulence, Cambridge University Press, Cambridge.

\bibitem{cf:Fr} Frisch, U.: \emph{Turbulence}, Cambridge University Press (1995).                                                                                                                                                                                                                                                                                                                                                                                                                                                                                                                                                                                                                                                                                                                                                              

\bibitem{disc:Froh}
Frohlich J., in Lecture Notes in Physics, vol. 216, ed.L. Garrido (Springer, 1985); Ambjorn J., Durhuus B. and Frohlich J., Nucl. Phys. B257 (1985) 433.

\bibitem{Fyo}
Fyodorov Y. and Bouchaud J.P.: Freezing and extreme-value statistics in a random energy model with logarithmically correlated potential, \emph{J.  Phys. A} \textbf{41} (2008) 372001.


\bibitem{rosso}
Fyodorov Y, Le Doussal P., Rosso A.:  Statistical Mechanics of Logarithmic REM: Duality, Freezing and Extreme Value Statistics of $1/f$ Noises generated by Gaussian Free Fields, 	 \emph{J. Stat. Mech.} (2009) P10005.

\bibitem{FLR}
Fyodorov Y, Le Doussal P., Rosso A.:  Freezing transition in decaying Burgers turbulence and random matrix dualities, \emph{Europhysics Letters}, 90,  (2010) 60004.

\bibitem{garban}
Garban C.: Quantum gravity and the KPZ formula, s\'eminaire Bourbaki, 64e annŽe, 2011-2012, no 1052. 

\bibitem{LBM}
Garban C, Rhodes R., Vargas V.: Liouville Brownian motion, arXiv:1301.2876.

\bibitem{LBM2}
Garban C, Rhodes R., Vargas V.: On the heat kernel and the Dirichlet form of Liouville Brownian motion, arXiv:1302.6050.

\bibitem{GM} Ginsparg P. and Moore G.: Lectures on 2D gravity and 2D string theory, in \textit{Recent direction in particle
theory}, Proceedings of the 1992 TASI, edited by J. Harvey
and J. Polchinski (World Scientific, Singapore, 1993).

\bibitem{GZ} Ginsparg P., Zinn-Justin J.: 2D gravity $+$ 1D matter, \emph{Physics Letters} \textbf{B240},  333-340 (1990).

\bibitem{glimm}
Glimm J., Jaffe A.: Quantum Physics: a functional integral point of view, Berlin-Heidelberg-New York, Springer-Verlag 1981.

\bibitem{cf:Gn} Gneiting, T.: Criteria of Polya type for radial positive definite functions, \emph{Proceedings of the American Mathematical Society}, \textbf{129} no. 8 (2001), 2309-2318.

\bibitem{golubov}
Gneiting, T.: Criteria of Polya type for radial positive definite functions, Proceedings of the
American Mathematical Society, 129 no. 8 (2001), 2309-2318.

\bibitem{cf:GrMacMa} Gray, Mathews, Macrobert: \emph{Bessel Functions}, Macmillan and co. (1922). 

 \bibitem{GrossKleban}
Gross, D.J. and Klebanov I.R.: One-dimensional string theory on a circle, \emph{Nuclear Physics} \textbf{B344} (1990) 475--498.


\bibitem{GuTo}
Guerra F., Toninelli F.L.: The Thermodynamic Limit in Mean Field Spin Glass Models,
\emph{Commun. Math. Phys.} \textbf{230}, 71-79 (2002).

\bibitem{Hoeg}
Hoegh-Krohn, R.: A general class of quantum fields without cut offs in two space-time dimensions.
Commun. Math. Phys. 21, 244-255 (1971). 


\bibitem{peres}
Hu X., Miller J., Peres Y.: Thick points of the Gaussian free field, Annals of Probability, 2010, vol. 38, 896-926.

\bibitem{HuShi}
Hu Y., Shi Z.: Minimal position and critical martingale convergence in branching random walks, and directed polymers on disordered trees, Annals of Probability, vol. 37, no 2, 2009, 742-789.


\bibitem{janson}
Janson S.: Gaussian Hilbert spaces, vol. 129, Cambridge Tracts in Mathematics. Cambridge University Press, Cambridge, 1997. 

\bibitem{cf:Kah} Kahane, J.-P.: Sur le chaos multiplicatif,
  \emph{Ann. Sci. Math. Qu{\'e}bec}, \textbf{9} no.2 (1985), 105-150.

\bibitem{kahane74} Kahane, J.-P.: Sur le mod\`ele de
  turbulence de Beno\^{\i}t Mandelbrot., \emph{C.R.\ Acad.\ Sci.\ Paris},
  \textbf{278} (1974), 567--569.
  
\bibitem{KP} Kahane, J.-P., Peyri\`ere, J., Sur certaines martingales
  de B.~Mandelbrot. {\it Adv.\ Math}.\ {\bf 22} (1976), 131--145. 

\bibitem{KKK} Kazakov V.,  Kostov I., and  Kutasov D.:  A Matrix Model for the 2d Black Hole,
in \textit{Nonperturbative Quantum Effects 2000}, JHEP Proceedings.

\bibitem{kaz1}
Kazakov V.A., Phys. Lett. 119A (1986) 140.

\bibitem{kaz2}
Kazakov V.A., Nucl. Phys. B (Proc. Supp.) 4 (1988) 93.
 

\bibitem{Kleb1}
Klebanov I.R.: Touching random surfaces and Liouville gravity, \emph{Phys. Rev. D 51}, 1836Ð1841 (1995).

\bibitem{Kleb2}
Klebanov I.R., Hashimoto A.: Non-perturbative Solution of Matrix Models Modified by Trace-squared Terms, \emph{Nucl. Phys}. \textbf{B434} (1995) 264-282.

\bibitem{Kleb3}
Klebanov I.R., Hashimoto A.: Wormholes, Matrix Models, and Liouville Gravity, \emph{Nucl. Phys.} (Proc. Suppl.) \textbf{45}B,C (1996) 135-148.

\bibitem{cf:KPZ} Knizhnik V.G., Polyakov A.M., Zamolodchikov A.B.: Fractal structure of 2D-quantum gravity, \emph{Modern Phys. Lett A} \textbf{3}(8) (1988), 819-826.  

\bibitem{Kol41}
Kolmogorov A. N. (1941), The local structure of turbulence in incompressible
viscous fluid for very large Reynolds numbers, \emph{Dokl. Akad. Nauk SSSR} \textbf{30}
  301 [in Russian]. English translation: \emph{Proc. R. Soc. London},
Ser. A 434 (1991), 9.


\bibitem{cf:Kol} Kolmogorov A.N.: A refinement of previous hypotheses
  concerning the local structure of turbulence,
  \emph{J. Fluid. Mech.}, \textbf{13} (1962), 83-85. 
  
\bibitem{korchemsky}
Korchemsky, Mod. Phys. Lett., A7 (1992) 3081; Phys. Lett. 296B (1992), 323.

\bibitem{Kyp}
Kyprianou, A.E.: Slow variation and uniqueness of solutions to the
functional equation in the branching random walk. \emph{J. Appl. Probab.} \textbf{35} (1998) 795Ð802.

\bibitem{lalley}
Lalley S.P., Sellke T.: A conditional limit theorem for frontier of a branching Brownian motion, \emph{Ann. Probab.} 15, 1987, 1052-1061.

\bibitem{Lawler}
Lawler G., Limic, V.: Random Walk, a modern introduction.



\bibitem{LSW1}
Lawler G.F., Schramm O., and Werner W.: Values of {B}rownian intersection exponents. {I}. {H}alf-plane
  exponents.
\emph {Acta Math.}, 187(2):237--273, 2001, arXiv:math.PR/9911084.

\bibitem{LSW2}
Lawler G.F., Schramm O., and Werner W.: Values of {B}rownian intersection exponents. {II}. {P}lane exponents.
\emph{Acta Math.}, 187(2):275--308, 2001, arXiv:math.PR/0003156.

\bibitem{LSW3}
Lawler G.F., Schramm O., and Werner W.: Values of {B}rownian intersection exponents. {III}. {T}wo-sided
  exponents.
\emph{\em Ann. Inst. H. Poincar\'e Probab. Statist.}, 38(1):109--123,
 2002, arXiv:math.PR/0005294.


\bibitem{cf:LeTa} Ledoux, M., Talagrand, M.: \emph{Probability in Banach Spaces}, Springer-Verlag (1991). 

\bibitem{LeGallS}
Le Gall, J-F.:
The topological structure of scaling limits of large planar maps. 
Invent. Math. 169 (2007), no. 3, 621Ð670. 

\bibitem{legall}
Le Gall J.F.: Uniqueness and universality of the Brownian map, arXiv:1105.4842, 2011.

\bibitem{LGM}
Le Gall J.F., Miermont G.: Scaling limits of random trees and planar maps. arXiv:1101.4856, 2011.

\bibitem{Liu}
Liu Q.S., Fixed points of a generalized smoothing transformation and applications to
the branching random walk, \emph{Adv. Appl. Probab.} {\bf 30} (1998) 85-112.

 \bibitem{madaule}
Madaule T.: Convergence in law for the branching random walk seen from its tip, arXiv:1107.2543v2.

\bibitem{makarov}
Makarov N., Smirnov S.:Off-critical lattice models and massive SLEs, in Exner, Pavel (ed.), XVIth International Congress on Mathematical Physics, Prague, 3 - 8 August 2009. 362-371, World Sci. Publ., Singapore. arXiv:0909.5377 

\bibitem{cf:Man} Mandelbrot B.B.: A possible refinement of the
  lognormal hypothesis concerning the distribution of energy in
  intermittent turbulence, \emph{Statistical Models and Turbulence},
  La Jolla, CA, Lecture Notes in Phys. no. 12, Springer, (1972), 333-351.


\bibitem{mandelbrot}
Mandelbrot B.B.: Intermittent turbulence in self-similar cascades, divergence of high moments and dimension of the carrier, \emph{J. Fluid. Mech.} \textbf{62} (1974), 331-358.

\bibitem{maru} Maruyama, G: Infinitely divisible processes,
  \emph{Theory Probab. Appl.}, \textbf{9} no.15 (1970), 3-23.
  
%
\bibitem{miermont}  
Miermont G.: The Brownian map is the scaling limit of uniform random plane quadrangulations, arXiv:1104.1606, 2011.  
  
\bibitem{Molchan}
 Molchan, G.M.: Scaling exponents and multifractal dimensions for independent random cascades, \emph{Comm. Math Physics} \textbf{179} (1996), 681-702.  
   
 \bibitem{neveu}
Neveu, J.: Multiplicative martingales for spatial branching processes. In Seminar on Stochastic Processes, 1987, (eds: E. Cinlar, K.L. Chung, R.K. Getoor). Progress in Probability and Statistics, 15, 223-241. Birkhauser, Boston (1988).

\bibitem{Nienhuis} Nienhuis, B.: Coulomb gas formulation of two-dimensional phase transitions, in \emph{Phase Transitions and Critical Phenomena},
edited by C. Domb and J. L. Lebowitz, (Academic, London, 1987), Vol. 11.
 
\bibitem{cf:Obu} Obukhov A.M.: Some specific features of atmospheric
  turbulence, \emph{J. Fluid. Mech.}, \textbf{13} (1962), 77-81. 

\bibitem{Parisi}
Parisi G.: On the one dimensional discretized string, \emph{Physics Letters} \textbf{B238},  209-212
(1990).

\bibitem{ParisiFrisch}
Parisi G., Frisch U.:  On the singularity structure
of fully developed turbulence, \emph{Proceed. Turbulence and predictability
in geophysical fluid dynamics and climate
dynamics}, (M. Ghil, R. Benzi and G. Parisi, eds, 1985)
Ed. Varenna, 1983, pp. 84-87. 

\bibitem{cf:PaYu} Pasenchenko, O. Yu.: Sufficient conditions for the characteristic function of a two-dimensional isotropic distribution, \emph{Theory Probab. Math. Statist.}, \textbf{53} (1996), 149-152.

\bibitem{pey74} Peyri\`ere, J., Turbulence et dimension de Hausdorff.
  \textit{C.~R.\ Acad.\ Sc.\ Paris} 278 (1974), 567--569.
  

\bibitem{Polch} Polchinski J., Critical behavior of random surfaces in one dimension, \emph{Nuclear Physics} {\bf B346} (1990) 253--263.

\bibitem{polch}
Polchinski J.: String theory, Volume 1, An introduction to the Bosonic String. Cambridge Monographs on Mathematical Physics. Cambridge university press 1998.

\bibitem{Pol}
Polyakov A.M., Phys. Lett. 163B (1981) 207.

\bibitem{Rnew8} Rhodes R., Sohier J., Vargas, V.: $\star$-scale invariant random measures, to appear in Annals of Probability, available on arxiv.

\bibitem{Rnew9} Rhodes R., Vargas, V.: Multidimensional multifractal random measures, \emph{Electronic Journal of Probability}, \textbf{15} (2010) 241-258.

\bibitem{Rnew10} Rhodes, R. Vargas, V.: KPZ formula for log-infinitely divisible multifractal random measures,  \emph{ESAIM Probability and Statistics}, \textbf{15} (2011) 358-371.

\bibitem{Rnew11} Rhodes, R. Vargas, V.: Optimal transportation for multifractal random measures and applications, \emph{Ann. Inst. H. PoincarŽ Probab. Statist.} \textbf{49}, Number 1 (2013), 119-137. 

\bibitem{spectral} Rhodes, R. Vargas, V.: Spectral dimension of Liouville quantum gravity, arxiv.

 \bibitem{cf:RoVa1} Robert, R. Vargas, V.: Hydrodynamic Turbulence and Intermittent Random Fields, \emph{Communications in Mathematical Physics}, \textbf{284} (3) (2008), 649-673.


\bibitem{cf:RoVa} Robert, R., Vargas, V.: Gaussian multiplicative chaos revisited, \emph{Annals of Probability}, \textbf{38} 2 (2010) 605-631.

\bibitem{cf:Sch} Schmitt, F., Lavallee, D., Schertzer, D., Lovejoy, S.: Empirical determination of universal
multifractal exponents in turbulent velocity fields, \emph{Phys. Rev. Lett.} \textbf{68} (1992), 305-308.

\bibitem{She07}
Sheffield S.: Gaussian free fields for mathematicians, \emph{Probab. Th. Rel. Fields}, \textbf{139} (2007) 
521--541.

\bibitem{QZ}
Sheffield, S.: Conformal weldings of random surfaces: SLE and the quantum gravity zipper. preprint, arXiv:1012.4797.

\bibitem{Smir}
Smirnov, S.: Conformal invariance in random cluster models. I. Holomorphic fermions in the Ising model, \emph{Ann. Math.}, 172 (2010), 1435-1467. 

\bibitem{cf:Sto} Stolovitzky, G., Kailasnath, P., Sreenivasan, K.R.: Kolmogorov's Refined Similarity Hypotheses, \emph{Phys. Rev. Lett.} \textbf{69}(8) (1992), 1178-1181.  

\bibitem{Webb} C. Webb, Exact asymptotics of the freezing transition of a logarithmically correlated random energy model, J. Stat. Phys., {\bf 145} (2011), 1595--1619.


\end{thebibliography}
\end{document}